\pdfoutput=1

\documentclass[12pt,a4paper]{article}

\usepackage{amsfonts,amssymb,amsthm,amsmath,latexsym}
\usepackage{enumerate}
\usepackage{mathtools}
\usepackage{subeqnarray,eqsection,indent,eepicemu,url,cite,bm}
\usepackage{multirow}  
\usepackage{tensor}  
\usepackage{graphicx,color}
\usepackage{float}  

\usepackage{tikz,tkz-graph}
\usetikzlibrary{shapes}

\date{February 8, 2023}   

\begin{document}

\title{\vspace*{-2cm}Total positivity of some polynomial matrices \\[0.5mm]
       that enumerate labeled trees and forests \\[5mm]
       \hspace*{-8mm}
       II.~Rooted labeled trees and partial functional digraphs
       \hspace*{-4mm}
      }

\author{ \\
     \hspace*{-1cm}
      {\large Xi Chen${}^{1}$ and Alan D.~Sokal${}^{2,3}$}
   \\[5mm]
     \hspace*{-1.05cm}
      \normalsize
           ${}^1$School of Mathematical Sciences, Dalian University of
                   Technology, Dalian 116024, CHINA  \\[1mm]
     \hspace*{-2.5cm}
      \normalsize
           ${}^2$Department of Mathematics, University College London,
                    London WC1E 6BT, UK   \\[1mm]
     \hspace*{-4.1cm}
      \normalsize
           ${}^3$Department of Physics, New York University,
                    New York, NY 10003, USA
       \\[5mm]
     \hspace*{-7mm}
     {\tt chenxi@dlut.edu.cn},\, {\tt sokal@nyu.edu}  \\[1cm]
}

\maketitle
\thispagestyle{empty}   

\begin{abstract}
We study three combinatorial models for the lower-triangular matrix
with entries $t_{n,k} = \binom{n}{k} n^{n-k}$:
two involving rooted trees on the vertex set $[n+1]$,
and one involving partial functional digraphs on the vertex set $[n]$.
We show that this matrix is totally positive
and that the sequence of its row-generating polynomials
is coefficientwise Hankel-totally positive.
We then generalize to polynomials $t_{n,k}(y,z)$
that count improper and proper edges,
and further to polynomials $t_{n,k}(y,\bm{\phi})$
in infinitely many indeterminates
that give a weight $y$ to each improper edge
and a weight $m! \, \phi_m$ for each vertex with $m$~proper children.
We show that if the weight sequence~$\bm{\phi}$ is Toeplitz-totally positive,
then the two foregoing total-positivity results continue to hold.
Our proofs use production matrices and exponential Riordan arrays.
\end{abstract}


\bigskip
\bigskip
\noindent
{\bf Mathematics Subject Classification (MSC 2010) codes:}
05A15 (Primary);
05A19, 05A20, 05C05, 05C30, 15B05, 15B36, 15B48, 30E05, 44A60 (Secondary).

\clearpage

\newtheorem{theorem}{Theorem}[section]
\newtheorem{proposition}[theorem]{Proposition}
\newtheorem{lemma}[theorem]{Lemma}
\newtheorem{corollary}[theorem]{Corollary}
\newtheorem{definition}[theorem]{Definition}
\newtheorem{conjecture}[theorem]{Conjecture}
\newtheorem{question}[theorem]{Question}
\newtheorem{problem}[theorem]{Problem}
\newtheorem{openproblem}[theorem]{Open Problem}
\newtheorem{example}[theorem]{Example}
\newtheorem{remark}[theorem]{Remark}

\renewcommand{\theenumi}{\alph{enumi}}
\renewcommand{\labelenumi}{(\theenumi)}
\def\eop{\hbox{\kern1pt\vrule height6pt width4pt
depth1pt\kern1pt}\medskip}
\def\prf{\par\noindent{\bf Proof.\enspace}\rm}
\def\rmk{\par\medskip\noindent{\bf Remark\enspace}\rm}

\newcommand{\textbfit}[1]{\textbf{\textit{#1}}}

\newcommand{\bigdash}{%
\smallskip\begin{center} \rule{5cm}{0.1mm} \end{center}\smallskip}

\newcommand{\safepar}{ {\protect\hfill\protect\break\hspace*{5mm}} }

\newcommand{\be}{\begin{equation}}
\newcommand{\ee}{\end{equation}}
\newcommand{\<}{\langle}
\renewcommand{\>}{\rangle}
\newcommand{\widebar}{\overline}
\def\reff#1{(\protect\ref{#1})}
\def\spose#1{\hbox to 0pt{#1\hss}}
\def\ltapprox{\mathrel{\spose{\lower 3pt\hbox{$\mathchar"218$}}
    \raise 2.0pt\hbox{$\mathchar"13C$}}}
\def\gtapprox{\mathrel{\spose{\lower 3pt\hbox{$\mathchar"218$}}
    \raise 2.0pt\hbox{$\mathchar"13E$}}}
\def\textprime{${}^\prime$}
\def\proof{\par\medskip\noindent{\sc Proof.\ }}
\def\firstproof{\par\medskip\noindent{\sc First Proof.\ }}
\def\secondproof{\par\medskip\noindent{\sc Second Proof.\ }}
\def\alternateproof{\par\medskip\noindent{\sc Alternate Proof.\ }}
\def\algebraicproof{\par\medskip\noindent{\sc Algebraic Proof.\ }}
\def\combinatorialproof{\par\medskip\noindent{\sc Combinatorial Proof.\ }}
\def\proofof#1{\bigskip\noindent{\sc Proof of #1.\ }}
\def\firstproofof#1{\bigskip\noindent{\sc First Proof of #1.\ }}
\def\secondproofof#1{\bigskip\noindent{\sc Second Proof of #1.\ }}
\def\thirdproofof#1{\bigskip\noindent{\sc Third Proof of #1.\ }}
\def\fourthproofof#1{\bigskip\noindent{\sc Fourth Proof of #1.\ }}
\def\fifthproofof#1{\bigskip\noindent{\sc Fifth Proof of #1.\ }}
\def\algebraicproofof#1{\bigskip\noindent{\sc Algebraic Proof of #1.\ }}
\def\combinatorialproofof#1{\bigskip\noindent{\sc Combinatorial Proof of #1.\ }}
\def\sketchofproof{\par\medskip\noindent{\sc Sketch of proof.\ }}
\renewcommand{\qed}{ $\square$ \bigskip}
\newcommand{\myendremark}{ $\blacksquare$ \bigskip}
\def\half{ {1 \over 2} }
\def\third{ {1 \over 3} }
\def\twothird{ {2 \over 3} }
\def\smfrac#1#2{{\textstyle{#1\over #2}}}
\def\smhalf{ {\smfrac{1}{2}} }
\newcommand{\real}{\mathop{\rm Re}\nolimits}
\renewcommand{\Re}{\mathop{\rm Re}\nolimits}
\newcommand{\imag}{\mathop{\rm Im}\nolimits}
\renewcommand{\Im}{\mathop{\rm Im}\nolimits}
\newcommand{\sgn}{\mathop{\rm sgn}\nolimits}
\newcommand{\tr}{\mathop{\rm tr}\nolimits}
\newcommand{\supp}{\mathop{\rm supp}\nolimits}
\newcommand{\disc}{\mathop{\rm disc}\nolimits}
\newcommand{\diag}{\mathop{\rm diag}\nolimits}
\newcommand{\tridiag}{\mathop{\rm tridiag}\nolimits}
\newcommand{\AZ}{\mathop{\rm AZ}\nolimits}
\newcommand{\EAZ}{\mathop{\rm EAZ}\nolimits}
\newcommand{\perm}{\mathop{\rm perm}\nolimits}
\def\hboxscript#1{ {\hbox{\scriptsize\em #1}} }
\renewcommand{\emptyset}{\varnothing}
\newcommand{\eqdef}{\stackrel{\rm def}{=}}

\newcommand{\restrict}{\upharpoonright}

\newcommand{\compinv}{{\langle -1 \rangle}}   

\newcommand{\scra}{{\mathcal{A}}}
\newcommand{\scrb}{{\mathcal{B}}}
\newcommand{\scrc}{{\mathcal{C}}}
\newcommand{\scrd}{{\mathcal{D}}}
\newcommand{\scre}{{\mathcal{E}}}
\newcommand{\scrf}{{\mathcal{F}}}
\newcommand{\scrg}{{\mathcal{G}}}
\newcommand{\scrh}{{\mathcal{H}}}
\newcommand{\scri}{{\mathcal{I}}}
\newcommand{\scrj}{{\mathcal{J}}}
\newcommand{\scrk}{{\mathcal{K}}}
\newcommand{\scrl}{{\mathcal{L}}}
\newcommand{\scrm}{{\mathcal{M}}}
\newcommand{\scrn}{{\mathcal{N}}}
\newcommand{\scro}{{\mathcal{O}}}
\newcommand\scroo{
  \mathchoice
    {{\scriptstyle\mathcal{O}}}
    {{\scriptstyle\mathcal{O}}}
    {{\scriptscriptstyle\mathcal{O}}}
    {\scalebox{0.6}{$\scriptscriptstyle\mathcal{O}$}}
  }
\newcommand{\scrp}{{\mathcal{P}}}
\newcommand{\scrq}{{\mathcal{Q}}}
\newcommand{\scrr}{{\mathcal{R}}}
\newcommand{\scrs}{{\mathcal{S}}}
\newcommand{\scrt}{{\mathcal{T}}}
\newcommand{\scrv}{{\mathcal{V}}}
\newcommand{\scrw}{{\mathcal{W}}}
\newcommand{\scrz}{{\mathcal{Z}}}

\newcommand{\bfa}{{\mathbf{a}}}
\newcommand{\bfb}{{\mathbf{b}}}
\newcommand{\bfc}{{\mathbf{c}}}
\newcommand{\bfd}{{\mathbf{d}}}
\newcommand{\bfe}{{\mathbf{e}}}
\newcommand{\bfh}{{\mathbf{h}}}
\newcommand{\bfj}{{\mathbf{j}}}
\newcommand{\bfi}{{\mathbf{i}}}
\newcommand{\bfk}{{\mathbf{k}}}
\newcommand{\bfl}{{\mathbf{l}}}
\newcommand{\bfm}{{\mathbf{m}}}
\newcommand{\bfn}{{\mathbf{n}}}
\newcommand{\bfp}{{\mathbf{p}}}
\newcommand{\bfx}{{\mathbf{x}}}
\newcommand{\bfy}{{\mathbf{y}}}
\renewcommand{\k}{{\mathbf{k}}}
\newcommand{\n}{{\mathbf{n}}}
\newcommand{\vv}{{\mathbf{v}}}
\newcommand{\bv}{{\mathbf{v}}}
\newcommand{\w}{{\mathbf{w}}}
\newcommand{\x}{{\mathbf{x}}}
\newcommand{\y}{{\mathbf{y}}}
\newcommand{\cc}{{\mathbf{c}}}
\newcommand{\zero}{{\mathbf{0}}}
\newcommand{\one}{{\mathbf{1}}}
\newcommand{\bmm}{{\mathbf{m}}}

\newcommand{\ahat}{{\widehat{a}}}
\newcommand{\Zhat}{{\widehat{Z}}}

\newcommand{\C}{{\mathbb C}}
\newcommand{\D}{{\mathbb D}}
\newcommand{\Z}{{\mathbb Z}}
\newcommand{\N}{{\mathbb N}}
\newcommand{\Q}{{\mathbb Q}}
\newcommand{\PP}{{\mathbb P}}
\newcommand{\R}{{\mathbb R}}
\newcommand{\RR}{{\mathbb R}}
\newcommand{\E}{{\mathbb E}}

\newcommand{\Sym}{{\mathfrak{S}}}
\newcommand{\SymB}{{\mathfrak{B}}}
\newcommand{\Alt}{{\mathrm{Alt}}}

\newcommand{\germanA}{{\mathfrak{A}}}
\newcommand{\germanB}{{\mathfrak{B}}}
\newcommand{\germanQ}{{\mathfrak{Q}}}
\newcommand{\germanh}{{\mathfrak{h}}}

\newcommand{\myle}{\preceq}
\newcommand{\myge}{\succeq}
\newcommand{\mygt}{\succ}

\newcommand{\B}{{\sf B}}
\newcommand{\OB}{B^{\rm ord}}
\newcommand{\OS}{{\sf OS}}
\newcommand{\OO}{{\sf O}}
\newcommand{\SP}{{\sf SP}}
\newcommand{\OSP}{{\sf OSP}}
\newcommand{\Eu}{{\sf Eu}}
\newcommand{\ERR}{{\sf ERR}}
\newcommand{\sfB}{{\sf B}}
\newcommand{\sfD}{{\sf D}}
\newcommand{\sfE}{{\sf E}}
\newcommand{\sfF}{{\sf F}}
\newcommand{\sfG}{{\sf G}}
\newcommand{\sfJ}{{\sf J}}
\newcommand{\sfP}{{\sf P}}
\newcommand{\sfQ}{{\sf Q}}
\newcommand{\sfS}{{\sf S}}
\newcommand{\sfT}{{\sf T}}
\newcommand{\sfW}{{\sf W}}
\newcommand{\sfMV}{{\sf MV}}
\newcommand{\AMV}{{\sf AMV}}
\newcommand{\BM}{{\sf BM}}
\newcommand{\NC}{{\sf NC}}

\newcommand{\PFD}{\mathbf{PFD}}
\newcommand{\prope}{\mathrm{prope}}
\newcommand{\imprope}{\mathrm{imprope}}
\newcommand{\reg}{\mathrm{reg}}
\newcommand{\irreg}{\mathrm{irreg}}
\newcommand{\pc}{{\rm pc}}
\newcommand{\pdeg}{{\rm pdeg}}
\newcommand{\pindeg}{{\rm pindeg}}
\newcommand{\tre}{\textcolor{red}}
\newcommand{\relab}[2] {\rho_{+#1}^{-#2}}

\newcommand{\emIB}{B^{\rm irr}}
\newcommand{\emIP}{P^{\rm irr}}
\newcommand{\emOB}{B^{\rm ord}}
\newcommand{\emCB}{B^{\rm cyc}}
\newcommand{\emSC}{P^{\rm cyc}}

\newcommand{\stat}{{\rm stat}}
\newcommand{\cyc}{{\rm cyc}}
\newcommand{\Asc}{{\rm Asc}}
\newcommand{\asc}{{\rm asc}}
\newcommand{\Des}{{\rm Des}}
\newcommand{\des}{{\rm des}}
\newcommand{\Exc}{{\rm Exc}}
\newcommand{\exc}{{\rm exc}}
\newcommand{\Wex}{{\rm Wex}}
\newcommand{\wex}{{\rm wex}}
\newcommand{\Fix}{{\rm Fix}}
\newcommand{\fix}{{\rm fix}}
\newcommand{\lrmax}{{\rm lrmax}}
\newcommand{\rlmax}{{\rm rlmax}}
\newcommand{\Rec}{{\rm Rec}}
\newcommand{\rec}{{\rm rec}}
\newcommand{\Arec}{{\rm Arec}}
\newcommand{\arec}{{\rm arec}}
\newcommand{\ERec}{{\rm ERec}}
\newcommand{\erec}{{\rm erec}}
\newcommand{\EArec}{{\rm EArec}}
\newcommand{\earec}{{\rm earec}}
\newcommand{\recarec}{{\rm recarec}}
\newcommand{\nonrec}{{\rm nonrec}}
\newcommand{\Cpeak}{{\rm Cpeak}}
\newcommand{\cpeak}{{\rm cpeak}}
\newcommand{\Cval}{{\rm Cval}}
\newcommand{\cval}{{\rm cval}}
\newcommand{\Cdasc}{{\rm Cdasc}}
\newcommand{\cdasc}{{\rm cdasc}}
\newcommand{\Cddes}{{\rm Cddes}}
\newcommand{\cddes}{{\rm cddes}}
\newcommand{\cdrise}{{\rm cdrise}}
\newcommand{\cdfall}{{\rm cdfall}}
\newcommand{\Peak}{{\rm Peak}}
\newcommand{\peak}{{\rm peak}}
\newcommand{\Val}{{\rm Val}}
\newcommand{\val}{{\rm val}}
\newcommand{\Dasc}{{\rm Dasc}}
\newcommand{\dasc}{{\rm dasc}}
\newcommand{\Ddes}{{\rm Ddes}}
\newcommand{\ddes}{{\rm ddes}}
\newcommand{\inv}{{\rm inv}}
\newcommand{\maj}{{\rm maj}}
\newcommand{\rs}{{\rm rs}}
\newcommand{\cross}{{\rm cr}}
\newcommand{\crosshat}{{\widehat{\rm cr}}}
\newcommand{\nest}{{\rm ne}}
\newcommand{\rodd}{{\rm rodd}}
\newcommand{\reven}{{\rm reven}}
\newcommand{\lodd}{{\rm lodd}}
\newcommand{\leven}{{\rm leven}}
\newcommand{\sg}{{\rm sg}}
\newcommand{\bl}{{\rm bl}}
\newcommand{\tran}{{\rm tr}}
\newcommand{\area}{{\rm area}}
\newcommand{\ret}{{\rm ret}}
\newcommand{\peaks}{{\rm peaks}}
\newcommand{\hl}{{\rm hl}}
\newcommand{\sll}{{\rm sl}}
\newcommand{\negg}{{\rm neg}}
\newcommand{\imp}{{\rm imp}}
\newcommand{\osg}{{\rm osg}}
\newcommand{\ons}{{\rm ons}}
\newcommand{\isg}{{\rm isg}}
\newcommand{\ins}{{\rm ins}}
\newcommand{\LL}{{\rm LL}}
\newcommand{\height}{{\rm ht}}
\newcommand{\as}{{\rm as}}

\newcommand{\ba}{{\bm{a}}}
\newcommand{\bahat}{{\widehat{\bm{a}}}}
\newcommand{\sfa}{{{\sf a}}}
\newcommand{\bb}{{\bm{b}}}
\newcommand{\bc}{{\bm{c}}}
\newcommand{\bchat}{{\widehat{\bm{c}}}}
\newcommand{\bd}{{\bm{d}}}
\newcommand{\bee}{{\bm{e}}}
\newcommand{\beh}{{\bm{eh}}}
\newcommand{\bff}{{\bm{f}}}
\newcommand{\bg}{{\bm{g}}}
\newcommand{\bh}{{\bm{h}}}
\newcommand{\bll}{{\bm{\ell}}}
\newcommand{\bp}{{\bm{p}}}
\newcommand{\br}{{\bm{r}}}
\newcommand{\bs}{{\bm{s}}}
\newcommand{\bt}{{\bm{t}}}
\newcommand{\bu}{{\bm{u}}}
\newcommand{\bw}{{\bm{w}}}
\newcommand{\bx}{{\bm{x}}}
\newcommand{\by}{{\bm{y}}}
\newcommand{\bz}{{\bm{z}}}
\newcommand{\bA}{{\bm{A}}}
\newcommand{\bB}{{\bm{B}}}
\newcommand{\bC}{{\bm{C}}}
\newcommand{\bE}{{\bm{E}}}
\newcommand{\bF}{{\bm{F}}}
\newcommand{\bG}{{\bm{G}}}
\newcommand{\bH}{{\bm{H}}}
\newcommand{\bI}{{\bm{I}}}
\newcommand{\bJ}{{\bm{J}}}
\newcommand{\bM}{{\bm{M}}}
\newcommand{\bN}{{\bm{N}}}
\newcommand{\bP}{{\bm{P}}}
\newcommand{\bQ}{{\bm{Q}}}
\newcommand{\bR}{{\bm{R}}}
\newcommand{\bS}{{\bm{S}}}
\newcommand{\bT}{{\bm{T}}}
\newcommand{\bW}{{\bm{W}}}
\newcommand{\bX}{{\bm{X}}}
\newcommand{\bY}{{\bm{Y}}}
\newcommand{\bIB}{{\bm{B}^{\rm irr}}}
\newcommand{\bOB}{{\bm{B}^{\rm ord}}}
\newcommand{\bOS}{{\bm{OS}}}
\newcommand{\bERR}{{\bm{ERR}}}
\newcommand{\bSP}{{\bm{SP}}}
\newcommand{\bMV}{{\bm{MV}}}
\newcommand{\bBM}{{\bm{BM}}}
\newcommand{\balpha}{{\bm{\alpha}}}
\newcommand{\bbeta}{{\bm{\beta}}}
\newcommand{\bgamma}{{\bm{\gamma}}}
\newcommand{\bdelta}{{\bm{\delta}}}
\newcommand{\bkappa}{{\bm{\kappa}}}
\newcommand{\bmu}{{\bm{\mu}}}
\newcommand{\bomega}{{\bm{\omega}}}
\newcommand{\bphi}{{\bm{\phi}}}
\newcommand{\bsigma}{{\bm{\sigma}}}
\newcommand{\btau}{{\bm{\tau}}}
\newcommand{\bpsi}{{\bm{\psi}}}
\newcommand{\bzeta}{{\bm{\zeta}}}
\newcommand{\bone}{{\bm{1}}}
\newcommand{\bzero}{{\bm{0}}}

\newcommand{\phihat}{{\widehat{{\rule{0ex}{1.45ex}\mkern-2mu\smash{\phi}}}}}

\newcommand{\Cbar}{{\overline{C}}}
\newcommand{\Dbar}{{\overline{D}}}
\newcommand{\dbar}{{\overline{d}}}
\def\Ctilde{{\widetilde{C}}}
\def\Ftilde{{\widetilde{F}}}
\def\Gtilde{{\widetilde{G}}}
\def\Htilde{{\widetilde{H}}}
\def\Ptilde{{\widetilde{P}}}
\def\Phat{{\widehat{P}}}
\def\Chat{{\widehat{C}}}
\def\ctilde{{\widetilde{c}}}
\def\zbar{{\overline{Z}}}
\def\pitilde{{\widetilde{\pi}}}

\newcommand{\sech}{{\rm sech}}

%
%
\newcommand{\sn}{{\rm sn}}
\newcommand{\cn}{{\rm cn}}
\newcommand{\dn}{{\rm dn}}
\newcommand{\sm}{{\rm sm}}
\newcommand{\cm}{{\rm cm}}

%
%
\newcommand{\zfz}{ {{}_0 \! F_0} }
\newcommand{\zfo}{ {{}_0  F_1} }
\newcommand{\ofz}{ {{}_1 \! F_0} }
\newcommand{\ofo}{ {{}_1 \! F_1} }
\newcommand{\oft}{ {{}_1 \! F_2} }

%
%
\newcommand{\FHyper}[2]{ {\tensor[_{#1 \!}]{F}{_{#2}}\!} }
\newcommand{\FHYPER}[5]{ {\FHyper{#1}{#2} \!\biggl(
   \!\!\begin{array}{c} #3 \\[1mm] #4 \end{array}\! \bigg|\, #5 \! \biggr)} }
\newcommand{\tfo}{ {\FHyper{2}{1}} }
\newcommand{\tfz}{ {\FHyper{2}{0}} }
\newcommand{\threefz}{ {\FHyper{3}{0}} }
\newcommand{\FHYPERbottomzero}[3]{ {\FHyper{#1}{0} \hspace*{-0mm}\biggl(
   \!\!\begin{array}{c} #2 \\[1mm] \hbox{---} \end{array}\! \bigg|\, #3 \! \biggr)} }
\newcommand{\FHYPERtopzero}[3]{ {\FHyper{0}{#1} \hspace*{-0mm}\biggl(
   \!\!\begin{array}{c} \hbox{---} \\[1mm] #2 \end{array}\! \bigg|\, #3 \! \biggr)} }

\newcommand{\phiHyper}[2]{ {\tensor[_{#1}]{\phi}{_{#2}}} }
\newcommand{\psiHyper}[2]{ {\tensor[_{#1}]{\psi}{_{#2}}} }
\newcommand{\PhiHyper}[2]{ {\tensor[_{#1}]{\Phi}{_{#2}}} }
\newcommand{\PsiHyper}[2]{ {\tensor[_{#1}]{\Psi}{_{#2}}} }
\newcommand{\phiHYPER}[6]{ {\phiHyper{#1}{#2} \!\left(
   \!\!\begin{array}{c} #3 \\ #4 \end{array}\! ;\, #5, \, #6 \! \right)\!} }
\newcommand{\psiHYPER}[6]{ {\psiHyper{#1}{#2} \!\left(
   \!\!\begin{array}{c} #3 \\ #4 \end{array}\! ;\, #5, \, #6 \! \right)} }
\newcommand{\PhiHYPER}[5]{ {\PhiHyper{#1}{#2} \!\left(
   \!\!\begin{array}{c} #3 \\ #4 \end{array}\! ;\, #5 \! \right)\!} }
\newcommand{\PsiHYPER}[5]{ {\PsiHyper{#1}{#2} \!\left(
   \!\!\begin{array}{c} #3 \\ #4 \end{array}\! ;\, #5 \! \right)\!} }
\newcommand{\zerophizero}{ {\phiHyper{0}{0}} }
\newcommand{\ophizero}{ {\phiHyper{1}{0}} }
\newcommand{\zphio}{ {\phiHyper{0}{1}} }
\newcommand{\ophio}{ {\phiHyper{1}{1}} }
\newcommand{\tphio}{ {\phiHyper{2}{1}} }
\newcommand{\tphiz}{ {\phiHyper{2}{0}} }
\newcommand{\tPhio}{ {\PhiHyper{2}{1}} }
\newcommand{\opsio}{ {\psiHyper{1}{1}} }

%
%
\newcommand{\stirlingsubset}[2]{\genfrac{\{}{\}}{0pt}{}{#1}{#2}}
\newcommand{\stirlingcycleold}[2]{\genfrac{[}{]}{0pt}{}{#1}{#2}}
\newcommand{\stirlingcycle}[2]{\left[\! \stirlingcycleold{#1}{#2} \!\right]}
\newcommand{\assocstirlingsubset}[3]{{\genfrac{\{}{\}}{0pt}{}{#1}{#2}}_{\! \ge #3}}
\newcommand{\genstirlingsubset}[4]{{\genfrac{\{}{\}}{0pt}{}{#1}{#2}}_{\! #3,#4}}
\newcommand{\irredstirlingsubset}[2]{{\genfrac{\{}{\}}{0pt}{}{#1}{#2}}^{\!\rm irr}}
\newcommand{\euler}[2]{\genfrac{\langle}{\rangle}{0pt}{}{#1}{#2}}
\newcommand{\eulergen}[3]{{\genfrac{\langle}{\rangle}{0pt}{}{#1}{#2}}_{\! #3}}
\newcommand{\eulersecond}[2]{\left\langle\!\! \euler{#1}{#2} \!\!\right\rangle}
\newcommand{\eulersecondgen}[3]{{\left\langle\!\! \euler{#1}{#2} \!\!\right\rangle}_{\! #3}}
\newcommand{\binomvert}[2]{\genfrac{\vert}{\vert}{0pt}{}{#1}{#2}}
\newcommand{\binomsquare}[2]{\genfrac{[}{]}{0pt}{}{#1}{#2}}


\newenvironment{sarray}{
             \textfont0=\scriptfont0
             \scriptfont0=\scriptscriptfont0
             \textfont1=\scriptfont1
             \scriptfont1=\scriptscriptfont1
             \textfont2=\scriptfont2
             \scriptfont2=\scriptscriptfont2
             \textfont3=\scriptfont3
             \scriptfont3=\scriptscriptfont3
           \renewcommand{\arraystretch}{0.7}
           \begin{array}{l}}{\end{array}}

\newenvironment{scarray}{
             \textfont0=\scriptfont0
             \scriptfont0=\scriptscriptfont0
             \textfont1=\scriptfont1
             \scriptfont1=\scriptscriptfont1
             \textfont2=\scriptfont2
             \scriptfont2=\scriptscriptfont2
             \textfont3=\scriptfont3
             \scriptfont3=\scriptscriptfont3
           \renewcommand{\arraystretch}{0.7}
           \begin{array}{c}}{\end{array}}


\newcommand*\circled[1]{\tikz[baseline=(char.base)]{
  \node[shape=circle,draw,inner sep=1pt] (char) {#1};}}
\newcommand{\ostar}{{\circledast}}
\newcommand{\ostarN}{{\,\circledast_{\vphantom{\dot{N}}N}\,}}
\newcommand{\ostarPsi}{{\,\circledast_{\vphantom{\dot{\Psi}}\Psi}\,}}
\newcommand{\starN}{{\,\ast_{\vphantom{\dot{N}}N}\,}}
\newcommand{\starpsi}{{\,\ast_{\vphantom{\dot{\bpsi}}\!\bpsi}\,}}
\newcommand{\starone}{{\,\ast_{\vphantom{\dot{1}}1}\,}}
\newcommand{\startwo}{{\,\ast_{\vphantom{\dot{2}}2}\,}}
\newcommand{\starinfty}{{\,\ast_{\vphantom{\dot{\infty}}\infty}\,}}
\newcommand{\starT}{{\,\ast_{\vphantom{\dot{T}}T}\,}}

\newcommand*{\Scale}[2][4]{\scalebox{#1}{$#2$}}

\newcommand*{\Scaletext}[2][4]{\scalebox{#1}{#2}} 

\clearpage

\tableofcontents

\clearpage

\section{Introduction and statement of results}

It is well known \cite{Moon_70,Stanley_99}
that the number of rooted trees on the vertex set
$[n+1] \eqdef \{1,\ldots,n+1\}$ is $t_n = (n+1)^n$;
and it is also known (though perhaps less well~so)
\cite{Chauve_99,Chauve_00,Sokal_trees_enumeration}
that the number of rooted trees on the vertex set $[n+1]$
in~which exactly $k$ children of the root are lower-numbered than the root is
\be
   t_{n,k} \;=\;  \binom{n}{k} \, n^{n-k} 
   \;.
 \label{def.tnk}
\ee
The first few $t_{n,k}$ and $t_n$ are
%
%
\vspace*{-5mm}
\begin{table}[H]
\centering
\footnotesize
\begin{tabular}{c|rrrrrrrrr|r}
$n \setminus k$ & 0 & 1 & 2 & 3 & 4 & 5 & 6 & 7 & 8 & $(n+1)^n$ \\
\hline
0 & 1 &  &  &  &  &  &  &  &  & 1  \\
1 & 1 & 1 &  &  &  &  &  &  &  & 2  \\
2 & 4 & 4 & 1 &  &  &  &  &  &  & 9  \\
3 & 27 & 27 & 9 & 1 &  &  &  &  &  & 64  \\
4 & 256 & 256 & 96 & 16 & 1 &  &  &  &  & 625  \\
5 & 3125 & 3125 & 1250 & 250 & 25 & 1 &  &  &  & 7776  \\
6 & 46656 & 46656 & 19440 & 4320 & 540 & 36 & 1 &  &  & 117649  \\
7 & 823543 & 823543 & 352947 & 84035 & 12005 & 1029 & 49 & 1 &  & 2097152  \\
8 & 16777216 & 16777216 & 7340032 & 1835008 & 286720 & 28672 & 1792 & 64 & 1 & 43046721  \\
\end{tabular}
\end{table}
\vspace*{-5mm}

\noindent
\!\!\cite[A071207 and A000169]{OEIS}.

There is a second combinatorial interpretation of the numbers $t_{n,k}$,
also in terms of rooted trees:
namely, $t_{n,k}$ is the number of rooted trees on the vertex set $[n+1]$
in~which some specified vertex $i$ has $k$ children.\footnote{
   This fact ought to be well known, but to our surprise we have been unable
   to find any published reference.  Let us therefore give two proofs:

   {\sc First proof.}  Let $\scrt^\bullet_n$ denote the set of rooted trees
   on the vertex set $[n]$,
   and let $\deg_T(i)$ denote the number of children of the vertex~$i$
   in the rooted tree $T$.
   Rooted trees $T \in \scrt^\bullet_{n+1}$
   are associated bijectively to {\em Pr\"ufer sequences}\/
   $(s_1,\ldots,s_n) \in [n+1]^n$,
   in~which each index $i \in [n+1]$ appears $\deg_T(i)$ times
   \cite[pp.~25--26]{Stanley_99}.
   There are $\binom{n}{k} n^{n-k}$ sequences in~which
   the index $i$ appears exactly $k$ times.

   Equivalently, by \cite[Theorem~5.3.4, eq.~(5.47)]{Stanley_99},
$$
   \sum\limits_{T \in \scrt^\bullet_n} \, \prod_{j=1}^n x_j^{\deg_T(j)}
   \;=\;
   (x_1 + \ldots + x_n)^{n-1}
   \;.
$$
   Replacing $n \to n+1$ and then setting
   $x_i = x$ and $x_j = 1$ for $j \neq i$,
   we have
$$
   \sum\limits_{T \in \scrt^\bullet_{n+1}} \!\! x^{\deg_T(i)}
   \;=\;
   (x + n)^n
   \;.
$$
   Extracting the coefficient of $x^k$ yields $\binom{n}{k} n^{n-k}$.

   {\sc Second proof.}  There are $f_{n,k} = \binom{n}{k} \, k \, n^{n-k-1}$
   $k$-component forests of rooted trees on $n$ labeled vertices
   (see the references cited in \cite[footnote~1]{forests_totalpos}).
   By adding a new vertex 0 and connecting it to the roots of all the trees,
   we see that $f_{n,k}$ is also the number of unrooted trees
   on $n+1$ labeled vertices in~which some specified vertex (here vertex~0)
   has degree $k$.
   Now choose a root: if this root is 0, then vertex~0 has $k$ children;
   otherwise vertex~0 has $k-1$ children.  It follows that the number of
   rooted trees on $n+1$ labeled vertices in~which some specified vertex
   has $k$ children is
   $f_{n,k} + n f_{n,k+1} = \binom{n}{k} n^{n-k}$.

   The second proof was found independently by Ira Gessel
   (private communication).
}

And finally, there is a third combinatorial interpretation
of the numbers $t_{n,k}$ \cite{Critzer_12_OEIS}
that~is even simpler than the preceding two.
Recall first that a {\em functional digraph}\/
is a directed graph in~which every vertex has out-degree 1;
the terminology comes from the fact that such digraphs
are in obvious bijection with functions $f$ from the vertex set to itself
[namely, $\overrightarrow{ij}$ is an edge if and only~if $f(i) = j$].
Let us now define a \textbfit{partial functional digraph}
to be a directed graph in~which every vertex has out-degree 0~or~1;
and let us write $\PFD_{n,k}$ for the set of partial functional digraphs
on the vertex set $[n]$ in~which exactly $k$ vertices have out-degree~0.
(So $\PFD_{n,0}$ is the set of functional digraphs.)
A digraph in $\PFD_{n,k}$ has $n-k$ edges.
It is easy to see that $|\PFD_{n,k}| = t_{n,k}$:
there are $\binom{n}{k}$ choices for the out-degree-0 vertices,
and $n^{n-k}$ choices for the edges emanating from the remaining vertices.

%
%

We will use all three combinatorial models at various points in this paper.

The unit-lower-triangular matrix $(t_{n,k})_{n,k \ge 0}$
has the exponential generating function
\be
   \sum_{n=0}^\infty \sum_{k=0}^n t_{n,k} \, {t^n \over n!} \, x^k
   \;=\;
   {e^{x T(t)} \over 1 - T(t)}
 \label{eq.tnk.egf}
\ee
where
\be
   T(t)  \;\eqdef\; \sum_{n=1}^\infty n^{n-1} \, {t^n \over n!}
 \label{def.treefn}
\ee
is the {\em tree function}\/ \cite{Corless_96}.\footnote{
   In the analysis literature, expressions involving the tree function
   are often written in terms of the {\em Lambert $W$ function}\/
   $W(t) = -T(-t)$, which is the inverse function to $w \mapsto w e^w$
   \cite{Corless_96,Kalugin_12b}.
}
An equivalent statement is that
the unit-lower-triangular matrix $(t_{n,k})_{n,k \ge 0}$
is the exponential Riordan array
\cite{Deutsch_04,Deutsch_09,Barry_16,Shapiro_22}
$\scrr[F,G]$ with
$F(t) = \sum_{n=0}^\infty n^n \, t^n/n! = 1/[1-T(t)]$
and $G(t) = T(t)$;
we will discuss this connection in Section~\ref{subsec.EGF.1}.

The principal purpose of this paper is to prove the total positivity
of some matrices related to (and generalizing) $t_n$ and $t_{n,k}$.
Recall first that a finite or infinite matrix of real numbers is called
{\em totally positive}\/ (TP) if all its minors are nonnegative,
and {\em strictly totally positive}\/ (STP)
if all its minors are strictly positive.\footnote{
   {\bf Warning:}  Many authors
   (e.g.\ \cite{Gantmakher_37,Gantmacher_02,Fomin_00,Fallat_11})
   use the terms ``totally nonnegative'' and ``totally positive''
   for what we have termed ``totally positive'' and
   ``strictly totally positive'', respectively.
   So it is very important, when seeing any claim about
   ``totally positive'' matrices, to ascertain~which sense of
   ``totally positive'' is being used!
   (This is especially important because many theorems in this subject
    require {\em strict}\/ total positivity for their validity.)
}
Background information on totally positive matrices can be found
in \cite{Karlin_68,Gantmacher_02,Pinkus_10,Fallat_11};
they have applications to many areas of pure and applied mathematics.\footnote{
   Including combinatorics
   \cite{Brenti_89,Brenti_95,Brenti_96,Fomin_00,Skandera_03},
   stochastic processes \cite{Karlin_59,Karlin_68},
   statistics \cite{Karlin_68},
   the mechanics of oscillatory systems \cite{Gantmakher_37,Gantmacher_02},
   the zeros of polynomials and entire functions
   \cite{Karlin_68,Asner_70,Kemperman_82,Holtz_03,Pinkus_10,Dyachenko_14},
   spline interpolation \cite{Schoenberg_53,Karlin_68,Gasca_96},
   Lie theory \cite{Lusztig_94,Lusztig_98,Fomin_99,Lusztig_08}
   and cluster algebras \cite{Fomin_10,Fomin_forthcoming},
   the representation theory of the infinite symmetric group
   \cite{Thoma_64,Borodin_17},
   the theory of immanants \cite{Stembridge_91},
   planar discrete potential theory \cite{Curtis_98,Fomin_01}
   and the planar Ising model \cite{Lis_17},
   and several other areas \cite{Gasca_96}.
}

Our first result is the following:

\begin{theorem}
   \label{thm1.1}
\quad\hfill\vspace*{-1mm}
\begin{itemize}
   \item[(a)]  The unit-lower-triangular matrix $\sfT = (t_{n,k})_{n,k \ge 0}$
       is totally positive.
   \item[(b)]  The Hankel matrix
       $H_\infty(\bt^{(0)}) = (t_{n+n',0})_{n,n' \ge 0}$
       is totally positive.
\end{itemize}
\end{theorem}

It is known \cite{Gantmakher_37,Pinkus_10}
that a Hankel matrix of real numbers is totally positive
if and only if the underlying sequence
is a Stieltjes moment sequence,
i.e.\ the moments of a positive measure on $[0,\infty)$.
And it is also known
that $(n^n)_{n \ge 0}$
is a Stieltjes moment sequence.\footnote{
   The integral representation
   \cite{Bouwkamp_86} \cite[Corollary~2.4]{Kalugin_12b}
   $$
      {n^n \over n!}
      \;=\;
      {1 \over \pi}
      \int\limits_0^\pi
          \biggl( {\sin\nu \over \nu} \, e^{\nu \cot\nu} \biggr) ^{\! n}
                 \: d\nu
   $$
   shows that $n^n/n!$
   is a Stieltjes moment sequence.
   Moreover, $n! = \int_0^\infty x^n \, e^{-x} \, dx$
   is a Stieltjes moment sequence.
   Since the entrywise product of two Stieltjes moment sequences
   is easily seen to be a Stieltjes moment sequence,
   it follows that $n^n$
   is a Stieltjes moment sequence.
   But we do not know any simple formula
   (i.e.\ one involving only a single integral over a real variable)
   for its Stieltjes integral representation.
}
So Theorem~\ref{thm1.1}(b) is equivalent to this known result.
But our proof here is combinatorial and linear-algebraic, not analytic.

However, this is only the beginning of the story,
because our main interest \cite{Sokal_flajolet,Sokal_OPSFA,Sokal_totalpos}
is not with sequences and matrices of real numbers,
but rather with sequences and matrices of {\em polynomials}\/
(with integer or real coefficients) in one or more indeterminates $\bfx$:
in applications they will typically be generating polynomials that enumerate
some combinatorial objects with respect to one or more statistics.
We equip the polynomial ring $\R[\bfx]$ with the coefficientwise
partial order:  that~is, we say that $P$ is nonnegative
(and write $P \myge 0$)
in~case $P$ is a polynomial with nonnegative coefficients.
We then say that a matrix with entries in $\R[\bfx]$ is
\textbfit{coefficientwise totally positive}
if all its minors are polynomials with nonnegative coefficients;
and we say that a sequence $\ba = (a_n)_{n \ge 0}$ with entries in $\R[\bfx]$
is \textbfit{coefficientwise Hankel-totally positive}
if its associated infinite Hankel matrix
$H_\infty(\ba) = (a_{n+n'})_{n,n' \ge 0}$
is coefficientwise totally positive.

Returning now to the matrix $\sfT = (t_{n,k})_{n,k \ge 0}$,
let us define its \textbfit{row-generating polynomials} in the usual way:
\be
   T_n(x)  \;=\;  \sum_{k=0}^n t_{n,k} \, x^k
   \;.
 \label{def.Tn}
\ee
{}From the definition \reff{def.tnk} we obtain the explicit formula
\be
   T_n(x)  \;=\;  (x+n)^n
   \;.
 \label{eq.Tn.explicit}
\ee
Our second result is then:

\begin{theorem}
   \label{thm1.2}
The polynomial sequence $\bT = \bigl( T_n(x) \bigr)_{n \ge 0}$
is coefficientwise Hankel-totally positive.
[That is, the Hankel matrix
 $H_\infty(\bT) = \bigl( T_{n+n'}(x) \bigr)_{n,n' \ge 0}$
 is coefficientwise totally positive.]
\end{theorem}

Theorem~\ref{thm1.2} strengthens Theorem~\ref{thm1.1}(b),
and reduces to it when $x=0$.
The proof of Theorem~\ref{thm1.2} will be based on studying the
\textbfit{binomial row-generating matrix} $\sfT B_x$,
where $B_x$ is the weighted binomial matrix
\be
   (B_x)_{ij}  \;=\;  \binom{i}{j} \, x^{i-j}
 \label{def.Bx}
\ee
(see Sections~\ref{subsec.rowgen} and \ref{subsec.expriordan}).

\bigskip

But this is not the end of the story, because we want to
generalize these polynomials further by adding further variables.
Given a rooted tree $T$ and two vertices $i,j$ of $T$,
we say that $j$ is a {\em descendant}\/ of $i$
if the unique path from the root of $T$ to $j$ passes through $i$.
(Note in particular that every vertex is a descendant of itself.)
Now suppose that the vertex set of $T$ is totally ordered
(for us it will be $[n+1]$),
and let $e = ij$ be an edge of $T$,
ordered so that $j$ is a descendant of $i$.
We say that the edge $e = ij$ is \textbfit{improper}
if there exists a descendant of $j$ (possibly $j$ itself)
that~is lower-numbered than $i$;
otherwise we say that $e = ij$ is \textbfit{proper}.
We denote by $\imprope(T)$ [resp.~$\prope(T)$]
the number of improper (resp.~proper) edges in the tree $T$.

We now introduce these statistics into our second combinatorial model.
Let $\scrt^{\<i;k\>}_n$ denote the set of rooted trees on the vertex set $[n]$
in~which the vertex $i$ has $k$ children.
For the identity $|\scrt^{\<i;k\>}_{n+1}| = t_{n,k}$,
we can use any $i \in [n+1]$;
but for the following we specifically want to take $i=1$.
With this choice we observe that the $k$ edges from the vertex~1
to its children are automatically proper.
We therefore define
\be
   t_{n,k}(y,z)
   \;=\;
   \sum_{T \in \scrt^{\<1;k\>}_{n+1}} y^{\imprope(T)} z^{\prope(T) -k}
   \;.
 \label{def.tnkyz}
\ee
Clearly $t_{n,k}(y,z)$ is a homogeneous polynomial of degree~$n-k$
with nonnegative integer coefficients;
it is a polynomial refinement of $t_{n,k}$
in the sense that $t_{n,k}(1,1) = t_{n,k}$.
(Of~course, it was redundant to introduce the two variables $y$ and $z$
 instead of just one of them; we did it because it makes the formulae
 more symmetric.)
The first few polynomials $t_{n,k}(y,1)$ are
\vspace*{-5mm}
\begin{table}[H]
\hspace*{-4mm}
\small
\begin{tabular}{c|lllll}
$n \setminus k$ & 0 & 1 & 2 & 3 & 4 \\
\hline
 0 &    $1$  &      &      &      &        \\
 1 &    $y$  &   $1$  &      &      &        \\
 2 &    $y + 3 y^2$  &  $1 + 3 y$  &   $1$  &      &        \\
 3 &    $2 y + 10 y^2 + 15 y^3$  &   $2 + 10 y + 15 y^2$  &   $3 + 6 y$  &   $1$  &        \\
 4 &    $6 y + 40 y^2 + 105 y^3 + 105 y^4$  &   $6 + 40 y + 105 y^2 + 105 y^3$  &   $11 + 40 y + 45 y^2$  &  $6 + 10 y$  &   1
\end{tabular}
\end{table}
\vspace{-5mm}
\noindent
The coefficient matrix of the zeroth-column polynomials $t_{n,0}(y,1)$
is \cite[A239098/A075856]{OEIS}.
This table also suggests the following result,
for which we will give a bijective proof:

\begin{proposition}
   \label{prop.tn01yz}
For $n \ge 1$, $t_{n,0}(y,z) = y \, t_{n,1}(y,z)$.
\end{proposition}

In Section~\ref{subsec.EGF.2} we will show that
the unit-lower-triangular matrix
$\sfT(y,z) =$ \linebreak $\bigl( t_{n,k}(y,z) \bigr)_{n,k \ge 0}$
is an exponential Riordan array $\scrr[F,G]$,
and we will compute $F(t)$ and $G(t)$.

We now generalize \reff{def.Tn} by defining the row-generating polynomials
\be
   T_n(x,y,z)  \;=\;  \sum_{k=0}^n t_{n,k}(y,z) \, x^k
 \label{def.Tnxyz}
\ee
or in other words
\be
   T_n(x,y,z)
   \;=\;
   \sum_{T \in \scrt^\bullet_{n+1}}
            x^{\deg_T(1)} y^{\imprope(T)} z^{\prope(T) - \deg_T(1)}
 \label{def.Tnxyz.bis}
\ee
where $\deg_T(1)$ is the number of children of the vertex~1
in the rooted tree $T$.
Note that $T_n(x,y,z)$ is a homogeneous polynomial of degree~$n$ in $x,y,z$,
with nonnegative integer coefficients;
it reduces to $T_n(x)$ when $y=z=1$.
Our third result is then:

\begin{theorem}
   \label{thm1.4}
\quad\hfill\vspace*{-1mm}
\begin{itemize}
   \item[(a)]  The unit-lower-triangular polynomial matrix
      $\sfT(y,z) = \bigl( t_{n,k}(y,z) \bigr)_{n,k \ge 0}$
      is coefficientwise totally positive (jointly in $y,z$).
   \item[(b)]  The polynomial sequence
       $\bT = \bigl( T_n(x,y,z) \bigr)_{n \ge 0}$
       is coefficientwise Hankel-totally positive (jointly in $x,y,z$).
\end{itemize}
\end{theorem}

Theorem~\ref{thm1.4} strengthens Theorems~\ref{thm1.1}(a) and \ref{thm1.2},
and reduces to them when $y=z=1$.
The proof of Theorem~\ref{thm1.4}(b) will be based on
studying the binomial row-generating matrix $\sfT(y,z) B_x$,
using the representation of $\sfT(y,z)$ as an exponential Riordan array.

%

\bigskip

Finally, let us consider our third combinatorial model,
which is based on partial functional digraphs.
Recall that a {\em functional digraph}\/
(resp.~{\em partial functional digraph}\/)
is a directed graph in~which every vertex has out-degree 1
(resp.~0 or 1).
Each weakly connected component of a functional digraph
consists of a directed cycle (possibly of length~1, i.e.~a loop)
together with a collection of (possibly trivial) directed trees
rooted at the vertices of the cycle (with edges pointing towards the root).
The weakly connected components of a partial functional digraph
are trees rooted at the out-degree-0 vertices
(with edges pointing towards the root)
together with components of the same form as in a functional digraph.
We say that a vertex of a partial functional digraph is
\textbfit{recurrent} (or \textbfit{cyclic}) if it lies on one of the cycles;
otherwise we call it \textbfit{transient} (or \textbfit{acyclic}).
If $j$ and $k$ are vertices of a digraph,
we say that $k$ is a \textbfit{predecessor} of $j$
if there exists a directed path from $k$ to~$j$
(in particular, every vertex is a predecessor of itself).\footnote{
   In a functional digraph,
   Dumont and Ramamonjisoa \cite[p.~11]{Dumont_96}
   use the term ``ascendance'', and the notation $A(j)$,
   to denote the set of all predecessors of $j$.
}
Note that ``predecessor'' in a digraph generalizes the notion of
``descendant'' in a rooted tree, if we make the convention that
all edges in the tree are oriented towards the root.
Indeed, if $j$ is a transient vertex in a partial functional digraph,
then the predecessors of $j$ are precisely the descendants of $j$
in the rooted tree
(rooted at either a recurrent vertex or an out-degree-0 vertex)
to which $j$ belongs.
On the other hand, if $j$ is a recurrent vertex,
then the predecessors of $j$ are all the vertices
in the weakly connected component containing $j$.

Now consider a partial functional digraph on a totally ordered vertex set
(which for us will be $[n]$).
We say that an edge $\overrightarrow{ji}$ (pointing from $j$ to $i$)
is \textbfit{improper}
if there exists a predecessor of $j$ (possibly $j$ itself)
that~is $\le i$;
otherwise we say that the edge $\overrightarrow{ji}$ is \textbfit{proper}.
When $j$ is a transient vertex, this coincides with
the notion of improper/proper edge in a rooted tree.
When $j$ is a recurrent vertex,
the edge $\overrightarrow{ji}$ is always improper,
because one of the predecessors of $j$ is $i$.
(This includes the case $i=j$: a loop is always an improper edge.)
We denote by $\imprope(G)$ [resp.~$\prope(G)$]
the number of improper (resp.~proper) edges
in the partial functional digraph $G$.
We then define the generating polynomial
\be
   \widetilde{t}_{n,k}(y,z)
   \;=\;
   \sum_{G \in \PFD_{n,k}} \!\! y^{\imprope(G)} z^{\prope(G)}
   \;.
 \label{def.tnkyz.tilde}
\ee
Since $G \in \PFD_{n,k}$ has $n-k$ edges,
$\widetilde{t}_{n,k}(y,z)$ is a homogeneous polynomial of degree~$n-k$
with nonnegative integer coefficients.
By bijection between our second and third combinatorial models, we will prove:

\begin{proposition}
   \label{prop.equivalence}
$\quad t_{n,k}(y,z)  \:=\:  \widetilde{t}_{n,k}(y,z) \:.$
\end{proposition}

The row-generating polynomials \reff{def.Tnxyz}/\reff{def.Tnxyz.bis}
thus have the alternate combinatorial interpretation
\be
   T_n(x,y,z)
   \;=\;
   \sum_{G \in \PFD_n} x^{\deg0(G)} y^{\imprope(G)} z^{\prope(G)}
 \label{def.Tnxyz.bisbis}
\ee
where $\deg0(G)$ is the number of out-degree-0 vertices in $G$.


We also have an interpretation
of the polynomials $t_{n,k}(y,z)$ in our first combinatorial model
(rooted trees in~which the root has $k$ lower-numbered children);
but since this interpretation is rather complicated,
we defer it to Appendix~\ref{app.first.tnkyz}.

\bigskip

But this is {\em still}\/ not the end of the story,
because we can add even more variables
into our second combinatorial model --- in fact, an infinite set.
Given a rooted tree $T$ on a totally ordered vertex set
and vertices $i,j \in T$ such that $j$ is a child of~$i$,
we say that $j$ is a \textbfit{proper child} of $i$
if the edge $e = ij$ is proper
(that~is, $j$ and all its descendants are higher-numbered than~$i$).
Now let $\bphi = (\phi_m)_{m \ge 0}$ be indeterminates,
and let $t_{n,k}(y,\bphi)$ be the generating polynomial for
rooted trees $T \in \scrt^{\<1;k\>}_{n+1}$
with a weight $y$ for each improper edge
and a weight $\phihat_m \eqdef m! \, \phi_m$
for each vertex $i \neq 1$ that has $m$ proper children:
\be
   t_{n,k}(y,\bphi)
   \;=\;
   \sum_{T \in \scrt^{\<1;k\>}_{n+1}} y^{\imprope(T)}
      \, \prod_{i=2}^{n+1} \phihat_{\pdeg_T(i)} 
 \label{def.tnkyphi}
\ee
where $\pdeg_T(i)$ denotes the number of {\em proper}\/ children
of the vertex~$i$ in the rooted tree $T$.
We will see later why it is convenient to introduce the factors~$m!$
in this definition.
Observe also that the variables~$z$ are now redundant and therefore omitted,
because they would simply scale $\phi_m \to z^m \phi_m$.
And note finally that, in conformity with \reff{def.tnkyz},
we have chosen to suppress the weight $\phihat_k$
that would otherwise be associated to the vertex~1.
We call the polynomials $t_{n,k}(y,\bphi)$
the \textbfit{generic rooted-tree polynomials},
and the lower-triangular matrix
$\sfT(y,\bphi) = \bigl( t_{n,k}(y,\bphi) \bigr)_{n,k \ge 0}$
the \textbfit{generic rooted-tree matrix}.
Here $\bphi = (\phi_m)_{m \ge 0}$ are in the first instance indeterminates,
so that $t_{n,k}(y,\bphi)$ belongs to the polynomial ring $\Z[y,\bphi]$;
but we can then, if we wish, substitute specific values for $\bphi$
in any commutative ring $R$, leading to values $t_{n,k}(y,\bphi) \in R[y]$.
(Similar substitutions can of course also be made for $y$.)
When doing this we will use the same notation $t_{n,k}(y,\bphi)$,
as the desired interpretation for $\bphi$ should be clear from the context.
The polynomial $t_{n,k}(y,\bphi)$ is homogeneous of degree~$n$ in $\bphi$;
it is also quasi-homogeneous of degree~$n-k$ in $y$ and $\bphi$
when $\phi_m$ is assigned weight $m$ and $y$ is assigned weight 1.
By specializing $t_{n,k}(y,\bphi)$ to
$\phi_m = z^m/m!$ and hence $\phihat_m = z^m$, we recover $t_{n,k}(y,z)$.

We remark that the matrix $\sfT(y,\bphi)$, unlike $\sfT(y,z)$,
is not {\em unit}\/-lower-triangular:
rather, it has diagonal entries $t_{n,n}(y,\bphi) = \phi_0^n$,
corresponding to the tree in~which 1~is the root
and has all the vertices $2,\ldots,n+1$ as children.
More generally, the polynomial $t_{n,k}(y,\bphi)$ is divisible by $\phi_0^k$,
since the vertex~1 always has at least $k$ leaf descendants.
So we could define a unit-lower-triangular matrix
$\sfT^\flat(y,\bphi) = \bigl( t_{n,k}^\flat(y,\bphi) \bigr)_{n,k \ge 0}$
by $t_{n,k}^\flat(y,\bphi) = t_{n,k}(y,\bphi)/\phi_0^k$.
(Alternatively, we could simply choose to normalize to $\phi_0 = 1$.)

In Section~\ref{subsec.EGF.3} we will show that $\sfT(y,\bphi)$
is an exponential Riordan array $\scrr[F,G]$,
and we will compute $F(t)$ and $G(t)$.

Also, generalizing Proposition~\ref{prop.tn01yz}, we will prove:

\begin{proposition}
   \label{prop.tn01yphi}
For $n \ge 1$, $t_{n,0}(y,\bphi) = y \, t_{n,1}(y,\bphi)$.
\end{proposition}

We can also define the corresponding polynomials $\widetilde{t}_{n,k}(y,\bphi)$
in the partial-functional-digraph model, as follows:
If $G$ is a partial functional digraph on a totally ordered vertex set,
and $i$ is a vertex of $G$,
we define the \textbfit{proper in-degree} of $i$, $\pindeg_G(i)$,
to be the number of proper edges $\overrightarrow{ji}$ in $G$.
We then define
\be
   \widetilde{t}_{n,k}(y,\bphi)
   \;=\;
   \sum_{G \in \PFD_{n,k}} \!\! y^{\imprope(G)}
      \, \prod_{i=1}^{n} \phihat_{\pindeg_G(i)} 
   \;.
 \label{def.tnkyphi.tilde}
\ee
Then, generalizing Proposition~\ref{prop.equivalence}, we will prove:

\begin{proposition}
   \label{prop.equivalence.phi}
$\quad t_{n,k}(y,\bphi)  \:=\:  \widetilde{t}_{n,k}(y,\bphi) \:.$
\end{proposition}

\medskip

Now define the row-generating polynomials
\be
   T_n(x,y,\bphi)  \;=\;  \sum_{k=0}^n t_{n,k}(y,\bphi) \, x^k
 \label{def.Fnyz.phi}
\ee
or in other words
\be
   T_n(x,y,\bphi)
   \;=\;
   \sum_{T \in \scrt^\bullet_{n+1}}
            x^{\deg_T(1)} y^{\imprope(T)}
            \prod_{i=2}^{n+1} \phihat_{\pdeg_T(i)} 
   \;.
 \label{def.Fnyz.phi.bis}
\ee
The main result of this paper is then the following:

\begin{theorem}
   \label{thm1.7}
Fix $1 \le r \le \infty$.
Let $R$ be a partially ordered commutative ring,
and let $\bphi = (\phi_m)_{m \ge 0}$ be a sequence in $R$
that~is Toeplitz-totally positive of order~$r$.
Then:
\begin{itemize}
      \item[(a)]  The lower-triangular polynomial matrix
      $\sfT(y,\bphi) = \bigl( t_{n,k}(y,\bphi) \bigr)_{n,k \ge 0}$
      is coefficientwise totally positive of order~$r$ (in $y$).
   \item[(b)]  The polynomial sequence
       $\bT = \bigl( T_{n}(x,y,\bphi) \bigr)_{n \ge 0}$
       is coefficientwise Hankel-totally positive of order~$r$
       (jointly in $x,y$).
\end{itemize}
\end{theorem}

\noindent
(The concept of Toeplitz-total positivity
 in a partially ordered commutative ring
 will be explained in detail in Section~\ref{subsec.totalpos.prelim}.
 Total positivity of order~$r$ means that the minors of size $\le r$
 are nonnegative.)
Specializing Theorem~\ref{thm1.7} to
$r = \infty$, $R = \Q$ and $\phi_m = z^m/m!$
(which is indeed Toeplitz-totally positive: see \reff{eq.thm.aissen} below),
we recover Theorem~\ref{thm1.4}.
The method of proof of Theorem~\ref{thm1.7}
will, in fact, be the same as that of Theorem~\ref{thm1.4},
suitably generalized.

\bigskip

We now give an overview of the contents of this paper.
The main tool in our proofs will be the theory of production matrices
\cite{Deutsch_05,Deutsch_09}
as applied to total positivity \cite{Sokal_totalpos},
combined with the theory of exponential Riordan arrays
\cite{Deutsch_04,Deutsch_09,Barry_16,Shapiro_22}.
Therefore, in Section~\ref{sec.prelim} we review some facts
about total positivity, production matrices and exponential Riordan arrays
that will play a central role in our arguments.
This development culminates in Corollary~\ref{cor.EAZ.hankel};
it is the fundamental theoretical result that underlies all our proofs.
In Section~\ref{sec.bijective} we give bijective proofs of
Propositions~\ref{prop.tn01yz}, \ref{prop.equivalence}, \ref{prop.tn01yphi}
and \ref{prop.equivalence.phi}.
In Section~\ref{sec.EGF} we show that the matrices
$\sfT$, $\sfT(y,z)$ and $\sfT(y,\bphi)$
are exponential Riordan arrays $\scrr[F,G]$,
and we compute their generating functions $F$ and $G$.
In Section~\ref{sec.proofs} we combine the results
of Sections~\ref{sec.prelim} and \ref{sec.EGF}
to complete the proofs of
Theorems~\ref{thm1.1}, \ref{thm1.2}, \ref{thm1.4} and \ref{thm1.7}.

This paper is a sequel to our paper \cite{forests_totalpos}
on the total positivity of matrices that enumerate
forests of rooted labeled trees.
The methods here are basically the same as in this previous paper,
but generalized nontrivially to handle
exponential Riordan arrays $\scrr[F,G]$ with $F \neq 1$.
Zhu \cite{Zhu_21a,Zhu_22} has employed closely related methods.
See also Gilmore \cite{Gilmore_21} for some total-positivity
results for $q$-generalizations of tree and forest matrices,
using very different methods.

\section{Preliminaries}   \label{sec.prelim}

Here we review some definitions and results from
\cite{Sokal_totalpos,latpath_lah}
that will be needed in the sequel.
We also include a brief review of
ordinary and exponential Riordan arrays
\cite{Shapiro_91,Sprugnoli_94,Deutsch_04,Deutsch_09,Barry_16,Shapiro_22}
and Lagrange inversion \cite{Gessel_16}.

The treatment of exponential Riordan arrays
in Section~\ref{subsec.expriordan} contains one novelty:
namely, the rewriting of the production matrix
in terms of new series $\Phi$ and $\Psi$
(see \reff{prop.EAZ.Phi.first} ff.\ and Proposition~\ref{prop.EAZ.Phi.first}).
This is the key step that leads to Corollary~\ref{cor.EAZ.hankel}.

\subsection{Partially ordered commutative rings and total positivity}
   \label{subsec.totalpos.prelim}

In this paper all rings will be assumed to have an identity element 1
and to be nontrivial ($1 \ne 0$).

A \textbfit{partially ordered commutative ring} is a pair $(R,\scrp)$ where
$R$ is a commutative ring and $\scrp$ is a subset of $R$ satisfying
\begin{itemize}
   \item[(a)]  $0,1 \in \scrp$.
   \item[(b)]  If $a,b \in \scrp$, then $a+b \in \scrp$ and $ab \in \scrp$.
   \item[(c)]  $\scrp \cap (-\scrp) = \{0\}$.
\end{itemize}
We call $\scrp$ the {\em nonnegative elements}\/ of $R$,
and we define a partial order on $R$ (compatible with the ring structure)
by writing $a \le b$ as a synonym for $b-a \in \scrp$.
Please note that, unlike the practice in real algebraic geometry
\cite{Brumfiel_79,Lam_84,Prestel_01,Marshall_08},
we do {\em not}\/ assume here that squares are nonnegative;
indeed, this property fails completely for our prototypical example,
the ring of polynomials with the coefficientwise order,
since $(1-x)^2 = 1-2x+x^2 \not\myge 0$.

Now let $(R,\scrp)$ be a partially ordered commutative ring
and let $\bfx = \{x_i\}_{i \in I}$ be a collection of indeterminates.
In the polynomial ring $R[\bfx]$ and the formal-power-series ring $R[[\bfx]]$,
let $\scrp[\bfx]$ and $\scrp[[\bfx]]$ be the subsets
consisting of polynomials (resp.\ series) with nonnegative coefficients.
Then $(R[\bfx],\scrp[\bfx])$ and $(R[[\bfx]],\scrp[[\bfx]])$
are partially ordered commutative rings;
we refer to this as the \textbfit{coefficientwise order}
on $R[\bfx]$ and $R[[\bfx]]$.

A (finite or infinite) matrix with entries in a
partially ordered commutative ring
is called \textbfit{totally positive} (TP) if all its minors are nonnegative;
it is called \textbfit{totally positive of order~$\bm{r}$} (TP${}_r$)
if all its minors of size $\le r$ are nonnegative.
It follows immediately from the Cauchy--Binet formula that
the product of two TP (resp.\ TP${}_r$) matrices is TP
(resp.\ TP${}_r$).\footnote{
   For infinite matrices, we need some condition to ensure that
   the product is well-defined.
   For instance, the product $AB$ is well-defined whenever
   $A$ is row-finite (i.e.\ has only finitely many nonzero entries in each row)
   or $B$ is column-finite.
}
This fact is so fundamental to the theory of total positivity
that we shall henceforth use it without comment.

We say that a sequence $\ba = (a_n)_{n \ge 0}$
with entries in a partially ordered commutative ring
is \textbfit{Hankel-totally positive} 
(resp.\ \textbfit{Hankel-totally positive of order~$\bm{r}$})
if its associated infinite Hankel matrix
$H_\infty(\ba) = (a_{i+j})_{i,j \ge 0}$
is TP (resp.\ TP${}_r$).
We say that $\ba$
is \textbfit{Toeplitz-totally positive} 
(resp.\ \textbfit{Toeplitz-totally positive of order~$\bm{r}$})
if its associated infinite Toeplitz matrix
$T_\infty(\ba) = (a_{i-j})_{i,j \ge 0}$
(where $a_n \eqdef 0$ for $n < 0$)
is TP (resp.\ TP${}_r$).\footnote{
   When $R = \R$, Toeplitz-totally positive sequences are traditionally called
   {\em P\'olya frequency sequences}\/ (PF),
   and Toeplitz-totally positive sequences of order $r$
   are called {\em P\'olya frequency sequences of order $r$}\/ (PF${}_r$).
   See \cite[chapter~8]{Karlin_68} for a detailed treatment.
}

When $R = \R$, Hankel- and Toeplitz-total positivity have simple
analytic characterizations.
A sequence $(a_n)_{n \ge 0}$ of real numbers
is Hankel-totally positive if and only if it is a Stieltjes moment sequence
\cite[Th\'eor\`eme~9]{Gantmakher_37} \cite[section~4.6]{Pinkus_10}.
And a sequence $(a_n)_{n \ge 0}$ of real numbers 
is Toeplitz-totally positive if and only if its ordinary generating function
can be written as
\be
   \sum_{n=0}^\infty a_n t^n
   \;=\;
   C e^{\gamma t} t^m \prod_{i=1}^\infty {1 + \alpha_i t  \over  1 - \beta_i t}
 \label{eq.thm.aissen}
\ee
with $m \in \N$, $C,\gamma,\alpha_i,\beta_i \ge 0$,
$\sum \alpha_i < \infty$ and $\sum \beta_i < \infty$:
this is the celebrated Aissen--Schoenberg--Whitney--Edrei theorem
\cite[Theorem~5.3, p.~412]{Karlin_68}.
However, in a general partially ordered commutative ring $R$,
the concepts of Hankel- and Toeplitz-total positivity are more subtle.

We will need a few easy facts about the total positivity of special matrices:

\begin{lemma}[Bidiagonal matrices]
  \label{lemma.bidiagonal}
Let $A$ be a matrix with entries in a partially ordered commutative ring,
with the property that all its nonzero entries belong to two consecutive
diagonals.
Then $A$ is totally positive if and only if all its entries are nonnegative.
\end{lemma}

\proof
The nonnegativity of the entries (i.e.\ TP${}_1$)
is obviously a necessary condition for TP.
Conversely, for a matrix of this type it is easy to see that
every nonzero minor is simply a product of some entries.
\qed

\begin{lemma}[Toeplitz matrix of powers]
   \label{lemma.toeplitz.power}
Let $R$ be a partially ordered commutative ring, let $x \in R$,
and consider the infinite Toeplitz matrix
\be
   T_x
   \;\eqdef\;
   T_\infty(x^\N)
   \;=\;
   \begin{bmatrix}
      1   &     &     &     &         \\
      x   &  1  &     &     &            \\
      x^2 &  x  &  1  &     &            \\
      x^3 &  x^2  &  x  &  1    &            \\
      \vdots  & \vdots   & \vdots  & \vdots   & \ddots
   \end{bmatrix}
   \;.
 \label{def.Tx}
\ee
Then every minor of $T_x$ is either zero or else a power of $x$.
Hence $T_x$ is TP $\iff$ $T_x$ is TP${}_1$ $\iff$ $x \ge 0$.

In particular, if $x$ is an indeterminate, then $T_x$ 
is totally positive in the ring $\Z[x]$ equipped with the
coefficientwise order.
\end{lemma}

\proof
Consider a submatrix $A = (T_x)_{IJ}$
with rows $I = \{i_1 < \ldots < i_k \}$
and columns $J = \{j_1 < \ldots < j_k \}$.
We will prove by induction on $k$ that
$\det A$ is either zero or a power of $x$.
It is trivial if $k=0$ or 1.
If $A_{12} = A_{22} = 0$,
then $A_{1s} = A_{2s} = 0$ for all $s \ge 2$ by definition of $T_x$,
and $\det A = 0$.
If $A_{12}$ and $A_{22}$ are both nonzero,
then the first column of $A$ is $x^{j_2 - j_1}$ times the second column,
and again $\det A = 0$.
Finally, if $A_{12} = 0$ and $A_{22} \ne 0$
(by definition of $T_x$ this is the only other possibility),
then $A_{1s} = 0$ for all $s \ge 2$;
we then replace the first column of $A$ by
the first column minus $x^{j_2 - j_1}$ times the second column,
so that the new first column has $x^{i_1-j_1}$ in its first entry
(or zero if $i_1 < j_1$) and zeroes elsewhere.
Then $\det A$ equals $x^{i_1-j_1}$ (or zero if $i_1 < j_1$)
times the determinant of its last $k-1$ rows and columns,
so the claim follows from the inductive hypothesis.
\qed

\noindent
See also Example~\ref{exam.toeplitz.power.TP} below
for a second proof of the total positivity of $T_x$,
using production matrices.

\begin{lemma}[Binomial matrix]
   \label{lemma.binomialmatrix.TP}
In the ring $\Z$, the binomial matrix
$B = \bigl( {\textstyle \binom{n}{k}} \bigr) _{n,k \ge 0}$
is totally positive.
More generally, the weighted binomial matrix
$B_{x,y}  = \bigl( x^{n-k} y^k {\textstyle \binom{n}{k}} \bigr) _{\! n,k \ge 0}$
is totally positive in the ring $\Z[x,y]$ equipped with the
coefficientwise order.
\end{lemma}

\proof
It is well known that the binomial matrix $B$ is totally positive,
and this can be proven by a variety of methods:
e.g.\ using production matrices
\cite[pp.~136--137, Example~6.1]{Karlin_68}
\cite[pp.~108--109]{Pinkus_10},
by diagonal similarity to a totally positive Toeplitz matrix
\cite[p.~109]{Pinkus_10},
by exponentiation of a nonnegative lower-subdiagonal matrix
\cite[p.~63]{Fallat_11},
or by an application of the Lindstr\"om--Gessel--Viennot lemma
\cite[p.~24]{Fomin_00}.

Then $B_{x,y} = D B D'$ where $D = \diag\bigl( (x^n)_{n \ge 0} \bigr)$
and $D' = \diag\bigl( (x^{-k} y^k)_{k \ge 0} \bigr)$.
By Cauchy--Binet, $B_{x,y}$ is totally positive in the ring $\Z[x,x^{-1},y]$
equipped with the coefficientwise order.
But because $B$ is lower-triangular,
the elements of $B_{x,y}$ actually lie in the subring $\Z[x,y]$.
\qed

\noindent
See also Example~\ref{exam.binomial.matrix.TP} below
for an {\em ab initio}\/ proof of Lemma~\ref{lemma.binomialmatrix.TP}
using production matrices.

Finally, let us show that the sufficiency half
of the Aissen--Schoenberg--Whitney--Edrei theorem holds
(with a slight modification to avoid infinite products)
in a general partially ordered commutative ring.
We give two versions, depending on whether or~not
it is assumed that the ring $R$ contains the rationals:

\begin{lemma}[Sufficient condition for Toeplitz-total positivity]
   \label{lemma.toeplitz.1}
Let $R$ be a partially ordered commutative ring,
let $N$ be a nonnegative integer,
and let $\alpha_1,\ldots,\alpha_N$, $\beta_1,\ldots,\beta_N$ and $C$
be nonnegative elements in $R$.
Define the sequence $\ba = (a_n)_{n \ge 0}$ in $R$ by
\be
   \sum_{n=0}^\infty a_n t^n
   \;=\;
   C \, \prod_{i=1}^N {1 + \alpha_i t  \over  1 - \beta_i t}
   \;.
 \label{eq.thm.toeplitz.1}
\ee
Then the Toeplitz matrix $T_\infty(\ba)$ is totally positive.
\end{lemma}

\noindent
Of course, it is no loss of generality to have the same number $N$
of alphas and betas,
since some of the $\alpha_i$ or $\beta_i$ could be zero.

\begin{lemma}[Sufficient condition for Toeplitz-total positivity, with rationals]
   \label{lemma.toeplitz.1a}
Let $R$ be a partially ordered commutative ring containing the rationals,
let $N$ be a nonnegative integer,
and let $\alpha_1,\ldots,\alpha_N$, $\beta_1,\ldots,\beta_N$, $\gamma$ and $C$
be nonnegative elements in $R$.
Define the sequence $\ba = (a_n)_{n \ge 0}$ in $R$ by
\be
   \sum_{n=0}^\infty a_n t^n
   \;=\;
   C \, e^{\gamma t} \, \prod_{i=1}^N {1 + \alpha_i t  \over  1 - \beta_i t}
   \;.
 \label{eq.thm.toeplitz.1a}
\ee
Then the Toeplitz matrix $T_\infty(\ba)$ is totally positive.
\end{lemma}

\proofof{Lemma~\ref{lemma.toeplitz.1}}
We make a series of elementary observations:

1) The sequence $\ba = (1,\alpha,0,0,0,\ldots)$, corresponding
to the generating function $A(t) = 1 + \alpha t$,
is Toeplitz-totally positive if and only if $\alpha \ge 0$.
The ``only if'' is trivial,
and the ``if'' follows from Lemma~\ref{lemma.bidiagonal}
because the Toeplitz matrix $T_\infty(\ba)$ is bidiagonal.

2) The sequence $\ba = (1,\beta,\beta^2,\beta^3,\ldots)$,
corresponding to the generating function $A(t) = 1/(1 - \beta t)$,
is Toeplitz-totally positive if and only if $\beta \ge 0$.
The ``only if'' is again trivial,
and the ``if'' follows from Lemma~\ref{lemma.toeplitz.power}.

3) If $\ba$ and $\bb$ are sequences
with ordinary generating functions $A(t)$ and $B(t)$,
then the convolution $\bc = \ba * \bb$,
defined by $c_n = \sum_{k=0}^n a_k b_{n-k}$,
has ordinary generating function $C(t) = A(t) \, B(t)$;
moreover, the Toeplitz matrix $T_\infty(\bc)$ is simply
the matrix product $T_\infty(\ba) \, T_\infty(\bb)$.
It thus follows from the Cauchy--Binet formula that if
$\ba$ and~$\bb$ are Toeplitz-totally positive, then so is $\bc$.

4) A Toeplitz-totally positive sequence can be multiplied
by a nonnegative constant $C$, and it is still Toeplitz-totally positive.

Combining these observations proves the lemma.
\qed

\proofof{Lemma~\ref{lemma.toeplitz.1a}}
We add to the proof of Lemma~\ref{lemma.toeplitz.1}
the following additional observation:

5) The sequence $\ba = (\gamma^n/n!)_{n \ge 0}$,
corresponding to the generating function $A(t) = e^{\gamma t}$,
is Toeplitz-totally positive if and only if $\gamma \ge 0$.
The ``only if'' is again trivial,
and the ``if'' follows from Lemma~\ref{lemma.binomialmatrix.TP}
because $\gamma^{n-k}/(n-k)! = \binom{n}{k} \gamma^{n-k} \times k!/n!$
and hence $T_\infty(\ba) = D^{-1} B_{\gamma,1} D$
where $D = \diag(\, (n!)_{n \ge 0})$.
\qed

\subsection{Production matrices}   \label{subsec.production}

The method of production matrices \cite{Deutsch_05,Deutsch_09}
has become in recent years an important tool in enumerative combinatorics.
In the special case of a tridiagonal production matrix,
this construction goes back to Stieltjes' \cite{Stieltjes_1889,Stieltjes_1894}
work on continued fractions:
the production matrix of a classical S-fraction or J-fraction is tridiagonal.
In~the present paper, by contrast,
we shall need production matrices that are lower-Hessenberg
(i.e.\ vanish above the first superdiagonal)
but are not in general tridiagonal.
We therefore begin by reviewing briefly
the basic theory of production matrices.
The important connection of production matrices with total positivity
will be treated in the next subsection.

Let $P = (p_{ij})_{i,j \ge 0}$ be an infinite matrix
with entries in a commutative ring $R$.
In~order that powers of $P$ be well-defined,
we shall assume that $P$ is either row-finite
(i.e.\ has only finitely many nonzero entries in each row)
or column-finite.

Let us now define an infinite matrix $A = (a_{nk})_{n,k \ge 0}$ by
\be
   a_{nk}  \;=\;  (P^n)_{0k}
 \label{def.iteration}
\ee
(in particular, $a_{0k} = \delta_{0k}$).
Writing out the matrix multiplications explicitly, we have
\be
   a_{nk}
   \;=\;
   \sum_{i_1,\ldots,i_{n-1}}
      p_{0 i_1} \, p_{i_1 i_2} \, p_{i_2 i_3} \,\cdots\,
        p_{i_{n-2} i_{n-1}} \, p_{i_{n-1} k}
   \;,
 \label{def.iteration.walk}
\ee
so that $a_{nk}$ is the total weight for all $n$-step walks in $\N$
from $i_0 = 0$ to $i_n = k$, in~which the weight of a walk is the
product of the weights of its steps, and a step from $i$ to $j$
gets a weight $p_{ij}$.
Yet another equivalent formulation is to define the entries $a_{nk}$
by the recurrence
\be
   a_{nk}  \;=\;  \sum_{i=0}^\infty a_{n-1,i} \, p_{ik}
   \qquad\hbox{for $n \ge 1$}
 \label{def.iteration.bis}
\ee
with the initial condition $a_{0k} = \delta_{0k}$.

We call $P$ the \textbfit{production matrix}
and $A$ the \textbfit{output matrix},
and we write $A = \scro(P)$.
Note that if $P$ is row-finite, then so is $\scro(P)$;
if $P$ is lower-Hessenberg, then $\scro(P)$ is lower-triangular;
if $P$ is lower-Hessenberg with invertible superdiagonal entries,
then $\scro(P)$ is lower-triangular with invertible diagonal entries;
and if $P$ is unit-lower-Hessenberg
(i.e.\ lower-Hessenberg with entries 1 on the superdiagonal),
then $\scro(P)$ is unit-lower-triangular.
In all the applications in this paper, $P$ will be lower-Hessenberg.

The matrix $P$ can also be interpreted as the adjacency matrix
for a weighted directed graph on the vertex set $\N$
(where the edge $ij$ is omitted whenever $p_{ij}  = 0$).
Then $P$ is row-finite (resp.\ column-finite)
if and only if every vertex has finite out-degree (resp.\ finite in-degree).

This iteration process can be given a compact matrix formulation.
Let us define the \textbfit{augmented production matrix}
\be
   \widetilde{P}
   \;\eqdef\;
   \left[
    \begin{array}{c}
         1 \;\; 0 \;\; 0 \;\; 0 \;\; \cdots \; \vphantom{\Sigma} \\
         \hline
         P
    \end{array}
    \right]
   \;.
 \label{def.prodmat_augmented}
\ee
Then the recurrence \reff{def.iteration.bis}
together with the initial condition $a_{0k} = \delta_{0k}$ can be written as
\be
   A
   \;=\;
   \left[
    \begin{array}{c}
         1 \;\; 0 \;\; 0 \;\; 0 \;\; \cdots \; \vphantom{\Sigma} \\
         \hline
         AP
    \end{array}
    \right]
   \;=\;
   \left[
    \begin{array}{c|c}
         1 & \bzero  \\
         \hline
         \bzero & A
    \end{array}
    \right]
   \left[
    \begin{array}{c}
         1 \;\; 0 \;\; 0 \;\; 0 \;\; \cdots \; \vphantom{\Sigma} \\
         \hline
         P
    \end{array}
    \right]
   \;=\;
   \left[
    \begin{array}{c|c}
         1 & \bzero  \\
         \hline
         \bzero & A
    \end{array}
    \right]
   \widetilde{P}
   \;.
 \label{eq.prodmat.u}
\ee
This identity can be iterated to give the factorization
\be
   A
   \;=\;
   \cdots\,
   \left[
    \begin{array}{c|c}
         I_3 & \bzero  \\
         \hline
                &   \\[-4mm]
         \bzero & \widetilde{P}
    \end{array}
    \right]
   \left[
    \begin{array}{c|c}
         I_2 & \bzero  \\
         \hline
                &   \\[-4mm]
         \bzero & \widetilde{P}
    \end{array}
    \right]
   \left[
    \begin{array}{c|c}
         I_1 & \bzero  \\
         \hline
                &   \\[-4mm]
         \bzero & \widetilde{P}
    \end{array}
    \right]
    \widetilde{P}
 \label{eq.prodmat.u.iterated}
\ee
where $I_k$ is the $k \times k$ identity matrix;
and conversely, \reff{eq.prodmat.u.iterated} implies \reff{eq.prodmat.u}.

Now let $\Delta = (\delta_{i+1,j})_{i,j \ge 0}$
be the matrix with 1 on the superdiagonal and 0 elsewhere.
Then for any matrix $M$ with rows indexed by $\N$,
the product $\Delta M$ is simply $M$ with its zeroth row removed
and all other rows shifted upwards.
(Some authors use the notation $\overline{M} \eqdef \Delta M$.)
The recurrence \reff{def.iteration.bis} can then be written as
\be
   \Delta \, \scro(P)  \;=\;  \scro(P) \, P
   \;.
 \label{def.iteration.bis.matrixform}
\ee
It follows that if $A$ is a row-finite matrix
that has a row-finite inverse $A^{-1}$
and has first row $a_{0k} = \delta_{0k}$,
then $P = A^{-1} \Delta A$ is the unique matrix such that $A = \scro(P)$.
This holds, in particular, if $A$ is lower-triangular with
invertible diagonal entries and $a_{00} = 1$;
then $A^{-1}$ is lower-triangular
and $P = A^{-1} \Delta A$ is lower-Hessenberg.
And if $A$ is unit-lower-triangular,
then $P = A^{-1} \Delta A$ is unit-lower-Hessenberg.

We shall repeatedly use the following easy fact:

\begin{lemma}[Production matrix of a product]
   \label{lemma.production.AB}
Let $P = (p_{ij})_{i,j \ge 0}$ be a row-finite matrix
(with entries in a commutative ring $R$),
with output matrix $A = \scro(P)$;
and let $B = (b_{ij})_{i,j \ge 0}$
be a lower-triangular matrix with invertible (in $R$) diagonal entries.
Then
\be
   AB \;=\;  b_{00} \, \scro(B^{-1} P B)
   \;.
\ee
That is, up to a factor $b_{00}$,
the matrix $AB$ has production matrix $B^{-1} P B$.
\end{lemma}

\proof
Since $P$ is row-finite, so is $A = \scro(P)$;
then the matrix products $AB$ and $B^{-1} P B$
arising in the lemma are well-defined.  Now
\be
   a_{nk}
   \;=\;
   \sum_{i_1,\ldots,i_{n-1}}
      p_{0 i_1} \, p_{i_1 i_2} \, p_{i_2 i_3} \,\cdots\,
        p_{i_{n-2} i_{n-1}} \, p_{i_{n-1} k}
   \;,
\ee
while
\be
   \scro(B^{-1} P B)_{nk}
   \;=\;
   \sum_{j,i_1,\ldots,i_{n-1},i_n}
      (B^{-1})_{0j} \,
      p_{j i_1} \, p_{i_1 i_2} \, p_{i_2 i_3} \,\cdots\,
        p_{i_{n-2} i_{n-1}} \, p_{i_{n-1} i_n} \, b_{i_n k}
   \;.
\ee
But $B$ is lower-triangular with invertible diagonal entries,
so $B$ is invertible and $B^{-1}$ is lower-triangular,
with $(B^{-1})_{0j} = b_{00}^{-1} \delta_{j0}$.
It follows that $AB = b_{00} \, \scro(B^{-1} P B)$.
\qed

\subsection{Production matrices and total positivity}
   \label{subsec.totalpos.prodmat}

Let $P = (p_{ij})_{i,j \ge 0}$ be a matrix with entries in a
partially ordered commutative ring $R$.
We will use $P$ as a production matrix;
let $A = \scro(P)$ be the corresponding output matrix.
As before, we assume that $P$ is either row-finite or column-finite.

When $P$ is totally positive, it turns out \cite{Sokal_totalpos}
that the output matrix $\scro(P)$ has {\em two}\/ total-positivity properties:
firstly, it is totally positive;
and secondly, its zeroth column is Hankel-totally positive.
Since \cite{Sokal_totalpos} is not yet publicly available,
we shall present briefly here (with proof) the main results
that will be needed in the sequel.

The fundamental fact that drives the whole theory is the following:

\begin{proposition}[Minors of the output matrix]
   \label{prop.iteration.homo}
Every $k \times k$ minor of the output matrix $A = \scro(P)$
can be written as a sum of products of minors of size $\le k$
of the production matrix $P$.
\end{proposition}

In this proposition the matrix elements $\bfp = \{p_{ij}\}_{i,j \ge 0}$
should be interpreted in the first instance as indeterminates:
for instance, we can fix a row-finite or column-finite set
$S \subseteq \N \times \N$
and define the matrix $P^S = (p^S_{ij})_{i,j \in \N}$ with entries
\be
   p^S_{ij}
   \;=\;
   \begin{cases}
       p_{ij}  & \textrm{if $(i,j) \in S$} \\[1mm]
       0       & \textrm{if $(i,j) \notin S$}
   \end{cases}
\ee
Then the entries (and hence also the minors) of both $P$ and $A$
belong to the polynomial ring $\Z[\bfp]$,
and the assertion of Proposition~\ref{prop.iteration.homo} makes sense.
Of course, we can subsequently specialize the indeterminates $\bfp$
to values in any commutative ring $R$.

\proofof{Proposition~\ref{prop.iteration.homo}}
%
For any infinite matrix $X = (x_{ij})_{i,j \ge 0}$,
let us write $X_N = (x_{ij})_{0 \le i \le N-1 ,\, j \ge 0}$
for the submatrix consisting of the first $N$ rows
(and {\em all}\/ the columns) of $X$.
Every $k \times k$ minor of $A$ is of course
a $k \times k$ minor of $A_N$ for some $N$,
so it suffices to prove that the claim about minors holds for all the $A_N$.
But this is easy: the fundamental identity \reff{eq.prodmat.u} implies
\be
   A_N
   \;=\;
   \left[
    \begin{array}{c|c}
         1 & \bzero  \\
         \hline
         \bzero & A_{N-1}
    \end{array}
    \right]
   \,
   \left[
    \begin{array}{c}
         1 \;\; 0 \;\; 0 \;\; 0 \;\; \cdots \; \vphantom{\Sigma} \\
         \hline
         P
    \end{array}
    \right]
   \;.
 \label{eq.proof.prop.iteration.homo}
\ee
So the result follows by induction on $N$, using the Cauchy--Binet formula.
\qed

If we now specialize the indeterminates $\bfp$
to values in some partially ordered commutative ring $R$,
we can immediately conclude:

\begin{theorem}[Total positivity of the output matrix]
   \label{thm.iteration.homo}
Let $P$ be an infinite matrix that~is either row-finite or column-finite,
with entries in a partially ordered commutative ring $R$.
If $P$ is totally positive of order~$r$, then so is $A = \scro(P)$.
\end{theorem}

\medskip

{\bf Remarks.}
1.  In the case $R = \R$, Theorem~\ref{thm.iteration.homo}
is due to Karlin \cite[pp.~132--134]{Karlin_68};
see also \cite[Theorem~1.11]{Pinkus_10}.
Karlin's proof is different from ours.

2.  Our quick inductive proof of Proposition~\ref{prop.iteration.homo}
follows an idea of Zhu \cite[proof of Theorem~2.1]{Zhu_13},
which was in turn inspired in part by Aigner \cite[pp.~45--46]{Aigner_99}.
The same idea recurs in recent work of several authors
\cite[Theorem~2.1]{Zhu_14}
\cite[Theorem~2.1(i)]{Chen_15a}
\cite[Theorem~2.3(i)]{Chen_15b}
\cite[Theorem~2.1]{Liang_16}
\cite[Theorems~2.1 and 2.3]{Chen_19}
\cite{Gao_non-triangular_transforms}.
However, all of these results concerned only special cases:
\cite{Aigner_99,Zhu_13,Chen_15b,Liang_16}
treated the case in~which the production matrix $P$ is tridiagonal;
\cite{Zhu_14} treated a (special) case in~which $P$ is upper bidiagonal;
\cite{Chen_15a} treated the case in~which
$P$ is the production matrix of a Riordan array;
\cite{Chen_19,Gao_non-triangular_transforms}
treated (implicitly) the case in~which $P$ is upper-triangular and Toeplitz.
But the argument is in fact completely general, as we have just seen;
there is no need to assume any special form for the matrix $P$.

3. A slightly different version of this proof
was presented in \cite{latpath_SRTR,latpath_lah}.
The simplified reformulation
given here,
using the augmented production matrix,
is due to Mu and Wang \cite{Mu_20}.
\myendremark

\medskip

\begin{example}[Toeplitz matrix of powers]
   \label{exam.toeplitz.power.TP}
\rm
Let $P = x {\bf e}_{00} + y \Delta$,
where $x$ and $y$ are indeterminates
(here ${\bf e}_{ij}$ denotes the matrix with an entry~1 in position~$ij$
 and 0 elsewhere).
By Lemma~\ref{lemma.bidiagonal}, $P$ is TP
in the ring $\Z[x,y]$ equipped with the coefficientwise order.
An easy computation shows that
$\scro(x {\bf e}_{00} + y\Delta)_{nk} = x^{n-k} y^k \, {\rm I}[k \le n]$.
When $y=1$, this is the Toeplitz matrix of powers \reff{def.Tx}.
So Theorem~\ref{thm.iteration.homo} implies that $T_x$ is TP
in the ring $\Z[x]$ equipped with the coefficientwise order.
This gives a second proof of the total positivity stated in
Lemma~\ref{lemma.toeplitz.power}.
\myendremark
\end{example}
\vspace*{-6mm}

\begin{example}[Binomial matrix]
   \label{exam.binomial.matrix.TP}
\rm
Let $P$ be the upper-bidiagonal Toeplitz matrix
$xI + y\Delta$, where $x$ and $y$ are indeterminates.
By Lemma~\ref{lemma.bidiagonal}, $P$ is TP
in the ring $\Z[x,y]$ equipped with the coefficientwise order.
An easy computation shows that $\scro(xI + y\Delta) = B_{x,y}$,
the weighted binomial matrix
with entries $(B_{x,y})_{nk} = x^{n-k} y^k \binom{n}{k}$.
So Theorem~\ref{thm.iteration.homo} implies that $B_{x,y}$ is TP
in the ring $\Z[x,y]$ equipped with the coefficientwise order.
This gives an {\em ab initio}\/ proof of Lemma~\ref{lemma.binomialmatrix.TP}.
\myendremark
\end{example}

\bigskip

Now define 
$\scroo_0(P)$ to be the zeroth-column sequence of $\scro(P)$, i.e.
\be
   \scroo_0(P)_n  \;\eqdef\;  \scro(P)_{n0}  \;\eqdef\;  (P^n)_{00}
   \;.
 \label{def.scroo0}
\ee
Then the Hankel matrix of $\scroo_0(P)$ has matrix elements
\begin{eqnarray}
   & &
   \!\!\!\!\!\!\!
   H_\infty(\scroo_0(P))_{nn'}
   \;=\;
   \scroo_0(P)_{n+n'}
   \;=\;
   (P^{n+n'})_{00}
   \;=\;
   \sum_{k=0}^\infty (P^n)_{0k} \, (P^{n'})_{k0}
   \;=\;
          \nonumber \\
   & &
   \sum_{k=0}^\infty (P^n)_{0k} \, ((P^{\rm T})^{n'})_{0k}
   \;=\;
   \sum_{k=0}^\infty \scro(P)_{nk} \, \scro(P^{\rm T})_{n'k}
   \;=\;
   \big[ \scro(P) \, {\scro(P^{\rm T})}^{\rm T} \big]_{nn'}
   \;.
   \qquad
\end{eqnarray}
(Note that the sum over $k$ has only finitely many nonzero terms:
 if $P$ is row-finite, then there are finitely many nonzero $(P^n)_{0k}$,
 while if $P$ is column-finite,
 there are finitely many nonzero $(P^{n'})_{k0}$.)
We have therefore proven:

\begin{lemma}[Identity for Hankel matrix of the zeroth column]
   \label{lemma.hankel.karlin}
Let $P$ be a row-finite or column-finite matrix
with entries in a commutative ring $R$.
Then
\be
   H_\infty(\scroo_0(P))
   \;=\;
   \scro(P) \, {\scro(P^{\rm T})}^{\rm T}
   \;.
\ee
\end{lemma}

{\bf Remark.}
If $P$ is row-finite, then $\scro(P)$ is row-finite;
$\scro(P^{\rm T})$ need not be row- or column-finite,
but the product $\scro(P) \, {\scro(P^{\rm T})}^{\rm T}$
is anyway well-defined.
Similarly, if $P$ is column-finite, then ${\scro(P^{\rm T})}^{\rm T}$
is column-finite;
$\scro(P)$ need not be row- or column-finite,
but the product $\scro(P) \, {\scro(P^{\rm T})}^{\rm T}$
is again well-defined.
\myendremark

\medskip

Combining Proposition~\ref{prop.iteration.homo}
with Lemma~\ref{lemma.hankel.karlin} and the Cauchy--Binet formula,
we obtain:

\begin{corollary}[Hankel minors of the zeroth column]
   \label{cor.iteration2}
Every $k \times k$ minor of the infinite Hankel matrix
$H_\infty(\scroo_0(P)) = ((P^{n+n'})_{00})_{n,n' \ge 0}$
can be written as a sum of products
of the minors of size $\le k$ of the production matrix $P$.
\end{corollary}

And specializing the indeterminates $\bfp$
to nonnegative elements in a partially ordered commutative ring,
in such a way that $P$ is row-finite or column-finite,
we deduce:

\begin{theorem}[Hankel-total positivity of the zeroth column]
   \label{thm.iteration2bis}
Let $P = (p_{ij})_{i,j \ge 0}$ be an infinite row-finite or column-finite
matrix with entries in a partially ordered commutative ring $R$,
and define the infinite Hankel matrix
$H_\infty(\scroo_0(P)) = ((P^{n+n'})_{00})_{n,n' \ge 0}$.
If $P$ is totally positive of order~$r$, then so is $H_\infty(\scroo_0(P))$.
\end{theorem}

One might hope that Theorem~\ref{thm.iteration2bis}
could be strengthened to show not only Hankel-TP of the zeroth column
of the output matrix $A = \scro(P)$,
but in fact Hankel-TP of the row-generating polynomials $A_n(x)$
for all $x \ge 0$ (at least when $R = \R$) ---
or even more strongly,
coefficientwise Hankel-TP of the row-generating polynomials.
Alas, this hope is vain, for these properties do not hold {\em in general}\/:

\begin{example}[Failure of Hankel-TP of the row-generating polynomials]
   \label{exam.hankel.karlin.rowgen}
\rm
Let $P = {\bf e}_{00} + \Delta$
be the upper-bidiagonal matrix with 1 on the superdiagonal
and $1,0,0,0,\ldots$ on the diagonal;
by Lemma~\ref{lemma.bidiagonal} it is TP.
Then $A = \scro(P)$ is the lower-triangular matrix will all entries 1
(see Example~\ref{exam.toeplitz.power.TP}),
so that $A_n(x) = \sum_{k=0}^n x^k$.
Since $A_0(x) \, A_2(x) - A_1(x)^2 = -x$,
the sequence $(A_n(x))_{n \ge 0}$ is not even log-convex
(i.e.\ Hankel-TP${}_2$) for any real number $x > 0$.
%
\myendremark
\end{example}

Nevertheless, in one important special case ---
namely, exponential Riordan arrays $\scrr[1,G]$ ---
the total positivity of the production matrix {\em does}\/
imply the coefficientwise Hankel-TP of the row-generating polynomials
of the output matrix:
this was shown \cite[Theorem~2.20]{forests_totalpos}.
That result will be generalized here, in Corollary~\ref{cor.EAZ.hankel},
to provide a more general {\em sufficient}\/ (but not necessary)
condition for the coefficientwise Hankel-TP of the row-generating polynomials
of the output matrix.

\subsection{Binomial row-generating matrices}  \label{subsec.rowgen}

Let $A = (a_{nk})_{n,k \ge 0}$ be a row-finite matrix
with entries in a commutative ring $R$.
(In most applications, including all those in the present paper,
 the matrix $A$ will be lower-triangular.)
We define its \textbfit{row-generating polynomials} in the usual way:
\be
   A_n(x)  \;\eqdef\;  \sum_{k=0}^\infty a_{nk} \, x^k
   \;,
 \label{def.An.0}
\ee
where the sum is actually finite because $A$ is row-finite.
More generally, let us define its
\textbfit{binomial partial row-generating polynomials}
\begin{subeqnarray}
   A_{n,k}(x)
   & \eqdef &
   \sum_{\ell=k}^\infty a_{n\ell} \, \binom{\ell}{k} \, x^{\ell-k}
         \\[2mm]
   & = &
   {1 \over k!} \, {d^k \over dx^k} \, A_n(x)
   \;.
 \label{def.Ank}
\end{subeqnarray}
(Note that the operator $(1/k!) \, d^k\!/\!dx^k$ has a well-defined action
 on the polynomial ring $R[x]$ even if $R$ does not contain the rationals,
 since $(1/k!) (d^k\!/\!dx^k) x^n = \binom{n}{k} x^{n-k}$.)
The polynomials $A_{n,k}(x)$ are the matrix elements of the
\textbfit{binomial row-generating matrix} $A B_x$:
\be
   (A B_x)_{nk}  \;=\;  A_{n,k}(x)  \;,
\ee
where $B_x = B_{x,1}$ is the weighted binomial matrix defined in \reff{def.Bx}.
The zeroth column of the matrix $A B_x$ consists
of the row-generating polynomials $A_n(x) = A_{n,0}(x)$.

In this paper the matrix $A$
will be either the matrix $\sfT = (t_{n,k})_{n,k \ge 0}$
or one of its polynomial generalizations.

We can now explain the method that we will use to prove
Theorems~\ref{thm1.2} and \ref{thm1.4}:

\begin{proposition}
   \label{prop.method}
Let $P$ be a row-finite matrix
with entries in a partially ordered commutative ring $R$,
and let $A = \scro(P)$.
\begin{itemize}
   \item[(a)]  If $P$ is totally positive of order~$r$, then so is $A$.
   \item[(b)]  If the matrix $B_x^{-1} P B_x$ is totally positive of order~$r$
in the ring $R[x]$ equipped with the coefficientwise order,
then the sequence $(A_n(x))_{n \ge 0}$ of row-generating polynomials
is Hankel-totally positive of order~$r$
in the ring $R[x]$ equipped with the coefficientwise order.
\end{itemize}
\end{proposition}

\noindent
Indeed, (a) is just a restatement of Theorem~\ref{thm.iteration.homo};
and (b) is an immediate consequence of Lemma~\ref{lemma.production.AB}
and Theorem~\ref{thm.iteration2bis}
together with the fact that the zeroth column of the matrix $A B_x$ consists
of the row-generating polynomials $A_n(x)$.

\subsection{Riordan arrays}   \label{subsec.riordan}

Let $R$ be a commutative ring,
and let $f(t) = \sum_{n=0}^\infty f_n t^n$
and $g(t) = \sum_{n=1}^\infty g_n t^n$ be formal power series
with coefficients in $R$; note that $g$ has zero constant term
(for~clarity we set $g_0 = 0$).
Then the (ordinary) \textbfit{Riordan array}
associated to the pair $(f,g)$
is the infinite lower-triangular matrix
$\scrr(f,g) = (\scrr(f,g)_{nk})_{n,k \ge 0}$ defined by
\be
   \scrr(f,g)_{nk}
   \;=\;
   [t^n] \, f(t) g(t)^k
   \;.
 \label{def.riordan}
\ee
That is, the $k$th column of $\scrr(f,g)$
has ordinary generating function $f(t) g(t)^k$.
Note that $\scrr(f,g)$ is invertible in the ring $R^{\N \times \N}_{\rm lt}$
of lower-triangular matrices
$\iff$
the diagonal elements $\scrr(f,g)_{nn} = f_0 g_1^n$
are invertible elements of the ring $R$
$\iff$
$f_0$ and $g_1$ are invertible elements of $R$
$\iff$
$f(t)$ has a multiplicative inverse $f(t)^{-1}$ in the ring $R[[t]]$
and $g(t)$ has a compositional inverse $\bar{g}(t)$ in the ring $R[[t]]$.

\medskip

\begin{quotation}
\small
{\bf Warning.}
We have interchanged the letters $f$ and $g$ compared to the notation
of Shapiro {\em et al.}\/ \cite{Shapiro_91,Shapiro_22}
and Barry \cite{Barry_16}.
This notation seems to us more natural, but the reader should be warned.
\end{quotation}

\medskip

We shall use an easy but important result that~is sometimes called
the \textbfit{fundamental theorem of Riordan arrays} (FTRA):

\begin{lemma}[Fundamental theorem of Riordan arrays]
   \label{lemma.FTRA}
Let $\bb = (b_n)_{n \ge 0}$ be a sequence with
ordinary generating function $B(t) = \sum_{n=0}^\infty b_n t^n$.
Considering $\bb$ as a column vector and letting $\scrr(f,g)$
act on it by matrix multiplication, we obtain a sequence $\scrr(f,g) \bb$
whose ordinary generating function is $f(t) \, B(g(t))$.
\end{lemma}

\proof
We compute
\begin{subeqnarray}
   \sum_{k=0}^n \scrr(f,g)_{nk} \, b_k
   & = &
   \sum_{k=0}^\infty [t^n] \, f(t) g(t)^k \, b_k
            \\[2mm]
   & = &
   [t^n] \: f(t) \sum_{k=0}^\infty b_k \, g(t)^k
            \\[2mm]
   & = &
   [t^n] \: f(t) \, B(g(t))
   \;.
\end{subeqnarray}
\qed

We can now determine the production matrix of a Riordan array $\scrr(f,g)$.
Let $\ba = (a_n)_{n \ge 0}$ and $\bz = (z_n)_{n \ge 0}$
be sequences in a commutative ring $R$,
with ordinary generating functions
$A(t) = \sum_{n=0}^\infty a_n t^n$
and $Z(t) = \sum_{n=0}^\infty z_n t^n$.
We then define the \textbfit{AZ matrix}
associated to the sequences $\ba$ and $\bz$ by
\be
   \AZ(\ba,\bz)_{ij}
   \;=\;
   \begin{cases}
       z_i        & \textrm{if $j=0$} \\[1mm]
       a_{i-j+1}  & \textrm{if $j \ge 1$}
   \end{cases}
 \label{def.AZ.1}
\ee
or in other words
\be
   \AZ(\ba,\bz)  \;=\;
   \begin{bmatrix}
      z_0 & a_0 &   0  &  0   &  0   &    \\
      z_1 & a_1 & a_0 &   0  &   0  &    \\
      z_2 & a_2 & a_1 & a_0 &   0  &    \\
      z_3 & a_3 & a_2 & a_1 & a_0 &    \\
      \vdots & \vdots  &    & \cdots &    & \ddots
   \end{bmatrix}
   \;.
 \label{eq.AZmatrix.riordan}
\ee
We also write $\AZ(A,Z)$ as a synonym for $\AZ(\ba,\bz)$.
It is convenient to define also
\be
   Y(t)  \;=\;  {A(t) \over A(t) - tZ(t)}
   \;,
 \label{def.riordan.Yt}
\ee
which is well-defined if $a_0$ is invertible in $R$.
We then have
\cite{Deutsch_09,He_15}
\cite[pp.~148--149]{Barry_16}
\cite[Theorems~4.15 and 6.29]{Shapiro_22}\footnote{
   This theorem is also essentially contained in
   \cite[Theorems~3.2, 3.6 and 3.7]{Merlini_00},
   though those authors do not use the terminology of production matrices.
}:

\begin{theorem}[Production matrices of Riordan arrays]
   \label{thm.riordan.production}
Let $L$ be a lower-triangular matrix (with entries in a commutative ring $R$)
with invertible diagonal entries and $L_{00} = 1$,
and let $P = L^{-1} \Delta L$ be its production matrix.
Then $L$ is a Riordan array if and only~if $P$ is an AZ-matrix.

More precisely, $L = \scrr(f,g)$ if and only~if $P = \AZ(\ba,\bz)$,
where the generating functions $\big( f(t), g(t) \big)$
and $\big( A(t), Z(t) \big)$ are connected by
\be
   g(t) \;=\; t \, A(g(t))  \;,\qquad
   f(t)  \;=\;  {1 \over 1 \,-\, t Z(g(t))}
         \;=\;  Y(g(t))
 \label{eq.thm.riordan.production}
\ee
or equivalently
\be
   A(t)  \;=\;  {t \over \bar{g}(t)}  \;,\qquad
   Z(t)  \;=\;  {f(\bar{g}(t)) \,-\, 1
                     \over  
                     \bar{g}(t) \: f(\bar{g}(t))
                    }
   \;.
 \label{eq.prop.riordan.production.2}
\ee
\end{theorem}

\par\bigskip\noindent{\sc Proof} \cite[p.~18]{He_15}.
Suppose that $L = \scrr(f,g)$.
The hypotheses on $L$ imply that $f_0 = 1$
and that $g_1$ is invertible in $R$;
so $g(t)$ has a compositional inverse $\bar{g}(t)$.
Now let $(p_k(t))_{k \ge 0}$ be the column generating functions
of $P = L^{-1} \Delta L$.
Applying the FTRA to each column of $P$,
we see that $\scrr(f,g) P$ is a matrix whose column generating functions
are $\big( f(t) \, p_k(g(t)) \big)_{k \ge 0}$.
On the other hand, $\Delta \, \scrr(f,g)$
is the matrix $\scrr(f,g)$ with its zeroth row removed
and all other rows shifted upwards,
so it has column generating functions $[f(t)-1]/t$ for column~0
and $f(t) g(t)^k/t$ for columns $k \ge 1$.
Comparing these two results, we see that
$\Delta \, \scrr(f,g) = \scrr(f,g) \, P$ if and only~if
\be
   f(t) \, p_0(g(t))  \;=\;  {f(t) - 1 \over t}
 \label{eq.proof.prop.riordan.production.1}
\ee
and
\be
   p_k(g(t))  \;=\; {g(t)^k \over t}  \qquad\hbox{for $k \ge 1$} \;.
\ee
The latter equation can be rewritten as
\be
   p_k(t)  \;=\;  {t^k \over \bar{g}(t)}
   \;,
\ee
which means that the columns $k \ge 1$ of the production matrix $P$
are identical with those of $\AZ(\ba,\bz)$,
when $\ba$ is given by \reff{eq.prop.riordan.production.2}.
And \reff{eq.proof.prop.riordan.production.1} then states that
column~0 of the production matrix $P$
is identical with that of $\AZ(\ba,\bz)$,
when $\bz$ is given by \reff{eq.prop.riordan.production.2}.
Therefore, $L = \scrr(f,g)$ implies that $L^{-1} \Delta L = \AZ(\ba,\bz)$
where $\ba$ and $\bz$ are given by \reff{eq.prop.riordan.production.2}.

Conversely, suppose that $P = \AZ(\ba,\bz)$.
Let $g(t)$ be the unique formal power series in $R[[t]]$
with $g(0) = 0$ that satisfies the functional equation $g(t) = t \, A(g(t))$,
and then let $f(t) = 1/[1 - tZ(g(t))]$.
Then running the foregoing computation backwards shows that
$\Delta \, \scrr(f,g) = \scrr(f,g) \, P$.
Since by hypothesis $L_{00} = 1$, it follows that $L = \scro(P) = \scrr(f,g)$.
\qed

\subsection{Exponential Riordan arrays}   \label{subsec.expriordan}

Let $R$ be a commutative ring containing the rationals,
and let $F(t) = \sum_{n=0}^\infty f_n t^n/n!$
and $G(t) = \sum_{n=1}^\infty g_n t^n/n!$ be formal power series
with coefficients in $R$; we set $g_0 = 0$.
Then the \textbfit{exponential Riordan array}
\cite{Deutsch_04,Deutsch_09,Barry_16,Shapiro_22}
associated to the pair $(F,G)$
is the infinite lower-triangular matrix
$\scrr[F,G] = (\scrr[F,G]_{nk})_{n,k \ge 0}$ defined by
\be
   \scrr[F,G]_{nk}
   \;=\;
   {n! \over k!} \:
   [t^n] \, F(t) G(t)^k
   \;.
 \label{def.RFG}
\ee
That is, the $k$th column of $\scrr[F,G]$
has exponential generating function $F(t) G(t)^k/k!$.
Equivalently, the bivariate exponential generating function of $\scrr[F,G]$ is
\be
   \sum_{n,k=0}^\infty \scrr[F,G]_{nk} \, {t^n \over n!} \, x^k
   \;=\;
   F(t) \, e^{x G(t)}
   \;.
 \label{eq.egf.RFG}
\ee
The diagonal elements of $\scrr[F,G]$ are $\scrr[F,G]_{nn} = f_0 g_1^n$,
so the matrix $\scrr[F,G]$ is invertible
in the ring $R^{\N \times \N}_{\rm lt}$ of lower-triangular matrices
if and only if $f_0$ and $g_1$ are invertible in $R$.

Please note that the exponential Riordan array $\scrr[F,G]$
is nothing other than a diagonal similarity transform
of the ordinary Riordan array $\scrr(F,G)$
associated to the same power series $F$ and $G$:  that~is,
\be
   \scrr[F,G]  \;=\; D \, \scrr(F,G) \, D^{-1}
\ee
where $D = \diag\bigl( (n!)_{n \ge 0} \bigr)$.


\begin{lemma}[Fundamental theorem of exponential Riordan arrays]
   \label{lemma.FETRA}
Let $\bb = (b_n)_{n \ge 0}$ be a sequence with
exponential generating function $B(t) = \sum_{n=0}^\infty b_n t^n/n!$.
Considering $\bb$ as a column vector and letting $\scrr[F,G]$
act on it by matrix multiplication, we obtain a sequence $\scrr[F,G] \bb$
whose exponential generating function is $F(t) \, B(G(t))$.
\end{lemma}

\proof
We compute
\begin{subeqnarray}
   \sum_{k=0}^n \scrr[F,G]_{nk} \, b_k
   & = &
   \sum_{k=0}^\infty {n! \over k!} \, [t^n] \, F(t) G(t)^k \, b_k
            \\[2mm]
   & = &
   n! \: [t^n] \: F(t) \sum_{k=0}^\infty b_k \, {G(t)^k \over k!}
            \\[2mm]
   & = &
   n! \: [t^n] \: F(t) \, B(G(t))
   \;.
\end{subeqnarray}
\qed

Let us now consider the product of two exponential Riordan arrays
$\scrr[F_1,G_1]$ and $\scrr[F_2,G_2]$.
Applying the FTERA to the $k$th column of $\scrr[F_2,G_2]$,
whose exponential generating function is $F_2(t) G_2(t)^k/k!$,
we readily obtain:

\begin{lemma}[Product of two exponential Riordan arrays]
   \label{lemma.exp_riordan.product}
We have
\be
   \scrr[F_1,G_1] \, \scrr[F_2,G_2]
   \;=\;
   \scrr[ (F_2 \circ G_1) F_1 ,\, G_2 \circ G_1]
   \;.
 \label{eq.exp_riordan.composition}
\ee
\end{lemma}

In particular, if we let $\scrr[F_2,G_2]$
be the weighted binomial matrix $B_\xi = \scrr[e^{\xi t}, t]$
defined by \reff{def.Bx}, we obtain:

\begin{corollary}[Binomial row-generating matrix of an exponential Riordan array]
   \label{cor.exp_riordan.Bx}
We have
\be
   \scrr[F,G] \, B_\xi
   \;=\;
   \scrr[e^{\xi G} F, G]
   \;.
 \label{eq.cor.exp_riordan.Bx}
\ee
\end{corollary}

Similarly, letting $\scrr[F_1,G_1]$
be the weighted binomial matrix $B_\xi$, we obtain:

\begin{corollary}[Left binomial transform of an exponential Riordan array]
   \label{cor.Bx.exp_riordan}
We have
\be
   B_\xi \, \scrr[F,G]
   \;=\;
   \scrr[e^{\xi t} F, G]
   \;.
 \label{eq.cor.Bx.exp_riordan}
\ee
\end{corollary}

We can now determine the production matrix of an exponential
Riordan array $\scrr[F,G]$.
Let $\ba = (a_n)_{n \ge 0}$ and $\bz = (z_n)_{n \ge 0}$
be sequences in a commutative ring $R$,
with ordinary generating functions
$A(s) = \sum_{n=0}^\infty a_n s^n$
and $Z(s) = \sum_{n=0}^\infty z_n s^n$.
We then define the \textbfit{exponential AZ matrix}
associated to the sequences $\ba$ and $\bz$ by
\be
   \EAZ(\ba,\bz)_{nk}
   \;=\;
   {n! \over k!} \: (z_{n-k} \,+\, k \, a_{n-k+1})
   \;,
 \label{def.EAZ.1}
\ee
or equivalently (if $R$ contains the rationals)
\be
   \EAZ(\ba,\bz)
   \;=\;
   D \, T_\infty(\bz) \, D^{-1} \:+\: D \, T_\infty(\ba) \, D^{-1} \, \Delta
 \label{def.EAZ.2}
\ee
where $D = \diag\bigl( (n!)_{n \ge 0} \bigr)$.
We also write $\EAZ(A,Z)$ as a synonym for $\EAZ(\ba,\bz)$.

\bigskip

{\bf Remark.}  We have the exponential generating functions
\be
   \sum_{n,k=0}^\infty \EAZ(\ba,\bz)_{nk} \: {s^n \over n!} \, u^k
   \;=\;
   e^{su} \, \bigl[ Z(s) \,+\, u A(s) \bigr]
\ee
and
\be
   \sum_{n,k=0}^\infty \EAZ(\ba,\bz)_{nk} \: {s^n \over n!} \: k! \, u^k
   \;=\;
   {Z(s) \over 1-su} \:+\: {u \, A(s) \over (1-su)^2}
   \;.
\ee
\myendremark

\medskip

\begin{theorem}[Production matrices of exponential Riordan arrays]
   \label{thm.riordan.exponential.production}
Let $L$ be a lower-triangular matrix
(with entries in a commutative ring $R$ containing the rationals)
with invertible diagonal entries and $L_{00} = 1$,
and let $P = L^{-1} \Delta L$ be its production matrix.
Then $L$ is an exponential Riordan array
if and only~if $P$ is an exponential AZ matrix.

More precisely, $L = \scrr[F,G]$ if and only~if $P = \EAZ(A,Z)$,
where the generating functions $\big( F(t), G(t) \big)$
and $\big( A(s), Z(s) \big)$ are connected by
\be
   G'(t) \;=\; A(G(t))  \;,\qquad
   {F'(t) \over F(t)} \;=\; Z(G(t))
 \label{eq.prop.riordan.exponential.production.1}
\ee
or equivalently
\be
   A(s)  \;=\;  G'(\bar{G}(s))  \;,\qquad
   Z(s)  \;=\;  {F'(\bar{G}(s)) \over F(\bar{G}(s))}
 \label{eq.prop.riordan.exponential.production.2}
\ee
where $\bar{G}(s)$ is the compositional inverse of $G(t)$.
\end{theorem}

\par\bigskip\noindent{\sc Proof}
(mostly contained in \cite[pp.~217--218]{Barry_16}).
Suppose that $L = \scrr[F,G]$.
The hypotheses on $L$ imply that $f_0 = 1$
and that $g_1$ is invertible in $R$;
so $G(t)$ has a compositional inverse $\bar{G}(s)$.
Now let $P = (p_{nk})_{n,k \ge 0}$ be a matrix;
its column exponential generating functions are, by definition,
$P_k(t) = \sum_{n=0}^\infty p_{nk} \, t^n/n!$.
Applying the FTERA to each column of $P$,
we see that $\scrr[F,G] P$ is a matrix
whose column exponential generating functions
are $\big( F(t) \, P_k(G(t)) \big)_{k \ge 0}$.
On~the other hand, $\Delta \, \scrr[F,G]$
is the matrix $\scrr[F,G]$ with its zeroth row removed
and all other rows shifted upwards,
so it has column exponential generating functions
\be
   {d \over dt} \, \big( F(t) \, G(t)^k/k! \big)
   \;=\;
   {1 \over k!} \: \Big[ F'(t) \, G(t)^k
                         \:+\: k \, F(t) \, G(t)^{k-1} \, G'(t) \Big]
   \;.
\ee
Comparing these two results, we see that
$\Delta \, \scrr[F,G] = \scrr[F,G] \, P$
if and only~if
\be
   P_k(G(t))
   \;=\;
   {1 \over k!} \:
   {F'(t) \, G(t)^k \:+\: k \, F(t) \, G(t)^{k-1} \, G'(t)
    \over
    F(t)}
   \;,
\ee
or in other words
\be
   P_k(t)
   \;=\;
   {1 \over k!}  \:
      \biggl[ {F'(\bar{G}(t)) \over F(\bar{G}(t))} \, t^k
              \:+\: k \, t^{k-1} \, G'(\bar{G}(t))
      \biggr]
   \;.
\ee
Therefore
\begin{subeqnarray}
   p_{nk}
   & = &
   {n! \over k!} \: [t^n] \,
      \biggl[ {F'(\bar{G}(t)) \over F(\bar{G}(t))} \, t^k
              \:+\: k \, t^{k-1} \, G'(\bar{G}(t))
      \biggr]
    \\[2mm]
   & = &
   {n! \over k!} \:
      \biggl[ [t^{n-k}] \: {F'(\bar{G}(t)) \over F(\bar{G}(t))}
              \:+\: k \, [t^{n-k+1}] \: G'(\bar{G}(t))
      \biggr]
    \\[2mm]
   & = &
   {n! \over k!} \: (z_{n-k} \,+\, k \, a_{n-k+1})
\end{subeqnarray}
where $\ba = (a_n)_{n \ge 0}$ and $\bz = (z_n)_{n \ge 0}$
are given by \reff{eq.prop.riordan.exponential.production.2}.

Conversely, suppose that $P = \EAZ(A,Z)$.
Define $F(t)$ and $G(t)$
as the unique solutions (in the formal-power-series ring $R[[t]]$)
of the differential equations \reff{eq.prop.riordan.exponential.production.1}
with initial conditions $F(0) = 1$ and $G(0) = 0$.
Then running the foregoing computation backwards
shows that $\Delta \, \scrr[F,G] = \scrr[F,G] \, P$.
Since by hypothesis $L_{00} = 1$, it follows that $L = \scrr[F,G]$.
\qed

We refer to
$A(s) = \sum_{n=0}^\infty a_n s^n$ and $Z(s) = \sum_{n=0}^\infty z_n s^n$
as the \textbfit{A-series} and \textbfit{Z-series}
associated to the exponential Riordan array $\scrr[F,G]$.

\medskip

{\bf Remark.}
The identity $A(s) = G'(\bar{G}(s))$
can equivalently be written as $A(s) = 1/(\bar{G})'(s)$.
This is useful in comparing our work with that of Zhu \cite{Zhu_21a,Zhu_22},
who uses the latter formulation.
\myendremark

\bigskip

Let us now show how to rewrite the production matrix \reff{def.EAZ.2}
in a new way, which will be useful in what follows.
Define
\be
   \Psi(s)  \;\eqdef\;  F(\bar{G}(s))
   \;,
 \label{def.Psi}
\ee
so that $F(t) = \Psi(G(t))$ and $\Psi(0) = F(0) = 1$.
Then a simple computation using
\reff{eq.prop.riordan.exponential.production.1}/\reff{eq.prop.riordan.exponential.production.2}
shows that
\be
   Z(s)
   \;=\;
   {\Psi'(s) \over \Psi(s)} \, A(s)
   \;.
 \label{eq.ZoverA.Psi}
\ee
And let us define $\Phi(s) \eqdef A(s)/\Psi(s)$.
Then the pair $(\Phi,\Psi)$ is related to the pair $(A,Z)$ by
\begin{subeqnarray}
   A(s)   & = &  \Phi(s) \: \Psi(s)
      \\[2mm]
   Z(s)   & = &  \Phi(s) \: \Psi'(s)
 \label{def.Phi.AZ.first}
\end{subeqnarray}
And conversely, given any pair $(A,Z)$ of formal power series
(over a commutative ring~$R$ containing the rationals)
such that $A(0)$ is invertible in $R$,
there is a unique pair $(\Phi,\Psi)$ satisfying \reff{def.Phi.AZ.first}
together with the normalization $\Psi(0) = 1$, namely
\begin{subeqnarray}
   \Psi(s)  & = &  \exp\biggl[ \int {Z(s) \over A(s)} \, ds \biggr]
       \\[2mm]
   \Phi(s)  & = &  A(s) \, \exp\biggl[ -\int {Z(s) \over A(s)} \, ds \biggr]
\end{subeqnarray}
[Here the integral of a formal power series is defined by
\be
   \int \biggl( \sum_{n=0}^\infty \alpha_n s^n \biggr) \, ds
   \;\eqdef\;
   \sum_{n=0}^\infty \alpha_n \: {s^{n+1} \over n+1}
   \;.
\ee
It is the unique formal power series with zero constant term
whose derivative is the given series.]
We refer to $\Phi(s)$ and $\Psi(s)$
as the \textbfit{$\bm{\Phi}$-series} and \textbfit{$\bm{\Psi}$-series}
associated to the exponential Riordan array $\scrr[F,G]$.

Rewriting the production matrix \reff{def.EAZ.2}
in terms of the pair $(\Phi,\Psi)$
provides a beautiful --- and as we shall see, very useful --- factorization.
For reasons that shall become clear shortly
(see Lemma~\ref{lemma.BxinvEAZBx} below),
it is convenient to study the more general quantity $\EAZ(A,Z+\xi A)$:

\begin{proposition}
   \label{prop.EAZ.Phi.first}
Let $R$ be a commutative ring containing the rationals,
let $\Phi(s) = \sum\limits_{n=0}^\infty \phi_n  s^n$
and $\Psi(s) = \sum\limits_{n=0}^\infty \psi_n  s^n$
be formal power series with coefficients in $R$,
and let $A(s)$ and $Z(s)$ be defined by \reff{def.Phi.AZ.first}.
Now let $\xi$ be any element of $R$ (or an indeterminate).
Then
\be
   \EAZ(A, Z + \xi A)
   \;=\;
   [ D \, T_\infty(\bphi) \, D^{-1} ] \: (\Delta \,+\, \xi I) \:
        [ D \, T_\infty(\bpsi) \, D^{-1} ]
\ee
where $D = \diag\bigl( (n!)_{n \ge 0} \bigr)$.
\end{proposition}

To prove Proposition~\ref{prop.EAZ.Phi.first}, we need a lemma.
Given a sequence $\bpsi = (\psi_n)_{n \ge 0}$ in $R$
with ordinary generating function
$\Psi(s) = \sum_{n=0}^\infty \psi_n  s^n$,
we define $\bpsi' = (\psi'_n)_{n \ge 0}$ by $\psi'_n = (n+1) \psi_{n+1}$,
so that $\Psi'(s) = \sum_{n=0}^\infty \psi'_n s^n$.
We then have:

\begin{lemma}
   \label{lemma.toep.first}
Let $\bpsi$ and $\bpsi'$ be as above,
and let $D = \diag\bigl( (n!)_{n \ge 0} \bigr)$.
Then
\be
   T_\infty(\bpsi') \:+\: T_\infty(\bpsi) \, D^{-1} \Delta D
     \;=\;  D^{-1} \Delta D \, T_\infty(\bpsi)
   \;.
 \label{eq.lemma.toep.first}
\ee
\end{lemma}

\proof
All three matrices in \reff{eq.lemma.toep.first} are lower-Hessenberg,
and their $(n,k)$ matrix elements are (for $0 \le k \le n+1$)
\be
  (n-k+1) \psi_{n-k+1} \:+\: k \psi_{n-(k-1)}  \;=\;  (n+1) \psi_{(n+1)-k} 
  \;.
\ee
\qed

{\bf Remarks.}
1.  The identity \reff{eq.lemma.toep.first} can also be written as
$[D^{-1} \Delta D ,\, T_\infty(\bpsi)] = T_\infty(\bpsi')$,
where $[A,B] \eqdef AB-BA$ is the matrix commutator.
Thus, $[D^{-1} \Delta D ,\:\cdot\:]$ is the ``differentiation operator''
for Toeplitz matrices.
Note that $D^{-1} \Delta D$ is the matrix with $1,2,3,\ldots\,$
on the superdiagonal and zeroes elsewhere.

2.  Lemma~\ref{lemma.toep.first} was found independently
by Ding, Mu and Zhu \cite[proof of Theorem~2.1]{Ding_22}.
\myendremark

\proofof{Proposition~\ref{prop.EAZ.Phi.first}}
{}From \reff{def.EAZ.2} we have
\be
   \EAZ(A,Z+ \xi A)
   \;=\;
   D \, T_\infty(\bz+\xi\ba) \, D^{-1}
            \:+\: D \, T_\infty(\ba) \, D^{-1} \, \Delta
       \;.
\ee
The definitions \reff{def.Phi.AZ.first} imply
\begin{subeqnarray}
   T_\infty(\ba)       & = &  T_\infty(\bphi) \, T_\infty(\bpsi)
       \\[2mm]
   T_\infty(\bz+\xi\ba)  & = &  T_\infty(\bphi) \, T_\infty(\bpsi')
                              \:+\: \xi \, T_\infty(\bphi) \, T_\infty(\bpsi)
\end{subeqnarray}
Hence
\begin{subeqnarray}
   & &
   \!\!\!\!\!\!
   \EAZ(A,Z+\xi A)
         \nonumber \\[2mm]
   & & \quad =\;
   D \, \big[ T_\infty(\bphi) \, T_\infty(\bpsi')
               \:+\: \xi  \, T_\infty(\bphi) \, T_\infty(\bpsi) \big] \, D^{-1}
      \:+\:  D \, T_\infty(\bphi) \, T_\infty(\bpsi)  \, D^{-1} \, \Delta
         \nonumber \\[-1mm] \\
   & & \quad =\;
   D \, T_\infty(\bphi) \:
           \big[ \xi T_\infty(\bpsi) \,+\, T_\infty(\bpsi')
                      \,+\, T_\infty(\bpsi) \, D^{-1} \Delta D \big] \: D^{-1}
         \\[2mm]
   & & \quad =\;
   D \, T_\infty(\bphi) \: \big[ \xi T_\infty(\bpsi) \,+\,
                          D^{-1} \Delta D \, T_\infty(\bpsi) \big] \: D^{-1} \,
         \\[2mm]
   & & \quad =\;
   [D \, T_\infty(\bphi) \, D^{-1}] \: (\Delta \,+\, \xi I) \:
           [D \, T_\infty(\bpsi) \, D^{-1}]
   \;,
\end{subeqnarray}
where the next-to-last step used Lemma~\ref{lemma.toep.first}.
%
%
\qed

As an immediate consequence of Proposition~\ref{prop.EAZ.Phi.first}, we have:
   
\begin{corollary}
   \label{cor.EAZ.Phi.first}
Fix $1 \le r \le \infty$.
Let $R$ be a partially ordered commutative ring containing the rationals,
and let $\bphi = (\phi_n)_{n \ge 0}$ and $\bpsi = (\psi_n)_{n \ge 0}$
be sequences in $R$ that are Toeplitz-totally positive of order~$r$.
Let $\xi$ be an indeterminate.
With the definitions \reff{def.Phi.AZ.first},
the matrix $\EAZ(A,Z +\xi A)$ is totally positive of order~$r$
in the ring $R[\xi]$ equipped with the coefficientwise order.
\end{corollary}

\proof
By Lemma~\ref{lemma.bidiagonal},
the matrix $\Delta + \xi I$ is totally positive (of order $\infty$)
in the ring $\Z[\xi]$ equipped with the coefficientwise order.
By hypothesis the matrices $T_\infty(\bphi)$ and $T_\infty(\bpsi)$
are totally positive of order~$r$ in the ring $R$;
so Lemma~\ref{lemma.diagmult.TP} implies that
also $D \, T_\infty(\bphi) \, D^{-1}$ and $D \, T_\infty(\bpsi) \, D^{-1}$
are totally positive of order~$r$ in $R$.
The result then follows from Proposition~\ref{prop.EAZ.Phi.first}
and the Cauchy--Binet formula.
\qed

{\bf Remark.}
The hypothesis that the ring $R$ contains the rationals can be removed,
by using Lemma~\ref{lemma.diagmult.TP}
(see Section~\ref{subsec.diagonal_scaling})
together with the reasoning used in the proof of Theorem~\ref{thm1.7}
(see Section~\ref{subsec.Typhi}).
\myendremark

\medskip


It is worth observing that the converse to Corollary~\ref{cor.EAZ.Phi.first}
is false:

\begin{example}
   \label{exam.converse.cor.EAZ.Phi.first}
\rm
Let $A(s) = 1+s$ and $Z(s) = (\lambda + \mu) + \mu s$.
Then $P = \EAZ(A,Z)$ is the tridiagonal matrix with
$p_{n,n+1} = 1$, $p_{n,n} = \lambda + \mu + n$ and $p_{n,n-1} = n\mu$,
which can be written in the form $P = LU + \lambda I$,
where $L$ is the lower-bidiagonal matrix
with 1 on the diagonal and $1,2,3,\ldots$ on the subdiagonal,
$U$ is the upper-bidiagonal matrix
with 1 on the superdiagonal and $\mu$ on the diagonal,
and $I$ is the identity matrix;
so by the tridiagonal comparison theorem
\cite{Sokal_totalpos} \cite[Proposition~3.1]{Zhu_21b}
$P$ is totally positive, coefficientwise in $\lambda$ and $\mu$.
[When $\mu=0$ the total positivity is even more elementary,
 by Lemma~\ref{lemma.bidiagonal}.]
Note also that, in this example, $\EAZ(A,Z + \xi A)$
is simply $P = \EAZ(A,Z)$ with $\mu$ replaced by $\mu +\xi$.

But this pair $(A,Z)$ corresponds to
\be
   \Phi(s) \:=\: e^{-\mu s} \, (1 + s)^{1 - \lambda}
   \,\qquad
   \Psi(s) \:=\: e^{\mu s} \, (1 + s)^\lambda
\ee
which are {\em not}\/ Toeplitz-TP coefficientwise in $\lambda$ and $\mu$.
Indeed, even for real $\lambda$ and $\mu$,
the sequence $\bphi$ (resp.~$\bpsi$) is Toeplitz-TP only for
$\lambda \in \{0,1\}$ and $\mu \le 0$
(resp.~$\lambda \in \{0,1\}$ and $\mu \ge 0$).
So all nonnegative $\lambda,\mu$
other than $\lambda \in \{0,1\}$ and $\mu = 0$
yield counterexamples to the converse to Corollary~\ref{cor.EAZ.Phi.first},
and even to its restriction to the case $\xi = 0$.

In this example $F(t) = e^{\lambda t + \mu(e^t - 1)}$ and $G(t) = e^t - 1$,
so the exponential Riordan array is
$\scrr[F,G] = \scrr[e^{\lambda t + \mu (e^t -1)}, e^t - 1]
 = B_\lambda \scrr[1, e^t - 1] B_\mu$
by Corollaries~\ref{cor.exp_riordan.Bx} and \ref{cor.Bx.exp_riordan};
here $\scrr[1, e^t - 1]$ is the Stirling subset matrix \cite[A008277]{OEIS}.
\myendremark
\end{example}

\noindent
So the condition of Corollary~\ref{cor.EAZ.Phi.first}
is sufficient but not necessary for its conclusion.


\bigskip

%
%


Finally, a central role will be played in this paper
by a simple but remarkable identity for $B_\xi^{-1} \EAZ(\ba,\bz) \, B_\xi$,
where $B_\xi$ is the $\xi$-binomial matrix defined in \reff{def.Bx}
and $\EAZ(\ba,\bz)$ is the exponential AZ matrix
defined in \reff{def.EAZ.1}/\reff{def.EAZ.2}.

\begin{lemma}[Identity for $B_\xi^{-1} \EAZ(\ba,\bz) \, B_\xi$]
   \label{lemma.BxinvEAZBx}
Let $\ba = (a_n)_{n \ge 0}$, $\bz = (z_n)_{n \ge 0}$ and $\xi$
be indeterminates.
Then
\be
   B_\xi^{-1} \EAZ(\ba,\bz) \, B_\xi
   \;=\;
   \EAZ(\ba,\bz+\xi\ba)
   \;.
 \label{eq.lemma.BxinvEAZBx}
\ee
\end{lemma}

The special case $\bz = \bzero$ of this lemma was proven in
\cite[Lemma~3.6]{latpath_lah};
a simpler proof was given in \cite[Lemma~2.16]{forests_totalpos}.
Here we give the easy generalization to include $\bz$.
We will give two proofs:
a first proof by direct computation
from the definition \reff{def.EAZ.1}/\reff{def.EAZ.2},
and a second proof using exponential Riordan arrays.

\firstproof
We use the matrix definition \reff{def.EAZ.2}:
\be
   \EAZ(\ba,\bz)
   \;=\;
   D \, T_\infty(\bz) \, D^{-1} \:+\: D \, T_\infty(\ba) \, D^{-1} \, \Delta
\ee
where $D = \diag\bigl( (n!)_{n \ge 0} \bigr)$.
Since $\EAZ(\ba,\bz) = \EAZ(\ba,\bzero) + \EAZ(\bzero,\bz)$,
it suffices to consider separately the two contributions.

The key observation is that
$B_\xi = D \, T_\infty\big( (\xi^n/n!)_{n \ge 0} \big) \, D^{-1}$.
Now two Toeplitz matrices always commute:
$T_\infty(\ba) \, T_\infty(\bb) =
 T_\infty(\ba \star \bb) =
 T_\infty(\bb) \, T_\infty(\ba)$.
It follows that $D T_\infty(\bz) D^{-1}$ and $D T_\infty(\ba) D^{-1}$
commute with $B_\xi$.
Therefore
\be
   B_\xi^{-1} \EAZ(\bzero,\bz) \, B_\xi
   \;=\;
   \EAZ(\bzero,\bz)
   \;.
 \label{eq.EAZ.bz}
\ee
On the other hand, the classic recurrence for binomial coefficients implies
\be
   \Delta B_\xi  \;=\; B_\xi \, (\xi I + \Delta)
\ee
(cf.\ Example~\ref{exam.binomial.matrix.TP}).  Therefore
\begin{subeqnarray}
   B_\xi^{-1} \EAZ(\ba,\bzero) B_\xi
   & = &
   B_\xi^{-1} \, D T_\infty(\ba) D^{-1} \, \Delta \, B_\xi
      \\[2mm]
   & = &
   B_\xi^{-1} \, D T_\infty(\ba) D^{-1} \, B_\xi \, (\xi I + \Delta)
      \\[2mm]
   & = &
   D T_\infty(\ba) D^{-1} \, (\xi I + \Delta)
      \\[2mm]
   & = &
   \EAZ(\ba,\xi\ba)
   \;.
 \label{eq.EAZ.ba}
\end{subeqnarray}
Adding \reff{eq.EAZ.bz} and \reff{eq.EAZ.ba} yields \reff{eq.lemma.BxinvEAZBx}.
\qed

\secondproof
Let $\scrr[F,G]$ be the exponential Riordan array
generated by the production matrix $\EAZ(A,Z)$
according to Theorem~\ref{thm.riordan.exponential.production},
so that
\be
   G'(t) \;=\; A(G(t))  \;,\qquad
   {F'(t) \over F(t)} \;=\; Z(G(t))
   \;.
\ee
By Corollary~\ref{cor.exp_riordan.Bx} we have
\be
   \scrr[F,G] \, B_\xi
   \;=\;
   \scrr[e^{\xi G} F, G]
   \;.
\ee
Since
\be
   {d \over dt} \, \log[e^{\xi G(t)} F(t)]
   \;=\;
   {F'(t) \over F(t)} \:+\: \xi G'(t)
   \;,
\ee
Theorem~\ref{thm.riordan.exponential.production} implies that
the exponential Riordan array $\scrr[e^{\xi G} F, G]$
has production matrix $\EAZ(\widehat{A},\widehat{Z})$ where
\be
   \widehat{A}(s,\xi) \;=\; A(s) \;,\qquad
   \widehat{Z}(s,\xi) \;=\; Z(s) \,+\, \xi A(s)
   \;.
\ee
On the other hand, by Lemma~\ref{lemma.production.AB}
the production matrix of $\scrr[e^{\xi G} F, G] = \scrr[F,G] \, B_\xi$
is $B_\xi^{-1} \EAZ(A,Z) \, B_\xi$.
\qed

{\bf Remark.}
A special case of the ideas in the second proof
can be found in \cite[Proposition~4]{Barry_17}.
\myendremark

Combining Proposition~\ref{prop.method}
with Corollary~\ref{cor.EAZ.Phi.first} and Lemma~\ref{lemma.BxinvEAZBx},
we obtain:

\begin{corollary}
   \label{cor.EAZ.hankel}
Fix $1 \le r \le \infty$.
Let $R$ be a partially ordered commutative ring containing the rationals,
and let $\bphi = (\phi_n)_{n \ge 0}$ and $\bpsi = (\psi_n)_{n \ge 0}$
be sequences in $R$ that are Toeplitz-totally positive of order~$r$.
Then the exponential Riordan array $\scrr[F,G]$ defined by
\reff{eq.prop.riordan.exponential.production.1}/\reff{def.Phi.AZ.first}
has the following two properties:
\begin{itemize}
      \item[(a)]  The lower-triangular matrix $\scrr[F,G]$
      is totally positive of order~$r$.
   \item[(b)]  The sequence of row-generating polynomials of $\scrr[F,G]$
      is coefficientwise Hankel-totally positive of order~$r$.
\end{itemize}
\end{corollary}

Corollary~\ref{cor.EAZ.hankel} will be the main theoretical tool in this paper.

\bigskip

{\bf Remarks.}
1. The special case $\Psi = 1$ (i.e., $F=1$)
of Corollary~\ref{cor.EAZ.hankel}(b)
was proven in \cite[Theorem~2.20]{forests_totalpos}
and was the fundamental theoretical tool of that paper.

2.  Another special case of Proposition~\ref{prop.EAZ.Phi.first}
and Corollaries~\ref{cor.EAZ.Phi.first} and \ref{cor.EAZ.hankel}
was employed recently by Ding, Mu and Zhu
\cite[proof of Theorem~2.1]{Ding_22}
to study some far-reaching generalizations of the Eulerian polynomials.

3. Example~\ref{exam.converse.cor.EAZ.Phi.first} shows that
the condition of Corollary~\ref{cor.EAZ.hankel} is sufficient
but not necessary for its two conclusions.
\myendremark

\bigskip

Finally, it is worth singling out a subclass of Riordan arrays
that will occur in the cases to be studied in the present paper:

\begin{lemma}
   \label{lemma.EAZ.zeroth_first}
Consider an exponential Riordan array $\scrr[F,G]$ with $F(0) = 1$
and corresponding series $A(s)$, $Z(s)$, $\Phi(s)$, $\Psi(s)$.
Then, for any constant $c$, the following are equivalent:
\begin{itemize}
   \item[(a)] $\scrr[F,G]_{n,0} = c \, \scrr[F,G]_{n,1}$ for all $n \ge 1$.
   \item[(b)] $\EAZ(A,Z)_{n,0} = c \EAZ(A,Z)_{n,1}$ for all $n \ge 0$.
   \item[(b${}'$)] $\EAZ(A,Z) = \EAZ(A,Z) \, \Delta^{\rm T} \,
                                           (c\, {\bf e}_{00} + \Delta)$
where ${\bf e}_{00}$ denotes the matrix with an entry~$1$ in position $(0,0)$
and all other entries zero.
   \item[(c)] $\Psi(s) = 1/(1-cs)$.
\end{itemize}
\end{lemma}

\proof
(a)$\iff$(c):
(a) holds if and only if $F(t) = 1 + c F(t) G(t)$,
or in other words $F(t) = 1 / [1 - cG(t)]$,
or in other words $\Psi(s) = 1/(1-cs)$.

(b)$\iff$(c):
By \reff{def.EAZ.1}, (b) holds
if and only if $z_n = c (z_{n-1} + a_n)$,
or in other words $Z(s) = c [sZ(s) + A(s)]$,
or in other words
\be
   {\Psi'(s) \over \Psi(s)}
   \;=\;
   {Z(s) \over A(s)}
   \;=\;
   {c \over 1-cs}
   \;.
\ee
Since $\Psi(0) = 1$, this is equivalent to $\Psi(s) = 1/(1-cs)$.

(b${}'$)$\implies$(b):
The zeroth column of the matrix $c\, {\bf e}_{00} + \Delta$
equals $c$ times its first column;
so for any matrix $M$,
the zeroth column of the matrix $M \, (c\, {\bf e}_{00} + \Delta)$
equals $c$ times its first column.

(b)$\implies$(b${}'$):
The matrix $\EAZ(A,Z) \, \Delta^{\rm T}$ is obtained from $\EAZ(A,Z)$
by removing its zeroth column; it is lower-triangular.
And since, by hypothesis, the zeroth column of $\EAZ(A,Z)$
is $c$ times its first column,
$\EAZ(A,Z)$ can be recovered from $\EAZ(A,Z) \, \Delta^{\rm T}$ by
by right-multiplying by $c\, {\bf e}_{00} + \Delta$.
\qed

The case $c=0$ (that~is, $\Psi = 1$ and hence $F=1$)
corresponds to the {\em associated subgroup}\/ (or {\em Lagrange subgroup}\/)
of exponential Riordan arrays;
it arose in our earlier work \cite{latpath_lah,forests_totalpos}
on generic Lah and rooted-forest polynomials.
Using criterion~(a), we can already see that
the matrix $\sfT$ defined in \reff{def.tnk} will correspond to $c=1$,
while the matrices $\sfT(y,z)$ and $\sfT(y,\bphi)$
defined in \reff{def.tnkyz}/\reff{def.tnkyphi} will correspond,
according to Propositions~\ref{prop.tn01yz} and \ref{prop.tn01yphi}, to $c=y$.
Of course, in order to apply Lemma~\ref{lemma.EAZ.zeroth_first}
we will first need to prove that these matrices are indeed
exponential Riordan arrays:  that will be done in Section~\ref{sec.EGF}.
But we can see now that, once we do this,
the $\Psi$-series will be $\Psi(s) = 1/(1-cs)$.

\subsection{A lemma on diagonal scaling}   \label{subsec.diagonal_scaling}

Given a lower-triangular matrix $A = (a_{nk})_{n,k \ge 0}$ 
with entries in a commutative ring $R$,
let us define the matrix $A^\sharp = (a^\sharp_{nk})_{n,k \ge 0}$ by
\be
   a^\sharp_{nk}  \;=\;  {n! \over k!} \: a_{nk}
   \;;
 \label{def.sharp}
\ee
this is well-defined since $a_{nk} \neq 0$ only when $n \ge k$,
in~which case $n!/k!$ is an integer.

If $R$ contains the rationals, we can of course write
$A^\sharp = D A D^{-1}$ where $D = \diag\big( (n!)_{n \ge 0} \big)$.
And if $R$ is a partially ordered commutative ring
that contains the rationals and $A$ is TP${}_r$,
then we deduce immediately from $A^\sharp = D A D^{-1}$
that also $A^\sharp$ is TP${}_r$.
The following simple lemma \cite[Lemma~3.7]{latpath_lah}
shows that this conclusion holds even when $R$ does not contain the rationals:

\begin{lemma}
   \label{lemma.diagmult.TP}
Let $A = (a_{ij})_{i,j \ge 0}$ be a lower-triangular matrix
with entries in a partially ordered commutative ring $R$,
and let $\bd = (d_i)_{i \ge 1}$.
Define the lower-triangular matrix
$A^{\sharp\bd} = (a^{\sharp\bd}_{ij})_{i,j \ge 0}$ by
\be
   a^{\sharp\bd}_{ij}  \;=\;  d_{j+1} d_{j+2} \cdots d_i \, a_{ij}
   \;.
\ee
Then:
\begin{itemize}
   \item[(a)] If $A$ is TP${}_r$ and $\bd$ are indeterminates,
      then $A^{\sharp\bd}$ is TP${}_r$ in the ring $R[\bd]$ equipped with
      the coefficientwise order.
   \item[(b)] If $A$ is TP${}_r$ and $\bd$ are nonnegative elements of $R$,
      then $A^{\sharp\bd}$ is TP${}_r$ in the ring $R$.
\end{itemize}
\end{lemma}

\proof
(a) Let $\bd = (d_i)_{i \ge 1}$ be commuting indeterminates,
and let us work in the ring $R[\bd,\bd^{-1}]$
equipped with the coefficientwise order.
Let $D = \diag(1,\, d_1,\, d_1 d_2,\, \ldots)$.
Then $D$ is invertible, and both $D$ and
$D^{-1} = \diag(1,\, d_1^{-1},\, d_1^{-1} d_2^{-1},\, \ldots)$
have nonnegative elements.
It follows that $A^{\sharp\bd} = D A D^{-1}$ is TP${}_r$
in the ring $R[\bd,\bd^{-1}]$ equipped with the coefficientwise order.
But the matrix elements $a^{\sharp\bd}_{ij}$
actually belong to the subring $R[\bd] \subseteq R[\bd,\bd^{-1}]$.
So $A^{\sharp\bd}$ is TP${}_r$ in the ring $R[\bd]$
equipped with the coefficientwise order.

(b) follows from (a) by specializing indeterminates.
\qed

\noindent
The special case $A^{\sharp\bd} = A^\sharp$ corresponds to taking $d_i = i$.

Lemma~\ref{lemma.diagmult.TP} will be important
to proving Theorem~\ref{thm1.7} in the case where the ring $R$
does not contain the rationals (see Section~\ref{subsec.Typhi}).

\subsection{Lagrange inversion}

We will use Lagrange inversion in the following form \cite{Gessel_16}:
If $\phi(u)$ is a formal power series
with coefficients in a commutative ring $R$ containing the rationals,
then there exists a unique formal power series $f(t)$
with zero constant term satisfying
\be
   f(t)  \;=\;  t \, \phi(f(t))
   \;,
\ee
and it is given by
\be
   [t^n] \, f(t)  \;=\;  {1 \over n} \, [u^{n-1}] \, \phi(u)^n
     \quad\hbox{for $n \ge 1$}
   \;;
\ee
and more generally, if $H(u)$ is any formal power series, then
\be
   [t^n] \, H(f(t))  \;=\;  {1 \over n} \, [u^{n-1}] \, H'(u) \, \phi(u)^n
     \quad\hbox{for $n \ge 1$}
   \;.
 \label{eq.lagrange.H}
\ee
In particular, taking $H(u) = u^k$ with integer $k \ge 0$, we have
\be
   [t^n] \, f(t)^k  \;=\;  {k \over n} \, [u^{n-k}] \, \phi(u)^n
     \quad\hbox{for $n \ge 1$}
   \;.
 \label{eq.lagrange.k}
\ee

\section{Bijective proofs}
   \label{sec.bijective}

In this section we give bijective proofs of Propositions~\ref{prop.tn01yz},
\ref{prop.equivalence}, \ref{prop.tn01yphi} and \ref{prop.equivalence.phi}.
This section can be skipped on a first reading,
as it is not needed for proving the main theorems of the paper.

\subsection{Proof of Propositions~\ref{prop.tn01yz} and \ref{prop.tn01yphi}}

Here we will prove Proposition~\ref{prop.tn01yz},
which asserts that the polynomials $t_{n,k}(y,z)$
defined in \reff{def.tnkyz}
satisfy $t_{n,0}(y,z) = y \, t_{n,1}(y,z)$ for all $n \ge 1$;
and more generally Proposition~\ref{prop.tn01yphi},
which asserts that the polynomials $t_{n,k}(y,\bphi)$
defined in \reff{def.tnkyphi}
satisfy $t_{n,0}(y,\bphi) = y \, t_{n,1}(y,\bphi)$ for all $n \ge 1$.

\begin{figure}[t]
\begin{center}
\begin{tikzpicture}[baseline=(current bounding box.base)]
	\node at (0,1) {$\bm{T \; (k=1)}$};
	\fill (0,0) circle(.1) node[above] {$r$};
	\fill (-1.5,-2) circle(.1) node[left] {$1$};
	\fill (-1.5,-3) circle(.1);
        \node at (-1.7,-2.85) {$a$};
	\foreach \x in {-1,0,1}
		{
		\draw (\x,-1) circle(.1);
		\draw (0,0)--(\x,-1);
		\draw[dashed] (\x,-1)--(\x+.5,-2);
		}
	\foreach \x in {-.5,.5}
		\node at (\x,-1) {$\cdots$};
	\draw (0,-1.2) ellipse(2 and 1.2);
	\node at (0,-2) {$T_r$};
	\draw (-1,-1)--(-1.5,-2);
	\draw [blue, ultra thick] (-1.5,-2)--(-1.5,-3);
        \node at (-2.2,-2.5) {$\textsl{proper}$};
	\foreach \x in {-2,-1}
		{
		\draw (\x,-4) circle(.1);
		\draw (-1.5,-3)--(\x,-4);
		\draw[dashed] (\x,-4)--(\x,-5);
		}
	\node at (-1.5,-4) {$\cdots$};
	\draw (-1.5,-4.2) ellipse(1 and 1.2);
	\node at (-1.5,-5) {$T_a$};
\end{tikzpicture}
\quad \raisebox{-70pt}{$\longleftrightarrow$} \quad
\begin{tikzpicture}[baseline=(current bounding box.base)]
	\node at (1.8,1) {$\bm{T' \; (k=0)}$};
	\fill (0,0) circle(.1) node[above] {$a$};
	\fill (3,-1) circle(.1) node[above] {$r$};
	\fill (1.5,-3) circle(.1) node[left] {$1$};
	\foreach \x in {-.5,.5}
		{
		\draw (\x,-1) circle(.1);
		\draw (0,0)--(\x,-1);
		\draw[dashed] (\x,-1)--(\x,-2);
		}
	\node at (0,-1) {$\cdots$};
	\draw (0,-1.2) ellipse(1 and 1.2);
	\node at (0,-2) {$T_a$};
	\draw [violet, ultra thick] (0,0)--(3,-1);
        \node at (2.2,-0.3) {$\textsl{improper}$};
	\foreach \x in {2,3,4}
		{
		\draw (\x,-2) circle(.1);
		\draw (3,-1)--(\x,-2);
		\draw[dashed] (\x,-2)--(\x+.5,-3);
		}
	\foreach \x in {2.5,3.5}
		\node at (\x,-2) {$\cdots$};
	\draw (3,-2.2) ellipse(2 and 1.2);
	\node at (3,-3) {$T_r$};
	\draw (2,-2)--(1.5,-3);
\end{tikzpicture}
\caption{Bijection between $T$ and $T'$.}
   \label{fig.prop.tn01yz}
\end{center}
\end{figure}
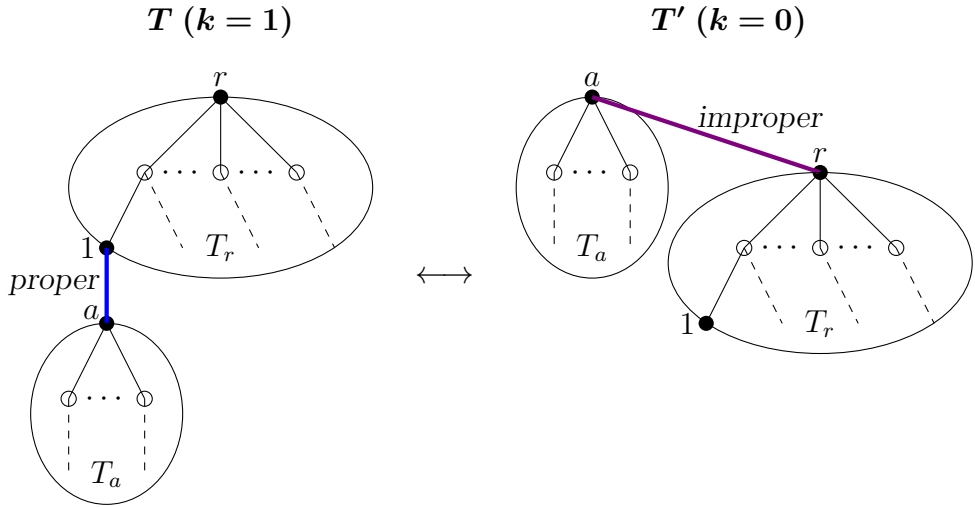

We will prove these results by constructing, for each $n \ge 1$,
a bijection from the set $\scrt^{\<1;1\>}_{n+1}$ 
of rooted trees on the vertex set $[n+1]$
in~which the vertex~1 has exactly one child,
to the set $\scrt^{\<1;0\>}_{n+1}$ of rooted trees on the vertex set $[n+1]$
in~which vertex~1 is a leaf,
with the properties that
\begin{itemize}
   \item[(a)] the number of improper edges is increased by 1, and
   \item[(b)] for each $m$, the number of vertices with $m$ proper children
        is preserved, provided that in $T \in \scrt^{\<1;1\>}_{n+1}$
        one ignores the vertex~1 (which has one child).
\end{itemize}
This construction is illustrated in Figure~\ref{fig.prop.tn01yz}.
Since the weight in \reff{def.tnkyphi} is $y$ for each improper edge
and $\phihat_m = m! \, \phi_m$ for each vertex $i \neq 1$
with $m$ proper children,
this proves $t_{n,0}(y,\bphi) = y \, t_{n,1}(y,\bphi)$.
Specializing to $\phi_m = z^m/m!$ then yields
$t_{n,0}(y,z) = y \, t_{n,1}(y,z)$.

\proofof{Proposition~\ref{prop.tn01yphi}}
Fix $n \ge 1$,
and let $T$ be a rooted tree on the vertex set $[n+1]$
in~which $r$ is the root and the vertex~$1$ has precisely one child $a$.
Let $T_a$ be the subtree rooted at $a$,
and let $T_r$ the subtree obtained from $T$
by removing $T_a$ and the edge $1a$.
The vertex~1 is a leaf in $T_r$.

Now we create a new tree $T'$, rooted at $a$, as follows:
we start with $T_a$ and then graft $T_r$ by making $r$ a child of $a$.
In the tree $T'$, the vertex~1 is a leaf.
The map $T \mapsto T'$ map is a bijection,
since this construction can be reversed.
(The vertex $r$ can be identified in $T'$
 as the child of $a$ that has 1 as a descendant.)

Clearly, all the proper (resp.~improper) edges in $T$
are still proper (resp.~improper) in $T'$, except that:
\begin{itemize}
   \item[(i)] The edge $1a$ in $T$ is proper, which is deleted in $T'$; and
   \item[(ii)] The edge $ar$ in $T'$ is new and improper,
      since the vertex~$1$ is a descendant of $r$.
\end{itemize}
In particular, the number of vertices with $m$ proper children
is the same in $T$ and $T'$, provided that in $T$ one ignores the vertex~1.
\qed

\subsection{Proof of Propositions~\ref{prop.equivalence} and
      \ref{prop.equivalence.phi}}
   \label{subsec.proof.prop.equivalence}

Now we will prove Proposition~\ref{prop.equivalence},
which asserts the equality of the polynomials $t_{n,k}(y,z)$
defined in \reff{def.tnkyz} using rooted trees
and the polynomials $\widetilde{t}_{n,k}(y,z)$
defined in \reff{def.tnkyz.tilde} using partial functional digraphs.
We will then show that the same argument proves the more
general Proposition~\ref{prop.equivalence.phi},
which asserts the equivalence of the polynomials $t_{n,k}(y,\bphi)$
defined in \reff{def.tnkyphi}
and the polynomials $\widetilde{t}_{n,k}(y,\bphi)$
defined in \reff{def.tnkyphi.tilde}.

We recall that $\scrt^\bullet_n$ denotes the set of rooted trees
on the vertex set $[n]$,
while $\scrt^{\<1;k\>}_n$ denotes the subset
in~which the vertex~1 has $k$ children.
Similarly, $\PFD_n$ denotes the set of partial functional digraphs
on the vertex set $[n]$,
while $\PFD_{n,k}$ denotes the subset in~which there are exactly $k$ vertices
of out-degree 0.

To prove Proposition~\ref{prop.equivalence},
we will construct, for each fixed~$n$,
a bijection $\phi\colon\, \scrt^\bullet_{n+1} \to \PFD_n$
with the following properties:
\begin{itemize}
   \item[(a)] $\phi$ maps $\scrt^{\<1;k\>}_{n+1}$ onto $\PFD_{n,k}$.
   \item[(b)] $\phi$ preserves the number of improper edges.
   \item[(c)] $\phi|_{\scrt^{\<1;k\>}_{n+1}}$ reduces
       the number of proper edges by $k$.
\end{itemize}
We observe that (c) is an immediate consequence of (a) and (b),
since trees in $\scrt^\bullet_{n+1}$ have $n$ edges,
while digraphs in $\PFD_{n,k}$ have $n-k$ edges.

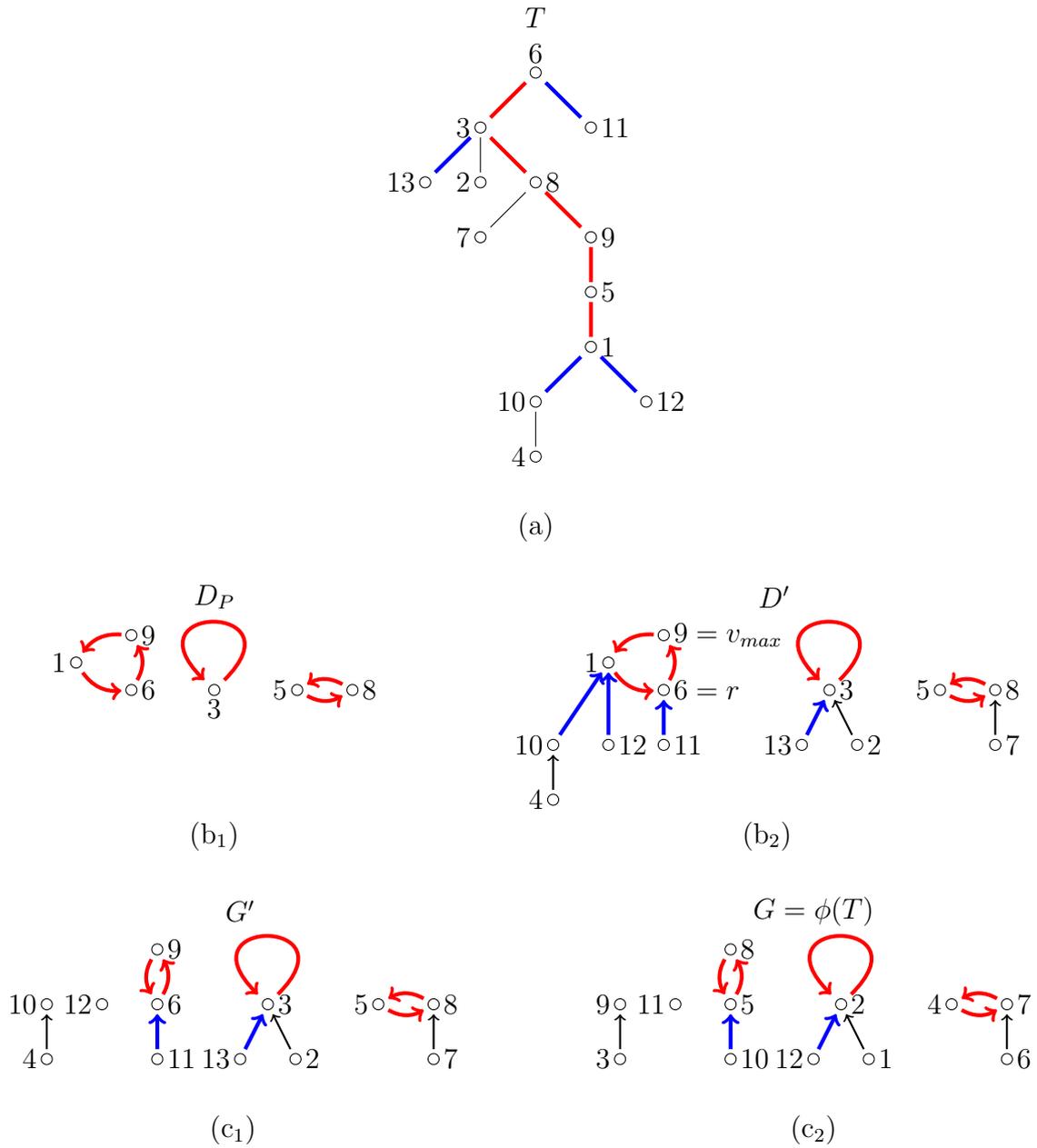
\begin{figure}[p]
  \begin{center}
  \begin{tikzpicture}[scale=0.8]
        \node at (0,1) {$T$};
	\node (1) at (1,-5){};
	\node (2) at (-1,-2){};
	\node (3) at (-1,-1){};
	\node (4) at (0,-7){};
	\node (5) at (1,-4){};
	\node (6) at (0,0){};
	\node (7) at (-1,-3){};
	\node (8) at (0,-2){};
	\node (9) at (1,-3){};
	\node (10) at (0,-6){};
	\node (11) at (1,-1){};
	\node (12) at (2,-6){};
	\node (13) at (-2,-2){};
	\draw (6) circle(.1) node[above] {$6$};
	\foreach \x in {2,3,4,7,10,13}
    		\draw (\x) circle(.1) node[left] {$\x$};
	\foreach \x in {1,5,8,9,11,12}
    		\draw (\x) circle(.1) node[right] {$\x$};
	\draw[ultra thick, red] (6)--(3)--(8)--(9)--(5)--(1);
	\draw[ultra thick, blue] (6)--(11)  (3)--(13)  (10)--(1)--(12);
	\draw (3)--(2)  (8)--(7)  (10)--(4);
  \end{tikzpicture}
  \\[10pt]
  (a)
  \end{center}
%
%
  \begin{center}
  \begin{tabular}{cp{25pt}c}
  \begin{tikzpicture}[baseline=(current bounding box.base),scale=0.8]
 	\tikzset{multiedge/.style={->, bend right}}
  	\node at (0,1.7) {$D_P$};
	\node (1) at (-2.5,.5){};
	\node (6) at (-1.5,0){};
	\node (9) at (-1.5,1){};
	\node (3) at (0,0){};
	\node (5) at (1.5,0){};
	\node (8) at (2.5,0){};
	\draw (3) circle(.1) node[below] {$3$};
	\foreach \x in {1,5}
    		\draw (\x) circle(.1) node[left] {$\x$};
	\foreach \x in {6,8,9}
    		\draw (\x) circle(.1) node[right] {$\x$};
	\Edge[style={multiedge,ultra thick,red}](1)(6)
	\Edge[style={multiedge,ultra thick,red}](6)(9)
	\Edge[style={multiedge,ultra thick,red}](9)(1)
	\Loop[dist=2cm,dir=NO,style={->,ultra thick,red}](3)
	\Edge[style={multiedge,ultra thick,red}](5)(8)
	\Edge[style={multiedge,ultra thick,red}](8)(5)
  \end{tikzpicture}
  &&
  \begin{tikzpicture}[baseline=(current bounding box.base),scale=0.8]
 	\tikzset{multiedge/.style={->, bend right}, diedge/.style={->}, pedge/.style={->, ultra thick, blue}}
  	\node at (-1,1.7) {$D'$};
	\node (1) at (-4,.5){};
	\node (6) at (-3,0){};
	\node (9) at (-3,1){};
	\node (3) at (0,0){};
	\node (5) at (2,0){};
	\node (8) at (3,0){};
	\node (2) at (.5,-1){};
	\node (4) at (-5,-2){};
	\node (7) at (3,-1){};
	\node (10) at (-5,-1){};
	\node (11) at (-3,-1){};
	\node (12) at (-4,-1){};
	\node (13) at (-.5,-1){};
	\foreach \x in {1,4,5,10,13}
    		\draw (\x) circle(.1) node[left] {$\x$};
	\foreach \x in {2,3,7,8,11,12}
    		\draw (\x) circle(.1) node[right] {$\x$};
	\draw (6) circle(.1) node[right] {$6=r$};
	\draw (9) circle(.1) node[right] {$9=v_{max}$};
	\Edge[style={multiedge,ultra thick,red}](1)(6)
	\Edge[style={multiedge,ultra thick,red}](6)(9)
	\Edge[style={multiedge,ultra thick,red}](9)(1)
	\Loop[dist=2cm,dir=NO,style={->,ultra thick,red}](3)
	\Edge[style={multiedge,ultra thick,red}](5)(8)
	\Edge[style={multiedge,ultra thick,red}](8)(5)
	\Edge[style=pedge](10)(1)
	\Edge[style=pedge](12)(1)
	\Edge[style=pedge](11)(6)
	\Edge[style=pedge](13)(3)
	\Edge[style=diedge](4)(10)
	\Edge[style=diedge](2)(3)
	\Edge[style=diedge](7)(8)
  \end{tikzpicture}
  \\[20pt]
  (b${}_1$)  &&  (b${}_2$) 
  \end{tabular}
  \end{center}
%
%
  \begin{center}
  \hspace*{-6mm}
  \begin{tabular}{cp{25pt}c}
  \begin{tikzpicture}[scale=0.8]
	\tikzset{multiedge/.style={->, bend right, ultra thick, red}, diedge/.style={->},  pedge/.style={->, ultra thick, blue}}
	\node at (-0.5,1.7) {$G'$};
	\node (6) at (-2,0){};
	\node (9) at (-2,1){};
	\node (3) at (0,0){};
	\node (5) at (2,0){};
	\node (8) at (3,0){};
	\node (2) at (.5,-1){};
	\node (4) at (-4,-1){};
	\node (7) at (3,-1){};
	\node (10) at (-4,0){};
	\node (11) at (-2,-1){};
	\node (12) at (-3,0){};
	\node (13) at (-.5,-1){};
	\foreach \x in {4,5,10,12,13}
    		\draw (\x) circle(.1) node[left] {$\x$};
	\foreach \x in {2,3,6,7,8,9,11}
    		\draw (\x) circle(.1) node[right] {$\x$};
	\Edge[style=multiedge](6)(9)
	\Edge[style=multiedge](9)(6)
	\Loop[dist=2cm,dir=NO,style={->, ultra thick, red}](3)
	\Edge[style=multiedge](5)(8)
	\Edge[style=multiedge](8)(5)
	\Edge[style=pedge](11)(6)
	\Edge[style=pedge](13)(3)
	\Edge[style=diedge](4)(10)
	\Edge[style=diedge](2)(3)
	\Edge[style=diedge](7)(8)
  \end{tikzpicture}
  &&
  \begin{tikzpicture}[scale=0.8]
	\tikzset{multiedge/.style={->, bend right, ultra thick, red}, diedge/.style={->},  pedge/.style={->, ultra thick, blue}}
        \node at (-0.5,1.7) {$G=\phi(T)$};
	\node (5) at (-2,0){};
	\node (8) at (-2,1){};
	\node (2) at (0,0){};
	\node (4) at (2,0){};
	\node (7) at (3,0){};
	\node (1) at (.5,-1){};
	\node (3) at (-4,-1){};
	\node (6) at (3,-1){};
	\node (9) at (-4,0){};
	\node (10) at (-2,-1){};
	\node (11) at (-3,0){};
	\node (12) at (-.5,-1){};
	\foreach \x in {3,4,9,11,12}
    		\draw (\x) circle(.1) node[left] {$\x$};
	\foreach \x in {1,2,5,6,7,8,10}
    		\draw (\x) circle(.1) node[right] {$\x$};
	\Edge[style=multiedge](5)(8)
	\Edge[style=multiedge](8)(5)
	\Loop[dist=2cm,dir=NO,style={->, ultra thick, red}](2)
	\Edge[style=multiedge](4)(7)
	\Edge[style=multiedge](7)(4)
	\Edge[style=pedge](10)(5)
	\Edge[style=pedge](12)(2)
	\Edge[style=diedge](3)(9)
	\Edge[style=diedge](1)(2)
	\Edge[style=diedge](6)(7)
  \end{tikzpicture}
   \\[10pt]
  (c$_1$) && (c$_2$)
  \end{tabular}
  \end{center}
%
%
  \caption{(a) Tree $T$ in the second model, where $r=v_1=6$,
     $v_{\rm max}=9$, $\sigma=638951=(169)(3)(58)$,
     and vertex~$1$ has two children.
     The backbone edges are shown in red and are improper;
     the other improper edges are shown in black;
     the proper edges are shown in blue.
\newline\hspace*{4mm}
        (b${}_1$,b${}_2$) Partial functional digraphs $D_P$ and $D'$.
         Improper edges arising from the cycles of the permutation $\sigma$
         are shown in red;
         the other improper edges are shown in black;
         the proper edges are shown in blue.
\newline\hspace*{4mm}
        (c${}_1$,c${}_2$) Partial functional digraphs $G'$ and $G$
         in the third model,
         where the two vertices $10$ and $12$ (resp.~$9$ and~$11$)
         have out-degree $0$.
         Improper edges arising from the cycles of the permutation $\sigma$
         are shown in red;
         the other improper edges are shown in black;
         the proper edges are shown in blue.
}
  \label{fig.proof.prop.equivalence}
\end{figure}

\proofof{Proposition~\ref{prop.equivalence}}
(The reader may wish to follow, along with this proof,
 the example shown in Figure~\ref{fig.proof.prop.equivalence}.)

Let $T$ be a rooted tree on the vertex set $[n+1]$
in~which the vertex~$1$ has $k$ children.
Note that the $k$ edges from vertex~1 to its children are all proper.
Now let $P=v_1 \cdots v_{\ell+1}$ ($\ell \ge 0$)
be the unique path in $T$
from the root $v_1 = r$ to the vertex $v_{\ell+1} = 1$;
we call it the ``backbone''.
(Here $\ell=0$ corresponds to the case in~which vertex~1 is the root.)
Removing from $T$ the edges of the path $P$,
we obtain a collection of (possibly trivial) trees $T_1, \ldots, T_{\ell+1}$
rooted at the vertices $v_1, \ldots, v_{\ell+1}$.

Now regard $P$ as a permutation $\sigma$ (written in word form)
of its elements written in increasing order.\footnote{
   That is, let $v'_1 < \ldots < v'_{\ell+1}$
   be the elements of the set $S = \{ v_1, \ldots, v_{\ell+1} \}$
   written in increasing order.
   Then $\sigma$ is the permutation of $S$
   defined by $\sigma(v'_i) = v_i$.
}
In particular, $\sigma(1) = r$ and $\sigma(v_{\rm max}) = 1$
where $v_{\rm max} = \max(v_1, \ldots, v_{\ell+1})$.
Let $D_P$ be the digraph whose vertex set is
$\{v_1, \ldots, v_{\ell+1}\}$,
with edges $\overrightarrow{ij}$ whenever $j= \sigma(i)$.
Then $D_P$ consists of disjoint directed cycles (possibly of length $1$);
it is the representation in cycle form of the permutation $\sigma$.

Now let $D'$ be the digraph obtained from $D_P$
by attaching the trees $T_1, \ldots, T_{\ell+1}$ to $D_P$
(identifying vertices with the same label)
and directing all edges of those trees towards the root.
Then $D'$ is a functional digraph on the vertex set $[n+1]$.
Furthermore, the map $T \mapsto D'$ is a bijection,
since all the above steps can be reversed.

Now let $G'$ be the digraph obtained from $D'$
by deleting the vertex~1 and the $k$ tree edges incident on vertex~1,
and contracting the edges
$\overrightarrow{v_{\rm max} 1}$ and $\overrightarrow{1r}$
into a single edge $\overrightarrow{v_{\rm max} r}$.
Then $G'$ is a digraph on the vertex set $\{2,\ldots,n+1\}$
in~which every vertex has out-degree $1$
except for the $k$ children of vertex~1 in $T$, which have out-degree $0$.
Relabeling all vertices $i \to i-1$,
we obtain a partial functional digraph $G = \phi(T) \in \PFD_{n,k}$.

The step from $D'$ to $G$ can also be reversed:
given a partial functional digraph $G = \PFD_{n,k}$,
we relabel the vertices $i \to i+1$
and then insert the vertex~1
immediately after the largest cyclic vertex of $G$
(if any; otherwise $1$ becomes a loop in $D'$);
all the vertices of out-degree 0 in $G$
are made to point to the vertex~1 in $D'$.

It follows that the map $\phi \colon T \mapsto G$
is a bijection from $\scrt^\bullet_{n+1}$ to $\PFD_n$
that maps $\scrt^{\<1;k\>}_{n+1}$ onto $\PFD_{n,k}$.

Clearly, in the rooted tree $T$,
all the edges in the path $P = v_1 \cdots v_{\ell+1}$ are improper,
since each vertex in $P$ has $v_{\ell+1} = 1$ as its descendant.
These $\ell$ edges correspond, after relabeling,
to $\ell+1$ cyclic edges in the functional digraph $D'$.
These latter edges in turn correspond,
after removal of vertex~1 and contraction of its edges,
to $\ell$ cyclic edges in the partial functional digraph $G'$
(and hence also $G$).
Because they are cyclic edges, they are necessarily improper.
All the other improper/proper edges in $T$ coincide
with improper/proper edges $\overrightarrow{ij}$
in the partial functional digraph $G'$ (and hence $G$)
where $i$ is a transient vertex.
\qed

\medskip

{\bf Remark.}
The first part of this proof (namely, the map $T \mapsto D'$)
is the well-known bijection from doubly-rooted trees
to functional digraphs on the same vertex set
\cite[pp.~224--225]{Labelle_81} \cite[p.~26]{Stanley_99}.
In our application we need the second step to remove the vertex~1
and thereby obtain a map from rooted trees on the vertex set $[n+1]$
to partial functional digraphs on the vertex set $[n]$.
\myendremark

\medskip

\proofof{Proposition~\ref{prop.equivalence.phi}}
In the preceding proof, each vertex $i \neq 1$ in the rooted tree $T$
corresponds to a vertex $i-1$ in the partial functional digraph $G = \phi(T)$.
And for each proper child $j$ of $i$ in $T$,
the proper edge $ij$ in $T$ corresponds to a proper edge
$\overrightarrow{j-1 \: i-1}$ in $G$;
and those are the only proper edges in $G$.
Therefore, if the vertex $i \neq 1$ in $T$ has $m$ proper children,
then the vertex $i-1$ in $G$ has $m$ proper incoming edges.
This proves that $t_{n,k}(y,\bphi) = \widetilde{t}_{n,k}(y,\bphi)$.
\qed

\section{The matrices $\bm{\sfT}$, $\bm{\sfT(y,z)}$ and $\bm{\sfT(y,\bphi)}$
            as exponential Riordan arrays}
   \label{sec.EGF}

In this section we show that the matrices
$\sfT$, $\sfT(y,z)$ and $\sfT(y,\bphi)$
are exponential Riordan arrays $\scrr[F,G]$,
and we compute their generating functions $F$ and $G$
as well as their $A$-, $Z$-, $\Phi$- and $\Psi$-series.

\subsection[The matrix $\sfT$]{The matrix $\bm{\sfT}$}  \label{subsec.EGF.1}

\begin{proposition}
   \label{prop.EGF.1}
Define
\be
   t_{n,k} \;=\;  \binom{n}{k} \, n^{n-k} 
   \;.
\ee
Then the unit-lower-triangular matrix $\sfT = (t_{n,k})_{n,k \ge 0}$
is the exponential Riordan array $\scrr[F,G]$ with
$F(t) = \sum_{n=0}^\infty n^n \, t^n/n!$
and $G(t) = \sum_{n=1}^\infty n^{n-1} \, t^n/n!$.
\end{proposition}

Before proving Proposition~\ref{prop.EGF.1},
let us use it to compute the $A$-, $Z$-, $\Phi$- and $\Psi$-series:

\begin{corollary}
   \label{cor.prop.EGF.1}
The exponential Riordan array $\sfT = \scrr[F,G]$ has
\be
   A(s)  \;=\;  {e^{s} \over 1-s} \,,\quad
   Z(s)  \;=\;  {e^{s} \over (1-s)^2}
 \label{eq.cor.prop.EGF.1.1}
\ee
and
\be
   \Phi(s)  \;=\;  e^{s} \,,\quad
   \Psi(s)  \;=\;  {1 \over 1-s}
   \;.
 \label{eq.cor.prop.EGF.1.2}
\ee
\end{corollary}

\proof
We observe that $G(t)$ is the tree function $T(t)$ \cite{Corless_96},
which satisfies the functional equation $T(t) = t e^{T(t)}$.
Furthermore, we have $F(t) = 1/[1 - T(t)]$:
this well-known fact can be proven using the Lagrange inversion formula
[see \reff{eq.RFG.lagrange} below specialized to $x=0$]
or by various other methods.\footnote{
   {\sc Algebraic proof.}
   $F(t) \:=\:  1 \,+\, t T'(t)
         \:=\:  1 \,+\, \displaystyle{t e^{T(t)} \over 1 - T(t)}
         \:=\:  1 \,+\, \displaystyle{T(t) \over 1 - T(t)}
         \:=\:  {1 \over 1 - T(t)}
   $\,,
   where the first equality used the power series defining $F(t)$ and $T(t)$,
   the second equality used the identity
   $T'(t) = \displaystyle{ e^{T(t)} \over 1 - T(t) }$
   arising from implicit differentiation of the functional equation,
   and the third equality used the functional equation.

   \quad
   {\sc Combinatorial proof.}
   This follows from the identity of combinatorial species:
   endofunctions = permutations $\circ$ rooted trees
   \cite[pp.~41, 43]{Bergeron_98}.
   See also \cite[Exercise~5.32(b)]{Stanley_99}
   for a related combinatorial proof.
}
We now apply Theorem~\ref{thm.riordan.exponential.production}
to determine the functions $A(s)$ and $Z(s)$.
Implicit differentiation of the functional equation yields
$T'(t) = e^{T(t)} / [1 - T(t)]$,
which implies that $A(s) = e^s/(1-s)$.
On the other hand, it follows immediately from the relation
between $F$ and $G$ that $\Psi(s) = 1/(1-s)$.
This implies that $\Phi(s) = e^s$ and $Z(s) = e^s/(1-s)^2$.
\qed

We will give five proofs of Proposition~\ref{prop.EGF.1}:
a direct algebraic proof using Lagrange inversion and an Abel identity;
an inductive algebraic proof, using a different Abel identity;
a third algebraic proof using the $A$- and $Z$-sequences
of an ordinary Riordan array;
a combinatorial proof using exponential generating functions
based on the interpretation of $t_{n,k}$
as counting partial functional digraphs;
and a bijective combinatorial proof based on the interpretation of $t_{n,k}$
as counting rooted labeled trees
according to the number of children of the root
that are lower-numbered than the root.
In Section~\ref{subsec.EGF.2}
we will give yet another combinatorial proof
(also using exponential generating functions),
this time based on the interpretation of $t_{n,k}$
as counting rooted labeled trees
according to the number of children of a specified vertex~$i$;
but this proof will be given in the more general context
of the polynomials $t_{n,k}(y,z)$.

\firstproofof{Proposition~\ref{prop.EGF.1}}
The tree function $T(t)$ satisfies the functional equation $T(t) = t e^{T(t)}$.
We use Lagrange inversion \reff{eq.lagrange.H}
with $\phi(u) = e^u$ and $H(u) = e^{xu}/(1-u)$: this gives
\begin{subeqnarray}
   [t^n] \: {e^{x T(t)} \over 1 - T(t)}
   & = &
   {1 \over n} \, [u^{n-1}]
               \, \biggl( {x \over 1-u} \,+\, {1 \over (1-u)^2} \biggr)
               \, e^{(x+n)u}
       \\[2mm]
   & = &
   {1 \over n} \, \sum_{k=0}^{n-1} (x+k+1) \: {(x+n)^{n-1-k} \over (n-1-k)!}
       \\[2mm]
   & = &
   {1 \over n!} \, \sum_{k=0}^{n-1} \binom{n-1}{k} \, k! \;
                                 (x+k+1) \: {(x+n)^{n-1-k} \over (n-1-k)!}
       \\[2mm]
   & = &
   {(x+n)^n \over n!}
   \;,
 \label{eq.RFG.lagrange}
\end{subeqnarray}
where the last step used an Abel identity
\cite[p.~21, eq.~(25) with $n \to n-1$ and $x \to x+1$]{Riordan_68}.
In view of \reff{eq.Tn.explicit}, this proves \reff{eq.tnk.egf},
which by \reff{eq.egf.RFG} proves that $\sfT = \scrr[F,G]$.
\qed

\secondproofof{Proposition~\ref{prop.EGF.1}}
It is immediate that the zeroth column of $\sfT$
has exponential generating function $F(t) = \sum_{n=0}^\infty n^n \, t^n/n!$.
We now show by induction on $k$ that
the $k$th column has egf $F(t) \, G(t)^k / k!$
where $G(t) = \sum_{n=0}^\infty n^{n-1} \, t^n/n!$:
that~is, we need to show that the $k$th column has egf
equal to $G(t)/k$ times the egf of the $(k-1)$st column,
or in other words
\be
   k! \: t_{n,k}
   \;=\;
   \sum_{j=1}^{n-k+1} \binom{n}{j} \: j^{j-1} \: (k-1)! \: t_{n-j,k-1}
 \label{eq.prop.EGF.1.proof1}
\ee
for $k \ge 1$.
We start from the Abel identity \cite[p.~18, eq.~(13a)]{Riordan_68}
\be
   x^{-1} \, (x+y+n)^n
   \;=\;
   \sum_{i=0}^n \binom{n}{i} \, (x+i)^{i-1} \, (y+n-i)^{n-i}
   \;.
 \label{eq.abel.13a}
\ee
Now substitute $x=1$ and $y \to y-n-1$, divide both sides by $n!$,
and relabel $i = j-1$:  the result is
\be
   {y^n \over n!}
   \;=\;
   \sum_{j=1}^{n+1} {j^{j-1} \over j!} \: {(y-j)^{n+1-j} \over (n-j+1)!}
   \;.
 \label{eq.abel.13a.bis}
\ee
Next substitute $n \to n-k$ and then set $y=n$:
\be
   {n^{n-k} \over (n-k)!}
   \;=\;
   \sum_{j=1}^{n-k+1} {j^{j-1} \over j!} \: {(n-j)^{n-j-k+1} \over (n-j-k+1)!}
   \;.
 \label{eq.abel.13a.bis2}
\ee
Multiplying this by $n!$ yields
\be
   {n! \over (n-k)!} \, n^{n-k}
   \;=\;
   \sum_{j=1}^{n-k+1} \binom{n}{j} \, j^{j-1} \, {(n-j)! \over (n-j-k+1)!}
       \, (n-j)^{n-j-k+1}
   \;,
 \label{eq.abel.13a.bis3}
\ee
which is \reff{eq.prop.EGF.1.proof1}.
\qed

\thirdproofof{Proposition~\ref{prop.EGF.1}}
Showing that $(t_{n,k})_{n,k \ge 0}$
equals the exponential Riordan array $\scrr[F,G]$
is equivalent to showing that
$((k!/n!) \, t_{n,k})_{n,k \ge 0}$
equals the ordinary Riordan array $\scrr(F,G)$.
We write $r_{n,k} \eqdef (k!/n!) \, t_{n,k} = n^{n-k}/(n-k)!$
and $R = (r_{n,k})_{n,k \ge 0}$.
By the binomial theorem we have
\be
   r_{n+1,k+1}
   \;\eqdef\;
   {(n+1)^{n-k} \over (n-k)!}
   \;=\;
   \sum_{j=0}^{n-k} {1 \over j!} \: {n^{n-k-j} \over (n-k-j)!}
   \;\eqdef\;
   \sum_{j=0}^{n-k} {1 \over j!} \: r_{n,k+j}
\ee
for all $k \ge 0$.
So the matrix $R$ satisfies the appropriate identities
to have the $A$-sequence $a_j = 1/j!$.
Similarly, by the binomial theorem we have
\be
   r_{n+1,0}
   \;\eqdef\;
   {(n+1)^{n+1} \over (n+1)!}
   \;=\;
   {(n+1)^n \over n!}
   \;=\;
   \sum_{j=0}^n {1 \over j!} \, {n^{n-j} \over (n-j)!}
   \;\eqdef\;
   \sum_{j=0}^n {1 \over j!} \: r_{n,j}
   \;.
\ee
So the matrix $R$ satisfies the appropriate identities
to have the $Z$-sequence $z_j = 1/j!$.
It follows from Theorem~\ref{thm.riordan.production}
that $R$ is an ordinary Riordan array $\scrr(F,G)$
where $F$ and $G$ are given by \reff{eq.thm.riordan.production}
with $A(t) = Z(t) = e^t$.
By the Lagrange inversion formula we find
\be
   [t^n] \, G(t)
   \;=\;
   {1 \over n} \, [t^{n-1}] \, A(t)^n
   \;=\;
   {n^{n-1} \over n!}
   \;.
\ee
And the ordinary generating function of the zeroth column of $R$
is obviously $F(t) = \sum_{n=0}^\infty n^n \, t^n/n!$.
\qed

{\bf Remark.}
The identity \reff{eq.abel.13a.bis2} in the second proof can be written
for $k \ge 1$ as
\be
   r_{n,k}  \;=\; \sum_{j=1}^{n-k+1} {j^{j-1} \over j!} \: r_{n-j,k-1}
   \;,
\ee
which shows that the ordinary generating function of the $k$th column of $R$
equals $G(t) = \sum_{j=1}^\infty j^{j-1} \, t^j/j!$
times the ordinary generating function of the $(k-1)$st column.
Combining this with the fact that
the ordinary generating function of the zeroth column of $R$ is
$F(t) = \sum_{n=0}^\infty n^n \, t^n/n!$
gives an alternate proof that $R = \scrr(F,G)$.
\myendremark

\fourthproofof{Proposition~\ref{prop.EGF.1}}
We begin from the fact that $t_{n,k} = \binom{n}{k} n^{n-k}$ counts
partial functional digraphs on $n$ labeled vertices
that have $k$ vertices of out-degree~0.
Such a partial functional digraph is the disjoint union
of $k$ rooted trees (rooted at the vertices of out-degree~0)
together with a functional digraph on the remaining vertices.
Standard enumerative arguments then imply that the
exponential generating function for the numbers $t_{n,k}$ is
\be
   \sum_{n=k}^\infty t_{n,k} \, {t^n \over n!}
   \;=\;
   F(t) \, {T(t)^k \over k!}
\ee
where $F(t) = \sum_{n=0}^\infty n^n \, t^n/n!$
is the exponential generating function for functional digraphs
and $T(t) = \sum_{n=1}^\infty n^{n-1} \, t^n/n!$
is the exponential generating function for rooted trees.
\qed

\fifthproofof{Proposition~\ref{prop.EGF.1}}
We begin from the fact \cite{Chauve_99,Chauve_00,Sokal_trees_enumeration}
that $t_{n,k}$ equals the number of rooted trees on the vertex set $[n+1]$
in~which exactly $k$ children of the root are lower-numbered than the root.
We will prove \reff{eq.prop.EGF.1.proof1} in the form
\be
   k \, t_{n,k}
   \;=\;
   \sum_{j=1}^{n-k+1} \binom{n}{j} \: j^{j-1} \: t_{n-j,k-1}
 \label{eq.prop.EGF.1.thirdproof1}
\ee
for $k \ge 1$.
We interpret $k \, t_{n,k}$ as the number of triplets $(T,r,v_\star)$
in~which $(T,r)$ is a rooted tree on the vertex set $[n+1]$
in~which exactly $k$ children of the root are lower-numbered than the root,
and $v_\star$ is one of those lower-numbered children
(we call it the ``marked vertex'').
See Figure~\ref{fig:T}.
We interpret $j^{j-1}$ as the number of rooted trees
on $j$ labeled vertices.
So the summand on the right-hand side of \reff{eq.prop.EGF.1.thirdproof1}
enumerates quintuplets $(A,T_1,r_1,T_2,r_2)$
where $A$ is a subset of $[n]$ of cardinality $j$,
$(T_1,r_1)$ is a rooted tree on the vertex set $A$,
and $(T_2,r_2)$ is a rooted tree on the vertex set $[n+1] \setminus A$
in~which exactly $k-1$ children of the root are lower-numbered than the root.
See Figure~\ref{fig:T1&T2}.

\begin{figure}[p]
\begin{center}
\begin{tikzpicture}
	\fill (0,0) circle(.1) node[above] {$r$};
	\draw (-2.5,-1) circle(.08) node[left] {$v_1$};
	\fill (-1,-1) circle(.1) node[above=1pt] {$v_{\star}$};
	\draw (.5,-1) circle(.08) node[right] {$v_k$};
	\foreach \x in {-1.75,-.25}
	{	\node at (\x,-1) {$\cdots$};
		\node at (\x,-3) {$\cdots$};
	}
	\draw (2,-1) ellipse(.8 and .5) node (H) {$H$};
	\draw (-2.5,-3) ellipse(.3 and .8) node (S1) {$S_1$};
	\draw (-1,-3) ellipse(.3 and .8) node (S*) {$S_{\star}$};
	\draw (.5,-3) ellipse(.3 and .8) node (Sk) {$S_k$};
	\draw (2,-3) circle(.8) node (H') {$H'$};
	\draw (2,-2.2)--(2,-1.5)
		  (2,-.5)--(0,0)--(.5,-1)--(.5,-2.2)
		  (-2.5,-2.2)--(-2.5,-1)--(0,0)--(-1,-1)--(-1,-2.2);
	\node at (0,1) {Tree $T$ on vertex set $[n+1]$};
\end{tikzpicture}
\caption{A triplet $(T,r,v_{\star})$, where $v_1,\ldots,v_k$ are
         the children of the root $r$ that are lower-numbered than $r$,
         and $H$ depicts the children of $r$ that are higher-numbered than $r$.
        }
\label{fig:T}
\end{center}
\end{figure}
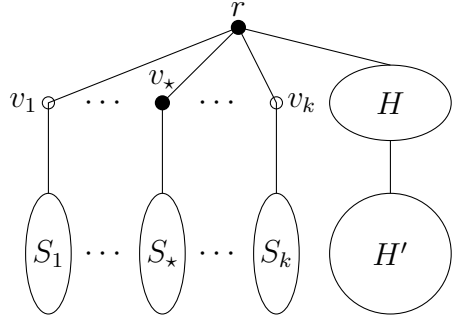

\begin{figure}[p]
\begin{center}
\vspace*{5mm}
\begin{tikzpicture}
	\fill (0,0) circle(.1) node[above] {$r_1$};
	\draw (-1,-1) ellipse(.8 and .5) node (L1) {$L_1$};
	\draw (1,-1) ellipse(.8 and .5) node (H1) {$H_1$};
	\draw (-1,-3) circle(.8) node (L1') {$L_1'$};
	\draw (1,-3) circle(.8) node (H1') {$H_1'$};
	\draw (-1,-2.2)--(-1,-1.5)
		  (-1,-.5)--(0,0)--(1,-.5)
		  (1,-1.5)--(1,-2.2);
	\node at (0,1) {Tree $T_1$ on vertex set $A \subseteq [n]$};
\end{tikzpicture}
\hspace{50pt}
\begin{tikzpicture}
	\fill[black!25!white] (0,.2) circle(.5);
	\fill (0,0) circle(.1) node[above] {$r_2$};
	\draw (-1,-1) ellipse(.8 and .5) node (L2) {$L_2$};
	\filldraw[fill=black!25!white] (1,-1) ellipse(.8 and .5) node (H2) {$H_2$};
	\filldraw[fill=black!25!white] (-1,-3) circle(.8) node (L2') {$L_2'$};
	\filldraw[fill=black!25!white] (1,-3) circle(.8) node (H2') {$H_2'$};
	\draw (-1,-2.2)--(-1,-1.5)
		  (-1,-.5)--(0,0)--(1,-.5)
		  (1,-1.5)--(1,-2.2);
	\node at (0,1) {Tree $T_2$ on vertex set $[n+1]\setminus A$};
\end{tikzpicture}
\caption{A quintuplet $(A,T_1,r_1,T_2,r_2)$.
         Here $L_i$ (resp.\ $H_i$) depicts the children of the root $r_i$
         that are lower-numbered (resp.\ higher-numbered) than $r_i$.
         In this figure and the following ones,
         the shaded parts are the possible locations of the vertex $n+1$.
        }
\label{fig:T1&T2}
\end{center}
\end{figure}
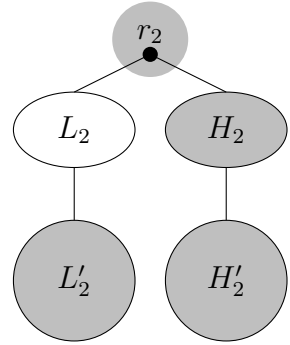

\begin{figure}[p]
\begin{center}
\vspace*{5mm}
\begin{tabular}{cp{1cm}c}
\begin{tikzpicture}[scale=0.8]
	\fill (-4,-1) circle(.1) node[above] {$r_1$};
	\draw (-5,-2) ellipse(.8 and .5) node (L1) {$L_1$};
	\draw (-3,-2) ellipse(.8 and .5) node (H1) {$H_1$};
	\draw (-5,-4) circle(.8) node (L1') {$L_1'$};
	\draw (-3,-4) circle(.8) node (H1') {$H_1'$};
	\draw (-5,-3.2)--(-5,-2.5)
		  (-5,-1.5)--(-4,-1)--(-3,-1.5)
		  (-3,-2.5)--(-3,-3.2);
	\fill[black!25!white] (0,.2) circle(.5);
	\fill (0,0) circle(.1) node[above] {$r_2$};
	\draw (-1,-1) ellipse(.8 and .5) node (L2) {$L_2$};
	\filldraw[fill=black!25!white] (1,-1) ellipse(.8 and .5) node (H2) {$H_2$};
	\filldraw[fill=black!25!white] (-1,-3) circle(.8) node (L2') {$L_2'$};
	\filldraw[fill=black!25!white] (1,-3) circle(.8) node (H2') {$H_2'$};
	\draw (-1,-2.2)--(-1,-1.5)
		  (-1,-.5)--(0,0)--(1,-.5)
		  (1,-1.5)--(1,-2.2);
	\draw[ultra thick] (0,0)--(-4,-1);
\end{tikzpicture}
& &
\begin{tikzpicture}[scale=0.8]
	\fill (0,0) circle(.1) node[above] {$r_1$};
	\draw (-1,-1) ellipse(.8 and .5) node (L2) {$L_2$};
	\draw (1,-1) ellipse(.8 and .5) node (H1) {$H_1$};
	\filldraw[fill=black!25!white] (-1,-3) circle(.8) node (L2') {$L_2'$};
	\draw (1,-3) circle(.8) node (H1') {$H_1'$};
	\draw (-1,-2.2)--(-1,-1.5)
		  (-1,-.5)--(0,0)--(1,-.5)
		  (1,-1.5)--(1,-2.2);
	\fill (4,-1) circle(.1) node[above] {$r_2$};
	\draw (3,-2) ellipse(.8 and .5) node (L1) {$L_1$};
	\filldraw[fill=black!25!white] (5,-2) ellipse(.8 and .5) node (H2) {$H_2$};
	\draw (3,-4) circle(.8) node (L1') {$L_1'$};
	\filldraw[fill=black!25!white] (5,-4) circle(.8) node (H2') {$H_2'$};
	\draw (3,-3.2)--(3,-2.5)
		  (3,-1.5)--(4,-1)--(5,-1.5)
		  (5,-2.5)--(5,-3.2);
	\draw[ultra thick] (0,0)--(4,-1);
	\fill (-1.5,-1) circle(.1);
	\draw[->] (-1.5,-1)--(-2,-1.5) node[below] {$v_{\star}$};
\end{tikzpicture}
\\[5pt]
(a) Case I: $r_1 < r_2$ && (b) Case II: $r_1 > r_2$
\end{tabular}
\caption{Bijection RHS $\Rightarrow$ LHS.}
\label{fig:bijection.RtoL}
\end{center}
\end{figure}
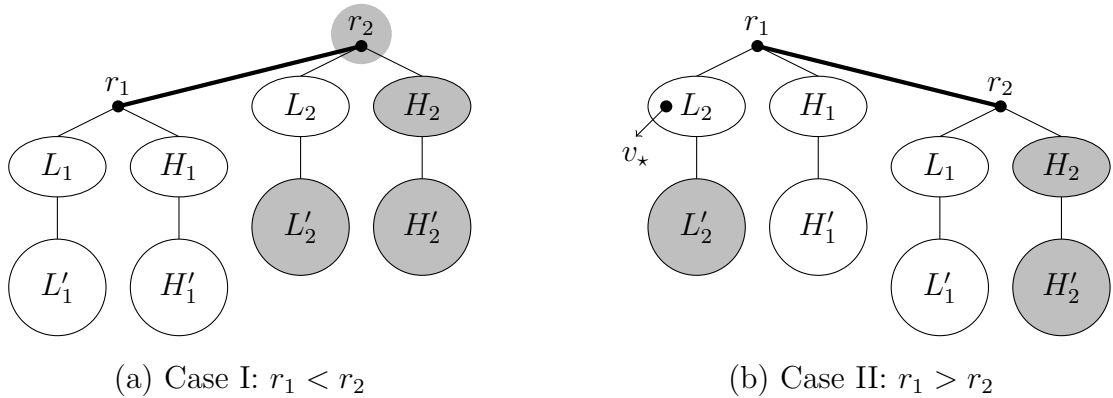

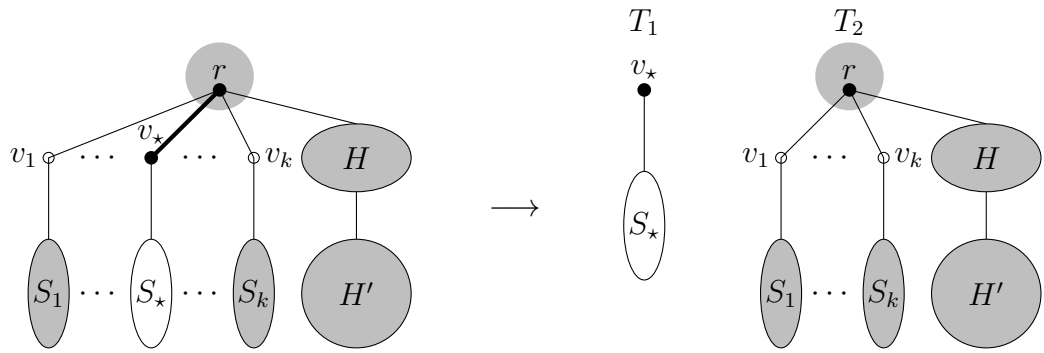
\begin{figure}[p]
\begin{center}
\begin{tikzpicture}[scale=0.9]
	\fill[black!25!white] (0,.2) circle(.5);
	\fill (0,0) circle(.1) node[above] {$r$};
	\draw (-2.5,-1) circle(.08) node[left] {$v_1$};
	\fill (-1,-1) circle(.1) node[above=1pt] {$v_{\star}$};
	\draw (.5,-1) circle(.08) node[right] {$v_k$};
	\foreach \x in {-1.75,-.25}
		\foreach \y in {-1,-3}
			\node at (\x,\y) {$\cdots$};
	\filldraw[fill=black!25!white] (2,-1) ellipse(.8 and .5) node (H) {$H$};
	\filldraw[fill=black!25!white] (-2.5,-3) ellipse(.3 and .8) node (S1) {$S_1$};
	\draw (-1,-3) ellipse(.3 and .8) node (S*) {$S_{\star}$};
	\filldraw[fill=black!25!white] (.5,-3) ellipse(.3 and .8) node (Sk) {$S_k$};
	\filldraw[fill=black!25!white] (2,-3) circle(.8) node (H') {$H'$};
	\draw (2,-2.2)--(2,-1.5)
		  (2,-.5)--(0,0)--(.5,-1)--(.5,-2.2)
		  (-2.5,-2.2)--(-2.5,-1)--(0,0)--(-1,-1)--(-1,-2.2);
	\draw[ultra thick] (0,0)--(-1,-1);
\end{tikzpicture}
\hspace{25pt}\raisebox{50pt}{$\longrightarrow$}\hspace{25pt}
\begin{tikzpicture}[scale=0.9]
	\fill (-3,0) circle(.1) node[above] {$v_{\star}$};
	\draw (-3,-2) ellipse(.3 and .8) node (S*) {$S_{\star}$};
	\draw (-3,0)--(-3,-1.2);
	\fill[black!25!white] (0,.2) circle(.5);
	\fill (0,0) circle(.1) node[above] {$r$};
	\draw (-1,-1) circle(.08) node[left] {$v_1$};
	\draw (.5,-1) circle(.08) node[right] {$v_k$};
	\foreach \y in {-1,-3} \node at (-.25,\y) {$\cdots$};
	\filldraw[fill=black!25!white] (2,-1) ellipse(.8 and .5) node (H) {$H$};
	\filldraw[fill=black!25!white] (-1,-3) ellipse(.3 and .8) node (S1) {$S_1$};
	\filldraw[fill=black!25!white] (.5,-3) ellipse(.3 and .8) node (Sk) {$S_k$};
	\filldraw[fill=black!25!white] (2,-3) circle(.8) node (H') {$H'$};
	\draw (2,-2.2)--(2,-1.5)
		  (2,-.5)--(0,0)--(.5,-1)--(.5,-2.2)
		  (-1,-2.2)--(-1,-1)--(0,0);
	\foreach \x in {1,2} \node at (3*\x-6,1) {$T_{\x}$};
\end{tikzpicture}
\caption{Bijection LHS $\Rightarrow$ RHS: Case I.}
\label{fig:bijection.LtoR1}
\end{center}
\end{figure}

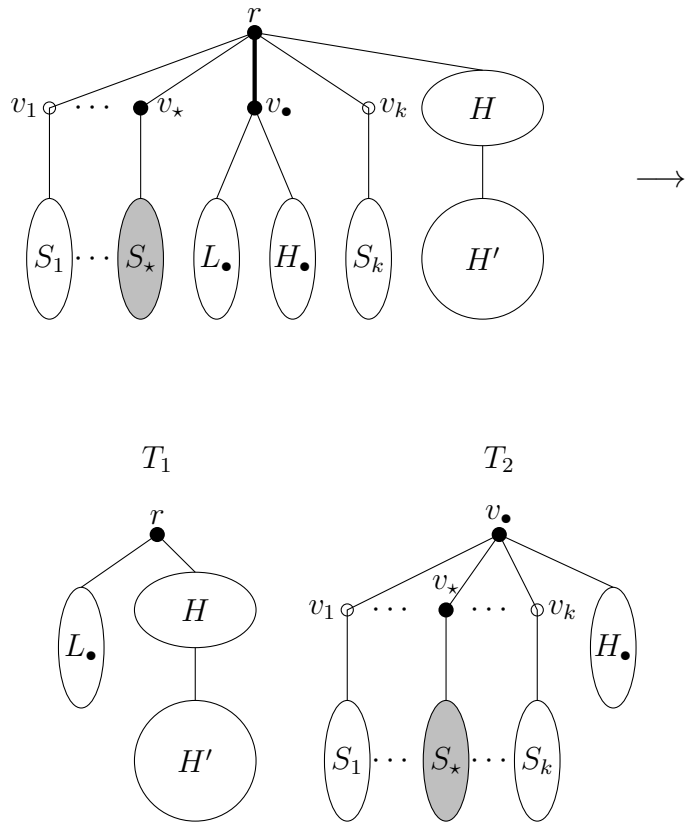
\begin{figure}[p]
\begin{center}
\begin{tikzpicture}
	\fill (0,0) circle(.1) node[above] {$r$};
	\draw (-2.7,-1) circle(.08) node[left] {$v_1$};
	\fill (-1.5,-1) circle(.1) node[right=2pt] {$v_{\star}$};
	\fill (0,-1) circle(.1) node[right] {$v_{\bullet}$};
	\draw (1.5,-1) circle(.08) node[right] {$v_k$};
	\draw (3,-1) ellipse(.8 and .5) node (H) {$H$};
	\draw (-2.7,-3) ellipse(.3 and .8) node (S1) {$S_1$};
	\foreach \y in {-1,-3} \node at (-2.1,\y) {$\cdots$};
	\filldraw[fill=black!25!white] (-1.5,-3) ellipse(.3 and .8) node (S*) {$S_{\star}$};
	\draw(-.5,-3) ellipse(.3 and .8) node (Lb) {$L_{\bullet}$};
	\draw(.5,-3) ellipse(.3 and .8) node (Hb) {$H_{\bullet}$};
	\draw(1.5,-3) ellipse(.3 and .8) node (Sk) {$S_k$};
	\draw (3,-3) circle(.8) node (H') {$H'$};
	\draw (3,-2.2)--(3,-1.5)
		 (3,-.5)--(0,0)--(1.5,-1)--(1.5,-2.2)
		 (-2.7,-2.2)--(-2.7,-1)--(0,0)--(-1.5,-1)--(-1.5,-2.2)
		 (-.5,-2.2)--(0,-1)--(.5,-2.2);
	\draw[ultra thick] (0,0)--(0,-1);
\end{tikzpicture}
\hspace*{30pt}\raisebox{50pt}{$\longrightarrow$}
\vspace*{1.5cm}

\begin{tikzpicture}
	\fill (-4.5,0) circle(.1) node[above] {$r$};
	\draw(-5.5,-1.5) ellipse(.3 and .8) node (Lb) {$L_{\bullet}$};
	\draw (-4,-1) ellipse(.8 and .5) node (H) {$H$};
	\draw (-4,-3) circle(.8) node (H') {$H'$};
	\draw (-5.5,-.7)--(-4.5,0)--(-4,-.5)
		  (-4,-1.5)--(-4,-2.2);
	\fill (0,0) circle(.1) node[above] {$v_{\bullet}$};
	\draw (-2,-1) circle(.08) node[left] {$v_1$};
	\fill (-.7,-1) circle(.1) node[above=2pt] {$v_{\star}$};
	\draw (.5,-1) circle(.08) node[right] {$v_k$};
	\draw (-2,-3) ellipse(.3 and .8) node (S1) {$S_1$};
	\foreach \x in {-1.4,-.1}
		\foreach \y in {-1,-3}
			\node at (\x,\y) {$\cdots$};
	\filldraw[fill=black!25!white] (-.7,-3) ellipse(.3 and .8) node (S*) {$S_{\star}$};
	\draw(.5,-3) ellipse(.3 and .8) node (Sk) {$S_k$};
	\draw(1.5,-1.5) ellipse(.3 and .8) node (Hb) {$H_{\bullet}$};
	\draw (-2,-2.2)--(-2,-1)--(0,0)--(-.7,-1)--(-.7,-2.2)
		 (.5,-2.2)--(.5,-1)--(0,0)--(1.5,-.7);
	\foreach \x in {1,2} \node at (4.5*\x-9,1) {$T_{\x}$};
\end{tikzpicture}
\caption{Bijection LHS $\Rightarrow$ RHS: Case II.}
\label{fig:bijection.LtoR2}
\end{center}
\end{figure}

\medskip

{\bf Bijection RHS $\Longrightarrow$ LHS.}
Given the quintuplet $(A,T_1,r_1,T_2,r_2)$,
we construct a triplet $(T,r,v_\star)$ as follows.
We distinguish two cases:
\begin{itemize}
   \item {\bf Case I: $\bm{r_1 < r_2}$.}  We let $r_2$ be the new root
and add an edge making $r_1$ a child of $r_2$;  this gives $(T,r)$.
We then mark the vertex $v_\star = r_1$.
Please note that in this case the vertex $n+1$
is not a descendant of $v_\star$
(see Figure~\ref{fig:bijection.RtoL}a).
   \item {\bf Case II: $\bm{r_1 > r_2}$.}  We let $r_1$ be the new root
and add an edge making $r_2$ a child of $r_1$;
we then interchange the lower-numbered children of $r_1$
(together with all their descendants)
with the lower-numbered children of $r_2$
(and their descendants).
This gives $(T,r)$.
We observe (see Figure~\ref{fig:bijection.RtoL}b)
that $r_2$ is the largest-numbered among all the lower-numbered
children of $r$ in $T$.
We observe also that the vertex $n+1$ must be a descendant
of some lower-numbered child of $r_1$ in $T$;
we set the marked vertex $v_\star$ to be this lower-numbered child.
Note that $v_\star$ must either belong to the set $S_2$
(consisting of the lower-numbered children of $r_2$ in $T_2$,
 which became lower-numbered children of $r_1$ in $T$)
or else be the vertex $r_2$.
\end{itemize}
In both cases, in the rooted tree $(T,r)$,
exactly $k$ children of the root are lower-numbered than the root.

\medskip

We now describe the inverse bijection:

\medskip

{\bf Bijection LHS $\Longrightarrow$ RHS.}
Given the triplet $(T,r,v_\star)$,
we reconstruct the quintuplet $(A,T_1,r_1,T_2,r_2)$ as follows.
We distinguish two cases:
\begin{itemize}
   \item {\bf Case I: $\bm{n+1}$ is not a descendant of $\bm{v_\star}$.}
We delete the edge between $r$ (the root of $T$) and $v_\star$.
Then $(T_1,r_1)$ is the tree whose root is $v_\star$,
and $A$ is its vertex set;
$(T_2,r_2)$ is the tree whose root is $r$,
and its vertex set is $[n+1] \setminus A$.
Since $n+1$ is not a descendant of $v_\star$,
it $n+1$ belongs to $T_2$, so that $A \subseteq [n]$.
And we have $r_1 = v_\star < r = r_2$.
See Figure~\ref{fig:bijection.LtoR1}.
   \item {\bf Case II: $\bm{n+1}$ is a descendant of $\bm{v_\star}$.}
The root $r$ of $T$ has $k$ ($\ge 1$) lower-numbered children;
let $v_\bullet$ be the largest-numbered of these.
We delete the edge between $r$ and $v_\bullet$;
then we interchange the lower-numbered children of $r$
(together with all their descendants)
with the lower-numbered children of $v_\bullet$
(and their descendants).
Then $(T_1,r_1)$ is the tree whose root is $r$,
and $A$ is its vertex set;
$(T_2,r_2)$ is the tree whose root is $v_\bullet$,
and its vertex set is $[n+1] \setminus A$.
Please observe that the marked vertex $v_\star$
was a lower-numbered child of $r$ in $T$;
therefore, it is either equal to $v_\bullet = r_2$
or else becomes a lower-numbered child of $v_\bullet = r_2$ in $T_2$.
Since the vertex $n+1$ was a descendant of $v_\star$, it must belong to $T_2$;
therefore $A \subseteq [n]$.
And we have $r_1 = r > v_\bullet = r_2$.
See Figure~\ref{fig:bijection.LtoR2}.
\end{itemize}
In both cases, in the rooted tree $(T_2,r_2)$,
exactly $k-1$ children of the root are lower-numbered than the root.
\qed

\subsection[The matrix $\sfT(y,z)$]{The matrix $\bm{\sfT(y,z)}$}
    \label{subsec.EGF.2}

We now prove that the matrix $\sfT(y,z) = (t_{n,k}(y,z))_{n,k \ge 0}$
is an exponential Riordan array $\scrr[F,G]$,
and we compute $F$ and $G$.
Most of this computation was done a quarter-century ago
by Dumont and Ramamonjisoa \cite{Dumont_96}:
their arguments handled the case $k=0$,
and we extend those arguments slightly to handle the case of general $k$.
Our presentation follows the notation of \cite{forests_totalpos}.

Let $\scrt^\bullet_n$ denote the set of rooted trees on the vertex set $[n]$;
let $\scrt^{[i]}_n$ denote the subset of $\scrt^\bullet_n$
in~which the root vertex is $i$;
and let $\scrt^{\<i;k\>}_n$ denote the subset of $\scrt^\bullet_n$
in~which the vertex $i$ has $k$ children.
Given a tree $T \in \scrt^\bullet_n$,
we write $\imprope(T)$ for the number of improper edges of $T$.
Now define the generating polynomials
\begin{eqnarray}
   R_n(y,z)
   & = &
   \sum_{T \in \scrt^\bullet_n} y^{\imprope(T)} z^{n-1-\imprope(T)}
   \qquad\qquad
        \label{def.Rn}  \\[2mm]
   S_n(y,z)
   & = &
   \sum_{T \in \scrt^{[1]}_{n+1}} y^{\imprope(T)} z^{n-\imprope(T)}
        \label{def.Sn}  \\[2mm]
   A_{n,k}(y,z)  \;=\;  t_{n,k}(y,z)
   & = &
   \sum_{T \in \scrt^{\<1;k\>}_{n+1}} y^{\imprope(T)} z^{n-k-\imprope(T)}
        \label{def.An}
\end{eqnarray}
in~which each improper (resp.~proper) edge gets a weight $y$ (resp.~$z$)
except that in $A_{n,k}$ the $k$ proper edges connecting
the vertex~1 to its children are unweighted.
And then define the exponential generating functions
\begin{eqnarray}
   \scrr(t;y,z)
   & = &  \sum\limits_{n=1}^\infty R_n(y,z) \: {t^n \over n!}
   \qquad\qquad\qquad
        \label{def.Rn.EGF}  \\[2mm]
   \scrs(t;y,z)  & = &  \sum\limits_{n=0}^\infty S_n(y,z) \: {t^n \over n!}
        \label{def.Sn.EGF}  \\[2mm]
   \scra_k(t;y,z)  & = &  \sum\limits_{n=0}^\infty A_{n,k}(y,z) \: {t^n \over n!}
        \label{def.An.EGF}
\end{eqnarray}
We will then prove the following key result,
which is a slight extension of \cite[Proposition~7]{Dumont_96}
to handle the case $k \neq 0$:

\begin{proposition}
   \label{prop.dumont}
The series $\scrr$, $\scrs$ and $\scra_k$ satisfy the following identities:
\begin{itemize}
   \item[(a)]  $\scrs(t;y,z) \:=\: \exp\big[ z \, \scrr(t;y,z) \bigr]$
   \item[(b)]  $\scra_k(t;y,z) \:=\:
       \displaystyle{\scrr(t;y,z)^k/k! \over 1 - y \scrr(t;y,z)}$
   \item[(c)]  $\displaystyle {d \over dt} \scrr(t;y,z)
                             \:=\: \scra_0(t;y,z) \, \scrs(t;y,z)$
\end{itemize}
and hence
\begin{itemize}
   \item[(d)] $\displaystyle {d \over dt} \scrr(t;y,z)
                             \:=\: 
               {\exp\big[ z \, \scrr(t;y,z) \bigr]
                \over
                1 - y \scrr(t;y,z)}$
\end{itemize}
\end{proposition}

Solving the differential equation of Proposition~\ref{prop.dumont}(d)
with the initial condition
$\scrr(0;y,z) = 0$, we obtain:

\begin{corollary}
   \label{cor.dumont}
The series $\scrr(t;y,z)$ satisfies the functional equation
\be
   y - z + yz \scrr  \;=\;  (y - z + z^2 t) \, e^{z\scrr}
 \label{eq.cor.dumont.1}
\ee
and hence has the solution
\be
   \scrr(t;y,z)
   \;=\;
   {1 \over z}
   \biggl[ T\Big( \Big(1 - {z \over y} + {z^2 \over y} t \Big) \:
                  e^{- \, \displaystyle \big(1 - {z \over y} \big)}
            \Big)
            \:-\: \Big(1 - {z \over y} \Big)
   \biggr]
 \label{eq.cor.dumont.2}
\ee
where $T(t)$ is the tree function \reff{def.treefn}.
\end{corollary}

Comparing Proposition~\ref{prop.dumont}(b)
with the definition \reff{def.RFG} of exponential Riordan arrays,
we conclude:

\begin{corollary}
   \label{cor.dumont2}
The matrix $\sfT(y,z)$
is the exponential Riordan array $\scrr[F,G]$ where
\be
   F(t)  \;=\;  {1 \over 1 - y \scrr(t;y,z)} \,,\quad
   G(t)  \;=\;  \scrr(t;y,z)
 \label{eq.cor.dumont2.1}
\ee
and $\scrr(t;y,z)$ is given by \reff{eq.cor.dumont.2}.
\end{corollary}

And comparing Proposition~\ref{prop.dumont}(b,d)
with the definitions
\reff{eq.prop.riordan.exponential.production.1}/\reff{def.Psi}/\reff{def.Phi.AZ.first}
of the $A$-series, $Z$-series, $\Phi$-series and $\Psi$-series
of an exponential Riordan array, we conclude:

\begin{corollary}
   \label{cor.dumont3}
The exponential Riordan array $\sfT(y,z)$ has
\be
   A(s)  \;=\;  {e^{zs} \over 1-ys} \,,\quad
   Z(s)  \;=\;  {y e^{zs} \over (1-ys)^2}
 \label{eq.cor.dumont3.1}
\ee
and
\be
   \Phi(s)  \;=\;  e^{zs} \,,\quad
   \Psi(s)  \;=\;  {1 \over 1-ys}
   \;.
 \label{eq.cor.dumont3.2}
\ee
\end{corollary}

The proof of Proposition~\ref{prop.dumont}
follows the elegant argument of Jiang Zeng
that was presented in \cite[section~7]{Dumont_96},
and extends it in part~(b) to handle $k \neq 0$:

\proofof{Proposition~\ref{prop.dumont}}
(a) Consider a tree $T \in \scrt^{[1]}_{n+1}$,
and suppose that the root vertex~1 has $k$ ($\ge 0$) children.
All $k$ edges emanating from the root vertex are proper
and thus get a weight $z$ each.
Deleting these edges and the vertex~1,
one obtains an {\em unordered}\/ partition
of $\{2,\ldots,n+1\}$ into blocks $B_1,\ldots,B_k$
and a rooted tree $T_j$ on each block $B_j$.
Standard enumerative arguments then yield the relation~(a)
for the exponential generating functions.

(b) Consider a tree $T \in \scrt^{\<1;k\>}_{n+1}$ with root $r$,
and let $r_1, \ldots, r_{l+1}$ ($l \ge 0$) be the path in $T$
from the root $r_1 = r$ to the vertex~$r_{l+1} = 1$.\footnote{
   Here $l = 0$ corresponds to the case in~which the vertex~1 is the root.
}
All $l$ edges of this path are improper,
and all $k$ edges from the vertex~1 to its children are proper (and unweighted).
Deleting these edges and the vertex~1,
one obtains a partition of $\{2,\ldots,n+1\}$
into an {\em ordered}\/ collection of blocks $B_1,\ldots,B_l$
and an {\em unordered}\/ collection of blocks $B'_1,\ldots,B'_k$,
together with a rooted tree on each block.
Standard enumerative arguments then yield the relation~(b)
for the exponential generating functions.

(c) In a tree $T \in \scrt^\bullet_n$, focus on the vertex~1
(which might be the root, a leaf, both or neither).
Let $T'$ be the subtree rooted at~1,
and let $T''$ be the tree obtained from $T$
by deleting all the vertices of $T'$ except the vertex~1
(it thus has the vertex~1 as a leaf).
The vertex set $[n]$ is then partitioned as $\{1\} \cup V' \cup V''$,
where $\{1\} \cup V'$ is the vertex set of $T'$
and $\{1\} \cup V''$ is the vertex set of $T''$;
and $T$ is obtained by joining $T'$ and $T''$ at the common vertex~1.
Standard enumerative arguments then yield the relation~(c)
for the exponential generating functions.
\qed

\medskip

{\bf Remarks.}
1. Dumont and Ramamonjisoa also gave \cite[sections~2--5]{Dumont_96}
a second (and very interesting) proof of the $k=0$ case of
Proposition~\ref{prop.dumont},
based on a context-free grammar \cite{Chen_93}
and its associated differential operator.

2. We leave it as an open problem to find a direct combinatorial proof of the
functional equation \reff{eq.cor.dumont.1},
without using the differential equation of Proposition~\ref{prop.dumont}(d).

3. The polynomials $R_n(y,z)$ enumerate rooted trees
according to the number of improper and proper edges;
they are homogenized versions of the celebrated
\textbfit{Ramanujan polynomials}
\cite{Shor_95,Dumont_96,Zeng_99,Guo_07,Lin_14,Josuat-Verges_15,Guo_18,Chen_21,Randazzo_21,forests_totalpos}
\cite[A054589]{OEIS}.

4. The polynomials $R_n$ and $A_{n,0}$ also arise \cite{Josuat-Verges_15}
as derivative polynomials for the tree function:
in the notation of \cite{Josuat-Verges_15} we have $R_n(y,1) = G_n(y-1)$
and $A_{n,0}(y,1) = y \, F_n(y-1)$ for $n \ge 1$.
The formula \reff{eq.cor.dumont.2} is then equivalent to
\cite[Theorem~4.2, equation for~$G_n$]{Josuat-Verges_15}.
\myendremark

\subsection[The matrix $\sfT(y,\bphi)$]{The matrix $\bm{\sfT(y,\bphi)}$}
    \label{subsec.EGF.3}

We now show how Proposition~\ref{prop.dumont} can be generalized
to incorporate the additional indeterminates $\bphi = (\phi_m)_{m \ge 0}$.
We define $\scrt^\bullet_n$, $\scrt^{[i]}_n$ and $\scrt^{\<i;k\>}_n$
as before, and then define the obvious generalizations of
\reff{def.Rn}--\reff{def.An}:
\begin{eqnarray}
   R_n(y,\bphi)
   & = &
   \sum_{T \in \scrt^\bullet_n} y^{\imprope(T)}
      \, \prod_{i=1}^{n+1} \phihat_{\pdeg_T(i)} 
   \qquad\qquad
        \label{def.Rn.phi}  \\[2mm]
   S_n(y,\bphi)
   & = &
   \sum_{T \in \scrt^{[1]}_{n+1}} y^{\imprope(T)}
      \, \prod_{i=1}^{n+1} \phihat_{\pdeg_T(i)} 
        \label{def.Sn.phi}  \\[2mm]
   A_{n,k}(y,\bphi)  \;=\;  t_{n,k}(y,\bphi)
   & = &
   \sum_{T \in \scrt^{\<1;k\>}_{n+1}} y^{\imprope(T)}
      \, \prod_{i=2}^{n+1} \phihat_{\pdeg_T(i)} 
        \label{def.An.phi}
\end{eqnarray}
where $\pdeg_T(i)$ denotes the number of proper children
of the vertex~$i$ in the rooted tree $T$,
and $\phihat_m = m! \, \phi_m$.
(Note that in $R_n$ and $S_n$ we give weights to all the vertices,
 while in $A_{n,k}$ we do {\em not}\/ give any weight
 to the vertex~1.\footnote{
     This differs from the convention used in
     \cite[eq.~(3.24)]{forests_totalpos},
     where $A_n = A_{n,0}$ included a factor $\phi_0 = \phihat_0$
     associated to the leaf vertex~1.
})
We then define the exponential generating functions
\begin{eqnarray}
   \scrr(t;y,\bphi)
   & = &  \sum\limits_{n=1}^\infty R_n(y,\bphi) \: {t^n \over n!}
   \qquad\qquad\qquad
        \label{def.Rn.EGF.phi}  \\[2mm]
   \scrs(t;y,\bphi)  & = &  \sum\limits_{n=0}^\infty S_n(y,\bphi) \: {t^n \over n!}
        \label{def.Sn.EGF.phi}  \\[2mm]
   \scra_k(t;y,\bphi)  & = &  \sum\limits_{n=0}^\infty A_{n,k}(y,\bphi) \: {t^n \over n!}
        \label{def.An.EGF.phi}
\end{eqnarray}
Let us also define the generating function
\be
   \Phi(s)
   \;\eqdef\;
   \sum_{m=0}^\infty \phi_m \, s^m
   \;=\;
   \sum_{m=0}^\infty \phihat_m \, {s^m \over m!}
   \;.
 \label{def.Phi}
\ee
We then have:

\begin{proposition}
   \label{prop.dumont.phi}
The series $\scrr$, $\scrs$ and $\scra_k$
defined in \reff{def.Rn.EGF.phi}--\reff{def.An.EGF.phi}
satisfy the following identities:
\begin{itemize}
   \item[(a)]  $\scrs(t;y,\bphi) \:=\:
                  \Phi\big( \scrr(t;y,\bphi) \bigr)$
   \item[(b)]  $\scra_k(t;y,\bphi) \:=\:
                  \displaystyle{\scrr(t;y,z)^k/k!
                                \over 1 - y \scrr(t;y,\bphi)}$
   \item[(c)]  $\displaystyle {d \over dt} \scrr(t;y,\bphi)  \:=\:
                  \scra_0(t;y,\bphi) \, \scrs(t;y,\bphi)$
\end{itemize}
and hence
\begin{itemize}
   \item[(d)] $\displaystyle {d \over dt} \scrr(t;y,\bphi)
                             \:=\: 
               {\Phi\big( \scrr(t;y,\bphi) \bigr)
                \over
                1 - y \scrr(t;y,\bphi)}$
\end{itemize}
\end{proposition}

\proof
The proof is identical to that of Proposition~\ref{prop.dumont},
with the following modifications:

(a) Consider a tree $T \in \scrt^{[1]}_{n+1}$
in which the root vertex~1 has $k$ children.
Since all $k$ edges emanating from the root vertex are proper,
we get here a factor $\phihat_k/k!$
in place of the $z^k/k!$ that was seen in Proposition~\ref{prop.dumont}.
Therefore, the function $e^{zs}$ in Proposition~\ref{prop.dumont}
is replaced here by the generating function $\Phi(s)$.

(b) No change is needed.

(c) No change is needed.
(The tree $T''$ has vertex~1 as a leaf,
 but in $A_{n,0}$ the vertex~1 is anyway unweighted.)
\qed

Comparing Proposition~\ref{prop.dumont.phi}(b)
with the definition \reff{def.RFG} of exponential Riordan arrays,
we conclude:

\begin{corollary}
   \label{cor.dumont2.phi}
The matrix $\sfT(y,\bphi)$ is the exponential Riordan array $\scrr[F,G]$ where
\be
   F(t)  \;=\;  {1 \over 1 - y \scrr(t;y,\bphi)} \,,\quad
   G(t)  \;=\;  \scrr(t;y,\bphi)
 \label{eq.cor.dumont2.phi.1}
\ee
and $\scrr(t;y,\bphi)$ is the solution of the differential equation
of Proposition~\ref{prop.dumont.phi}(d)
with initial condition $\scrr(0;y,\bphi) = 0$.
\end{corollary}

\noindent
We observe that \reff{eq.cor.dumont2.phi.1}
is identical in form to \reff{eq.cor.dumont2.1};
only $\scrr$ is different.

Comparing Proposition~\ref{prop.dumont.phi}(b,d)
with the definitions
\reff{eq.prop.riordan.exponential.production.1}/\reff{def.Psi}/\reff{def.Phi.AZ.first}
of the $A$-series, $Z$-series, $\Phi$-series and $\Psi$-series
of an exponential Riordan array, we conclude:

\begin{corollary}
   \label{cor.dumont3.phi}
The exponential Riordan array $\sfT(y,\bphi)$ has
\be
   A(s)  \;=\;  {\Phi(s) \over 1-ys} \,,\quad
   Z(s)  \;=\;  {y \, \Phi(s) \over (1-ys)^2}
 \label{eq.cor.dumont3.phi.1}
\ee
and
\be
   \Psi(s)  \;=\;  {1 \over 1-ys}
 \label{eq.cor.dumont3.phi.2}
\ee
where $\Phi(s)$ is given by \reff{def.Phi}.
\end{corollary}

We see that $\Psi(s)$ is the same here as in \reff{eq.cor.dumont3.2};
only $\Phi$ is different.
Proposition~\ref{prop.dumont.phi}
and Corollaries~\ref{cor.dumont2.phi}--\ref{cor.dumont3.phi}
reduce to
Proposition~\ref{prop.dumont}
and Corollaries~\ref{cor.dumont2}--\ref{cor.dumont3}
if we take $\phi_m = z^m/m!$ and hence $\phihat_m = z^m$, $\Phi(s) = e^{zs}$.

\section{Proof of Theorems~\ref{thm1.1}, \ref{thm1.2}, \ref{thm1.4} and \ref{thm1.7}}
   \label{sec.proofs}

In this section we will prove
Theorems~\ref{thm1.1}, \ref{thm1.2}, \ref{thm1.4} and \ref{thm1.7}.
The proofs are now very easy:
we combine the general theory of total positivity
in exponential Riordan arrays developed in Section~\ref{sec.prelim}
(culminating in Corollary~\ref{cor.EAZ.hankel})
with the specific computations of $\Phi$- and $\Psi$-series
carried out in Section~\ref{sec.EGF}.

It suffices of course to prove Theorem~\ref{thm1.7},
since Theorems~\ref{thm1.1}, \ref{thm1.2} and \ref{thm1.4}
are contained in it as special cases:
take $\phi_m = z^m/m!$ to get Theorem~\ref{thm1.4};
then take $y=z=1$ to get Theorems~\ref{thm1.1} and \ref{thm1.2}.
However, we find it instructive to work our way up,
starting with Theorems~\ref{thm1.1} and \ref{thm1.2}
and then gradually adding extra parameters.

\subsection[The matrix $\sfT$]{The matrix $\bm{\sfT}$}

\proofof{Theorems~\ref{thm1.1} and \ref{thm1.2}}
In order to employ the theory of exponential Riordan arrays,
we work here in the ring $\Q$,
even though the matrix elements actually lie in $\Z$.

By Corollary~\ref{cor.prop.EGF.1},
the exponential Riordan array $\sfT$ has
$\Phi(s) = e^s$ and $\Psi(s) = 1/(1-s)$.
By Lemma~\ref{lemma.toeplitz.1a},
the corresponding sequences $\bphi$ and $\bpsi$
(namely, $\phi_m = 1/m!$ and $\psi_m = 1$)
are Toeplitz-totally positive in $\Q$.
Corollary~\ref{cor.EAZ.hankel} then yields
Theorems~\ref{thm1.1}(a) and \ref{thm1.2}.
Theorem~\ref{thm1.1}(b) is obtained from Theorem~\ref{thm1.2}
by specializing to $x=0$.
\qed

Since this proof employed the production-matrix method
(hidden inside Corollary~\ref{cor.EAZ.hankel}),
it is worth making explicit what the production matrix is:

\begin{proposition}[Production matrix for $\sfT$]
   \label{prop.prodmat.T.0}
The production matrix $P = \sfT^{-1} \Delta \sfT$
is~the unit-lower-Hessenberg matrix
\be
   P  \;=\;  B_1 \, \Delta \, D T_1 D^{-1}
 \label{eq.lemma.prodmat.c1.0}
\ee
where $B_1$ is the binomial matrix [i.e.\ \reff{def.Bx} at $x=1$],
$T_1$ is the lower-triangular matrix of all ones [i.e.\ \reff{def.Tx} at $x=1$],
and $D = \diag\bigl( (n!)_{n \ge 0} \bigr)$.
More generally, we have
\be
   B_\xi^{-1} P \, B_\xi  \;=\;  B_1 \, (\Delta + \xi I) \, D T_1 D^{-1}
   \;.
 \label{eq.lemma.prodmat.c2.0}
\ee
\end{proposition}

\proof
Since $\phi_m = 1/m!$ and $\psi_m = 1$,
Proposition~\ref{prop.EAZ.Phi.first} implies
\be
   P  \;=\;  D T_\infty\big((1/m!)_{m \ge 0}\big) D^{-1}
             \, \Delta \, D T_1 D^{-1}
      \;=\;  B_1 \, \Delta \, D T_1 D^{-1}
   \;,
\ee
and Lemma~\ref{lemma.BxinvEAZBx} implies \reff{eq.lemma.prodmat.c2.0}.
\qed

\medskip

{\bf Remarks.}
1. The zeroth and first columns of the matrix $P$ are identical:
that~is, $p_{n,0} = p_{n,1}$.
This can be seen from Lemma~\ref{lemma.EAZ.zeroth_first} with $c=1$,
by noting either that $t_{n,0} = t_{n,1}$ for $n \ge 1$
or that $\Psi(s) = 1/(1-s)$.
Alternatively, it can be seen directly from \reff{eq.lemma.prodmat.c1.0}:
the zeroth and first columns of the matrix $\Delta \, D T_1 D^{-1}$
are identical (namely, they are both equal to $1/(n+1)!$);
so the zeroth and first columns of $M \, \Delta \, D T_1 D^{-1}$
are identical, for {\em any}\/ row-finite matrix $M$.
(Indeed, this would be the case if $D = \diag(\, (n!)_{n \ge 0})$
 were replaced by {\em any}\/ diagonal matrix $\diag(d_0,d_1,d_2,\ldots)$
 satisfying $d_0 = d_1$.)
We will also see that $p_{n,0} = p_{n,1}$
in the explicit formula \reff{eq.prop.prodmat.T}.

The equality $p_{n,0} = p_{n,1}$ implies,
by Lemma~\ref{lemma.EAZ.zeroth_first}(b)$\iff$(b${}'$),
the factorization
\be
   P  \;=\;  P \Delta^{\rm T} \, ({\bf e}_{00} + \Delta)
 \label{eq.P.e00}
\ee
where ${\bf e}_{00}$ denotes the matrix with an entry 1 in position $(0,0)$
and all other entries zero,
and $P \Delta^{\rm T}$ is the lower-triangular matrix obtained from $P$
by deleting its zeroth column.

2. Closely related to the production matrix $P = B_1 \, \Delta \, D T_1 D^{-1}$
are
\be
   \Phat  \:=\: B_1 \, D T_1 D^{-1} \, \Delta
   \qquad\hbox{and}\qquad
   \Phat'  \:=\: \Delta \, B_1 \, D T_1 D^{-1}
   \;.
\ee
It was shown in \cite[Section~4.1]{forests_totalpos}
that $\Phat$ is the production matrix for the forest matrix
$\sfF = (f_{n,k})_{n,k \ge 0}$
where $f_{n,k} = \binom{n}{k} \, k \, n^{n-k-1}$
counts $k$-component forests of rooted trees on $n$ labeled vertices;
and that $\Phat' = \Delta \Phat \Delta^{\rm T}$
is the production matrix for
$\sfF' = \Delta \sfF \Delta^{\rm T} = (f_{n+1,k+1})_{n,k \ge 0}$.
All three production matrices correspond to the same $A$-series
$A(s) = e^s/(1-s)$, but with different splittings into $\Phi$ and $\Psi$.
\myendremark

\medskip

We have more to say about this production matrix $P$,
but in order to avoid disrupting the flow of the argument
we defer it to Section~\ref{subsec.proofs.4}.

\subsection[The matrix $\sfT(y,z)$]{The matrix $\bm{\sfT(y,z)}$}

\proofof{Theorem~\ref{thm1.4}}
In order to employ the theory of exponential Riordan arrays,
we work here in the ring $\Q[y,z]$,
even though the matrix elements actually lie in $\Z[y,z]$.

By Corollary~\ref{cor.dumont3},
the exponential Riordan array $\sfT(y,z)$ has
$\Phi(s) = e^{zs}$ and $\Psi(s) = 1/(1-ys)$.
By Lemma~\ref{lemma.toeplitz.1a},
the corresponding sequences $\bphi$ and $\bpsi$
(namely, $\phi_m = z^m/m!$ and $\psi_m = y^m$)
are Toeplitz-totally positive in the ring $\Q[y,z]$
equipped with the coefficientwise order.
Corollary~\ref{cor.EAZ.hankel} then yields
Theorem~\ref{thm1.4}.
\qed

Analogously to Proposition~\ref{prop.prodmat.T.0}, we have:

\begin{proposition}[Production matrix for $\sfT(y,z)$]
   \label{prop.prodmat.T.yz}
The production matrix $P(y,z) = \sfT(y,z)^{-1} \Delta \sfT(y,z)$
is the unit-lower-Hessenberg matrix
\be
   P(y,z)  \;=\;  B_z \, \Delta \, D T_y D^{-1}
 \label{eq.lemma.prodmat.yz.1}
\ee
where $B_z$ is the weighted binomial matrix \reff{def.Bx},
$T_y$ is the Toeplitz matrix of powers \reff{def.Tx},
and $D = \diag\bigl( (n!)_{n \ge 0} \bigr)$.
More generally,
\be
   B_\xi^{-1} P(y,z) \, B_\xi  \;=\;  B_z \, (\Delta + \xi I) \, D T_y D^{-1}
   \;.
 \label{eq.lemma.prodmat.yz.2}
\ee
\end{proposition}

\bigskip

{\bf Remarks.}
1. The zeroth and first columns of the matrix $P(y,z)$ satisfy
$p_{n,0} = y p_{n,1}$.
This can be seen from Lemma~\ref{lemma.EAZ.zeroth_first} with $c=y$,
by noting either that $t_{n,0}(y,z) = y t_{n,1}(y,z)$ for $n \ge 1$
(Proposition~\ref{prop.tn01yz})
or that $\Psi(s) = 1/(1-ys)$.
Alternatively, it can be seen directly from \reff{eq.lemma.prodmat.c1.0}:
the zeroth column of the matrix $\Delta \, D T_y D^{-1}$
is $y$ times the first column
(they are, respectively, $y^{n+1}/(n+1)!$ and $y^n/(n+1)!$);
so the zeroth column of $M \, \Delta \, D T_y D^{-1}$
is $y$ times the first column, for {\em any}\/ row-finite matrix $M$.

The equality $p_{n,0} = y p_{n,1}$ implies,
by Lemma~\ref{lemma.EAZ.zeroth_first}(b)$\iff$(b${}'$),
the factorization
\be
   P(y,z)  \;=\;  P(y,z) \, \Delta^{\rm T} \, (y\, {\bf e}_{00} + \Delta)
   \;.
\ee

2.  Closely related to the production matrix
$P(y,z) = B_z \, \Delta \, D T_y D^{-1}$
are
\be
   \Phat(y,z)  \:=\: B_z \, D T_y D^{-1} \, \Delta
   \qquad\hbox{and}\qquad
   \Phat'(y,z)  \:=\: \Delta \, B_z \, D T_y D^{-1}
   \;.
\ee
It was shown in \cite[Section~4.3]{forests_totalpos}
that $\Phat(y,z)$ is the production matrix for
$\sfF(y,z) = \big(f_{n,k}(y,z)\bigr)_{n,k \ge 0}$
where $f_{n,k}(y,z)$ counts $k$-component forests of rooted trees
on the vertex set $[n]$ with a weight $y$ (resp.~$z$)
for each improper (resp.~proper) edge.
Likewise, $\Phat'(y,z) = \Delta \Phat(y,z) \Delta^{\rm T}$
is the production matrix for
$\sfF'(y,z) = \Delta \sfF(y,z) \Delta^{\rm T} =
              \big(f_{n+1,k+1}(y,z)\big)_{n,k \ge 0}$.
All three production matrices correspond to the same $A$-series
$A(s) = e^{zs}/(1-ys)$, but with different splittings into $\Phi$ and $\Psi$.
\myendremark

\subsection[The matrix $\sfT(y,\bphi)$]{The matrix $\bm{\sfT(y,\bphi)}$}
   \label{subsec.Typhi}

The proof is similar to that in the preceding subsections,
but a bit of care is needed to handle the case in which
the ring $R$ does not contain the rationals.

\proofof{Theorem~\ref{thm1.7}}
We start by letting $\bphi = (\phi_m)_{m \ge 0}$ be indeterminates,
and working in the ring $\Q[y,\bphi]$.

By Corollary~\ref{cor.dumont3.phi},
the exponential Riordan array $\sfT(y,\bphi)$ has
$\Phi(s) = \sum_{m=0}^\infty \phi_m s^m$
and $\Psi(s) = 1/(1-ys)$, so $\psi_m = y^m$.
We therefore have $\sfT(y,\bphi) = \scro(P)$
and more generally $\sfT(y,\bphi) B_x = \scro(B_x^{-1} P B_x)$,
where Proposition~\ref{prop.EAZ.Phi.first} and Lemma~\ref{lemma.BxinvEAZBx}
tell us that
\be
   B_x^{-1} P B_x
   \;=\;  
   [ D \, T_\infty(\bphi) \, D^{-1} ] \: (\Delta \,+\, x I) \:
        [ D \, T_\infty(\bpsi) \, D^{-1} ]
   \;.
\ee
We now use the definition \reff{def.sharp} to rewrite this as
\be
   B_x^{-1} P B_x
   \;=\;  
   T_\infty(\bphi)^\sharp  \: (\Delta \,+\, x I) \:  T_\infty(\bpsi)^\sharp
   \;.
\ee
Having done this, the equality $\sfT(y,\bphi) B_x = \scro(B_x^{-1} P B_x)$
is now a valid identity in the ring $\Z[y,\bphi]$.
We can therefore now substitute elements $\bphi$ in any commutative ring $R$
for the indeterminates $\bphi$,
and the identity still holds.

By hypothesis the sequence $\bphi$ is Toeplitz-totally positive
in the ring $R$.
By Lemma~\ref{lemma.toeplitz.1},
the sequence $\bpsi$ is Toeplitz-totally positive in the ring $\Z[y]$
equipped with the coefficientwise order.
By Lemma~\ref{lemma.diagmult.TP},
the matrices $T_\infty(\bphi)^\sharp$ and $T_\infty(\bpsi)^\sharp$
are also totally positive.
Therefore $B_x^{-1} P B_x$ is totally positive in the ring $R[x,y]$
equipped with the coefficientwise order.
Proposition~\ref{prop.method} then yields Theorem~\ref{thm1.7}.
\qed

\begin{proposition}[Production matrix for $\sfT(y,\bphi)$]
   \label{prop.prodmat.T.yphi}
The production matrix $P(y,\bphi) = \sfT(y,\bphi)^{-1} \Delta \sfT(y,\bphi)$
is the unit-lower-Hessenberg matrix
\be
   P(y,\bphi)  \;=\;  T_\infty(\bphi)^\sharp \, \Delta \, T_y^\sharp
 \label{eq.lemma.prodmat.yphi.1}
\ee
where $T_y$ is the Toeplitz matrix of powers \reff{def.Tx},
and ${}^\sharp$ is defined in \reff{def.sharp}.
More generally,
\be
   B_\xi^{-1} P(y,\bphi) \, B_\xi
   \;=\;
   T_\infty(\bphi)^\sharp \, (\Delta + \xi I) \, T_y^\sharp
   \;.
 \label{eq.lemma.prodmat.yphi.2}
\ee
\end{proposition}

\bigskip

{\bf Remark.}
1. The zeroth and first columns of the matrix $P(y,\bphi)$ satisfy
$p_{n,0} = y p_{n,1}$, for exactly the same reasons as were observed
for $P(y,z)$.  This implies the factorization
\be
  P(y,\bphi)  \;=\;  P(y,\bphi) \, \Delta^{\rm T} \, (y\, {\bf e}_{00} + \Delta)
   \;.
\ee

2.  Closely related to the production matrix
$P(y,\bphi) = T_\infty(\bphi)^\sharp \, \Delta \, T_y^\sharp$
are
\be
   \Phat(y,\bphi)  \:=\: T_\infty(\bphi)^\sharp \, T_y^\sharp \, \Delta
   \qquad\hbox{and}\qquad
   \Phat'(y,\bphi)  \:=\: \Delta \, T_\infty(\bphi)^\sharp \, T_y^\sharp
   \;.
\ee
It was shown in \cite[Section~4.4]{forests_totalpos}
that $\Phat(y,\bphi)$ is the production matrix for
$\sfF(y,\bphi) = \big(f_{n,k}(y,\bphi)\bigr)_{n,k \ge 0}$
where $f_{n,k}(y,\bphi)$ counts $k$-component forests of rooted trees
on the vertex set $[n]$ with a weight $y$ for each improper edge
and a weight $\phihat_m \eqdef m! \, \phi_m$
for each vertex with $m$ proper children.
Likewise, $\Phat'(y,\bphi) = \Delta \Phat(y,\bphi) \Delta^{\rm T}$
is the production matrix for
$\sfF'(y,\bphi) = \Delta \sfF(y,\bphi) \Delta^{\rm T} =
              \big(f_{n+1,k+1}(y,\bphi)\big)_{n,k \ge 0}$.
All three production matrices correspond to the same $A$-series
$A(s) = \Phi(s)/(1-ys)$, but with different splittings into $\Phi$ and $\Psi$.
\myendremark

\subsection[More on the production matrix for $\sfT$]{More on the production matrix for $\bm{\sfT}$}
   \label{subsec.proofs.4}

We now wish to say a bit more about the production matrix $P$
for the tree matrix $\sfT$.
We begin by giving an explicit formula:

\begin{proposition}
   \label{prop.prodmat.T}
The production matrix $P = \sfT^{-1} \Delta \sfT$
is the unit-lower-Hessenberg matrix with entries
\begin{subeqnarray}
   p_{n,k}
   & = &
   n \binom{n}{k} S_{n-k}  \;+\: \binom{n+1}{k}
      \slabel{eq.prop.prodmat.T.a}  \\[2mm]
   & = &
   {n! \over k! \, (n-k+1)!} \: (n S_{n-k+1} \,+\, 1)
      \slabel{eq.prop.prodmat.T.b}
 \label{eq.prop.prodmat.T}
\end{subeqnarray}
where $S_m$ denotes the \emph{ordered subset number \cite[A000522]{OEIS}}
\be
   S_m  \;\eqdef\;  \sum_{k=0}^m {m! \over k!}
   \;.
  \label{def.ordsubnum}
\ee
These matrix elements satisfy in particular
$p_{n,0} = p_{n,1} = nS_n + 1$ for all $n \ge 0$.
\end{proposition}

The formula \reff{eq.prop.prodmat.T} has a very easy proof,
based on the theory of exponential Riordan arrays
together with our formulae for $A(s)$ and $Z(s)$;
we begin by giving this proof.
On the other hand, it is also of some interest to see
that this production matrix can be found by ``elementary'' algebraic methods,
without relying on the machinery of exponential Riordan arrays
or on any combinatorial interpretation;
this will be our second proof.

\firstproofof{Proposition~\ref{prop.prodmat.T}}
{}From $A(s) = e^s/(1-s)$ we have
\be
   a_n  \;=\;  \sum_{j=0}^n {1 \over j!}
        \;=\;  {S_n \over n!}
   \;.
\ee
{}From $Z(s) = e^s/(1-s)^2$ we have
\be
   z_n  \;=\; \sum_{j=0}^n {n+1-j \over j!}
        \;=\;  n \sum_{j=0}^n {1 \over j!} \,+\, {1 \over n!}
        \;=\;  {n S_n + 1 \over n!}
   \;.
\ee
Theorem~\ref{thm.riordan.exponential.production} and \reff{def.EAZ.1} give
\be
   p_{n,k}
   \;=\;
   {n! \over k!} \: (z_{n-k} \,+\, k \, a_{n-k+1})
   \;,
\ee
and a little algebra leads to (\ref{eq.prop.prodmat.T}a,b).
It is then easy to see that $p_{n,0} = p_{n,1} = nS_n + 1$.
\qed

\secondproofof{Proposition~\ref{prop.prodmat.T}}
An Abel inverse relation
\cite[p.~96, unnumbered equation after (3b)]{Riordan_68}
says that the inverse matrix to
$\sfT = (t_{n,k})_{n,k \ge 0} = \big( \binom{n}{k} n^{n-k} \big)_{n,k \ge 0}$ is
\be
   (\sfT^{-1})_{n,k}
   \;=\;
   (-1)^{n-k} \, \binom{n}{k} \, n \, k^{n-k-1}
   \;.
\ee
It follows that $P = \sfT^{-1} \Delta \sfT$ has matrix elements
\be
   p_{n,k}
   \;=\;
   \sum_{j=k-1}^n (-1)^{n-j} \, \binom{n}{j} \, n \, j^{n-j-1}  \;
                  \binom{j+1}{k} \, (j+1)^{j+1-k}
   \;.
\ee
Setting $N = n-k+1$ and $j = k-1+\ell$ gives
\be
   p_{n,k}
   \;=\;
   \sum_{\ell=0}^N (-1)^{N-\ell} \, \binom{n}{k-1+\ell} \,
                                     n \, (k-1+\ell)^{N-\ell-1}  \;
                  \binom{k+\ell}{k} \, (k+\ell)^{\ell}
   \;,
\ee
which after a bit of playing with the binomial coefficients gives
\be
   p_{n,k}
   \;=\;
   - \: {n \,\cdot\, n! \over k! \, (n-k+1)!} \:
   \sum_{\ell=0}^N \binom{N}{\ell} \, (1-k-\ell)^{N-\ell-1} \, (k+\ell)^{\ell+1}
   \;.
\ee
We now use the Abel identity \cite[p.~22, eq.~(27)]{Riordan_68}
\be
   \sum_{\ell=0}^N \binom{N}{\ell} \, (x+\ell)^{\ell+1} \, (y+N-\ell)^{N-\ell-1}
   \;=\;
   y^{-1} \sum_{\ell=0}^N \binom{N}{\ell} \, \ell! \, (x+\ell) \, (x+y+N)^{N-\ell}
\ee
with $x = k$ and $y = 1-N-k = -n$:  this gives
\begin{subeqnarray}
   p_{n,k}
   & = &
   - \: {n \,\cdot\, n! \over k! \, (n-k+1)!} \: (-1/n) \:
        \sum_{\ell=0}^N \binom{N}{\ell} \, \ell! \, (k+\ell) \, 1^{N-\ell}
               \\[2mm]
   & = &
    {n! \over k! \, (n-k+1)!} \:
        \sum_{\ell=0}^N {N! \over (N-\ell)!} \, (k+\ell)
               \\[2mm]
   & = &
    {n! \over k! \, (n-k+1)!} \:
        \sum_{m=0}^N {N! \over m!} \, (k+N-m)
               \\[2mm]
   & = &
    {n! \over k! \, (n-k+1)!} \:
     \biggl[ (n+1) S_{n-k+1} \:-\: \sum_{m=0}^{n-k+1} {(n-k+1)! \over (m-1)!}
     \biggr]
               \\[2mm]
   & = &
    {n! \over k! \, (n-k+1)!} \: [n S_{n-k+1} \,+\, 1]
   \;,
\end{subeqnarray}
which is \reff{eq.prop.prodmat.T.b}.
\qed

\bigskip

{\bf Remarks.}
1.  The first few rows of this production matrix are
\be
   P
   \;=\;
\Scale[0.95]{
   \begin{bmatrix*}[r]
   1  &   1  &      &      &      &      &      &       \\
   3  &   3  &   1  &      &      &      &      &       \\
   11  &   11  &   5  &   1  &      &      &      &       \\
   49  &   49  &   24  &   7  &   1  &      &      &       \\
   261  &   261  &   130  &   42  &   9  &   1  &      &       \\
   1631  &   1631  &   815  &   270  &   65  &   11  &   1  &       \\
   11743  &   11743  &   5871  &   1955  &   485  &   93  &   13  &   1   \\
   95901  &   95901  &   47950  &   15981  &   3990  &   791  &   126  &   15 & \ddots \\
 \vdots & \vdots & \vdots & \vdots & \vdots & \vdots & \vdots & \vdots &  \ddots
    \end{bmatrix*}
}
    \,.
 \label{eq.prodmat}
\ee
This matrix $P$
(or its lower-triangular variant $P \Delta^{\rm T}$
 in~which the zeroth column is deleted)
 is not currently in \cite{OEIS}.
However, the zeroth and first columns are \cite[A001339]{OEIS},
and the second column $p_{n,2} = nS_n/2$ is \cite[A036919]{OEIS}.

2.  As mentioned earlier, it is not an accident that $p_{n,0} = p_{n,1}$:
by Lemma~\ref{lemma.EAZ.zeroth_first}
this reflects the fact that $\Psi(s) = 1/(1-s)$,
or equivalently that $t_{n,0} = t_{n,1}$.
For the same reason, the production matrices $P(y,z)$ and $P(y,\bphi)$
satisfy $p_{n,0} = y p_{n,1}$.

3.  The ordered subset numbers satisfy the recurrence $S_m = m S_{m-1} + 1$.
\myendremark

Let us now state some further properties of the matrix elements $p_{n,k}$:

\begin{proposition}
   \label{prop.prodmat.T.2}
Define the matrix $P = (p_{n,k})_{n,k \ge 0}$ by
\reff{eq.prop.prodmat.T}/\reff{def.ordsubnum}.
Then:
\begin{itemize}
   \item[(a)]  The $p_{n,k}$ are nonnegative integers that satisfy
the backward recurrence
\be
   p_{n,k}  \;=\;  (k+1) p_{n,k+1} \:+\: \binom{n}{k-1}
 \label{eq.pnk.recurrence}
\ee
with initial condition $p_{n,n+1} = 1$.
   \item[(b)]  The $p_{n,k}$ are also given by
\be
   p_{n,k}  \;=\;  {n S_n \,-\, Q_k(n)  \over  k!}
   \;,
 \label{eq.pnk.Qk}
\ee
where
\begin{subeqnarray}
   Q_k(n)
   & = &
   -1 \,+\, \sum_{j=2}^k (j-1)! \: \binom{n}{j-2}
          \\
   & = &
   -1 \,+\, \sum_{j=2}^k (j-1) \, n^{\underline{j-2}}
       \label{eq.Qkn}
       \slabel{eq.Qkn.b}
\end{subeqnarray}
are polynomials in $n$ with integer coefficients.
In particular, $Q_0(n) = Q_1(n) = -1$ and $Q_2(n) = 0$,
so that $p_{n,0} = p_{n,1} = nS_n + 1$ and $p_{n,2} = n S_n/2$.
\end{itemize}
\end{proposition}

\proof
(a) It is immediate from \reff{eq.prop.prodmat.T.a}/\reff{def.ordsubnum}
that the $p_{n,k}$ are nonnegative integers.
And it is easy to verify,
using the recurrence $S_m = m S_{m-1} + 1$,
that the quantities \reff{eq.prop.prodmat.T}
indeed satisfy the recurrence \reff{eq.pnk.recurrence}.

(b) Introducing the Ansatz \reff{eq.pnk.Qk},
a simple computation shows that
the recurrence \reff{eq.pnk.recurrence} for $p_{n,k}$
is equivalent to the recurrence
\be
   Q_{k+1}(n)  \;=\;  Q_k(n) \,+\, k! \, \binom{n}{k-1}
 \label{eq.Qk.recurrence}
\ee
for $Q_k(n)$.
Furthermore, simple computations show that
$p_{n,0} = p_{n,1} = n S_n + 1$,
so that $Q_0(n) = Q_1(n) = -1$.
It is then easy to see that the unique solution of the recurrence
\reff{eq.Qk.recurrence} with initial condition $Q_0(n) = -1$ is \reff{eq.Qkn}.
%
\qed


{\bf Remarks.}
1.  The first few polynomials $Q_k(n)$ are
\begin{subeqnarray}
   Q_0(n)  & = &  -1  \\[1mm]
   Q_1(n)  & = &  -1  \\[1mm]
   Q_2(n)  & = &  0   \\[1mm]
   Q_3(n)  & = &  2n  \\[1mm]
   Q_4(n)  & = &  3n^2 - n  \\[1mm]
   Q_5(n)  & = &  4n^3 - 9n^2 + 7n  \\[1mm]
   Q_6(n)  & = &  5n^4 - 26n^3 + 46n^2 - 23n  \\[1mm]
   Q_7(n)  & = &  6n^5 - 55n^4 + 184n^3 - 254n^2 + 121n
\end{subeqnarray}
This triangular array is apparently not in \cite{OEIS}.
In any case it follows immediately from \reff{eq.Qkn.b}
that for $k \ge 3$ the leading term in $Q_k(n)$ is $(k-1) n^{k-2}$.
And it also follows from \reff{eq.Qkn.b}
that for $k \ge 4$ the next-to-leading term in $Q_k(n)$ is
$-[(k-2)(k^2 - 4k + 1)/2] \, n^{k-3}$ \cite[A154560]{OEIS}.

2.  Before we found either of the two proofs
of Proposition~\ref{prop.prodmat.T},
we initially {\em guessed}\/ the formulae \reff{eq.prop.prodmat.T}
for $p_{n,k}$, as follows:
Comparison of successive columns of \reff{eq.prodmat}
suggested the backwards recurrence \reff{eq.pnk.recurrence}
for each row of \reff{eq.prodmat}, with initial condition $p_{n,n+1} = 1$.
On the other hand, by looking at the diagonals ($n-k =$ constant)
successively for $n-k = -1,0,1,2,\ldots\,$,
a little experimentation led to the formula \reff{eq.prop.prodmat.T.b}.

3.  The factorization \reff{eq.P.e00} implies that
the unit-lower-Hessenberg matrix $P$ is totally positive
if and only if
the unit-lower-triangular matrix $\widetilde{P} \eqdef P \Delta^{\rm T}$,
obtained from $P$ by deleting its zeroth column,
is totally positive.
Now, the production matrix of $\widetilde{P}$
--- namely, the unit-lower-Hessenberg matrix
$Q = \widetilde{P}^{-1} \Delta \widetilde{P}$ ---
appears to have a very simple form:
\be
   Q  \;\eqdef\;  \widetilde{P}^{-1} \Delta \widetilde{P}
   \;=\;
\Scale[0.95]{
   \begin{bmatrix*}[r]
    3  &   1  &      &      &      &            &    \\
    2  &   2  &   1  &      &      &            &    \\
    6  &   3  &   2  &   1  &      &            &    \\
    24  &   12  &   4  &   2  &   1  &            &    \\
    120  &   60  &   20  &   5  &   2  &   1        &    \\
    720  &   360  &   120  &   30  &   6  &   2     & \ddots    \\
 \vdots & \vdots & \vdots & \vdots & \vdots & \vdots &  \ddots
   \end{bmatrix*}
}
\ee
or in other words
\begin{subeqnarray}
   q_{n,k}  & = & {(n+1)! \over (k+1)!} \quad\hbox{for $k < n$} \\[1mm]
   q_{0,0}  & = & 3 \\[1mm]
   q_{n,n}  & = & 2 \quad\hbox{for $n \ge 1$}  \\[1mm]
   q_{n,n+1}  & = & 1 \\[1mm]
   q_{n,k}  & = & 0 \quad\hbox{for $k > n+1$}
\end{subeqnarray}
(We have not proven this formula for $Q$,
 but it is probably not too difficult.)
Alas, this matrix $Q$ is not even TP${}_2$
(for instance, $q_{10} q_{21} - q_{11} q_{20} = 2 \cdot 3 - 2 \cdot 6 = -6$),
so we cannot use this method to prove the total positivity of $P$.
Nor does it help to subtract of the identity matrix from $Q$
(which would correspond to factoring out a binomial matrix from $\widetilde{P}$
 on the left):
there does not exist any $c \in \R$
for which the leading $3 \times 3$ principal submatrix of $Q - cI$ is TP${}_2$.
\myendremark

\section*{Acknowledgments}

We wish to thank Tomack Gilmore for many helpful conversations.

This research was supported in part by
Engineering and Physical Sciences Research Council grant EP/N025636/1
and by National Natural Science Foundation of China grant 12271078.

\appendix
\section{Interpretation of $\bm{t_{n,k}(y,z)}$ in our first combinatorial model}
   \label{app.first.tnkyz}

In this appendix we give an interpretation
of the polynomials $t_{n,k}(y,z)$, which were defined in \reff{def.tnkyz},
in our first combinatorial model
(rooted trees in~which the root has $k$ lower-numbered children).
However, in order to make this interpretation most natural,
we modify the model slightly, by now considering
rooted trees in~which the root has $k$ {\em higher}\/-numbered children
(this is of course equivalent by reversing all the labels).
We~denote by $\scrt^\bullet_{n+1,k}$ the set of rooted trees
on the vertex set $[n+1]$
in~which exactly $k$ children of the root are higher-numbered than the root.

We will therefore be defining a bijection between two models
on the vertex set $[n+1]$:
\begin{enumerate}[\qquad\bf Model 1.]
  \item Rooted trees in which the root has $k$ higher-numbered children.
  \item Rooted trees in which the vertex~$1$ has $k$ children.
\end{enumerate}
We begin with some definitions.

Let $T$ be a tree on a totally ordered vertex set
(for us it will be $[n+1]$),
and let $e = ij$ be an edge of $T$,
where $i$ is the parent and $j$ is the child.
We say that the edge $e = ij$ is \textbfit{increasing} if $i < j$,
and \textbfit{decreasing} if $i > j$.
We recall that the edge $e = ij$ is \textbfit{improper}
if there exists a descendant of $j$ (possibly $j$ itself)
that~is lower-numbered than $i$;
otherwise it is \textbfit{proper}.
Clearly, every decreasing edge is necessarily improper;
an increasing edge can be either proper or improper,
depending on the behavior of the descendants of $j$.

We now classify edges in a tree $T \in \scrt^\bullet_{n+1,k}$
(that is, Model~1)
as either \textbfit{regular} or \textbfit{irregular}, as follows:

\begin{definition}
   \label{def-irregular}
Let $e = ij$ be an edge in a tree $T \in \scrt^\bullet_{n+1,k}$,
where $i$ is the parent and $j$ is the child.
We classify this edge as follows:
\begin{itemize}
  \item[(I1)] If $ij$ is decreasing, then it is irregular.
  \item[(I2)] If $ij$ is increasing and improper, and $i$ is not the root,
      then $ij$ is irregular.
  \item[(I3)] If $ij$ is increasing and $i$ is the root, then $ij$ is regular.
  	  (That is, the $k$ increasing edges emanating from the root
           are all regular.)
  \item[(I4)] Suppose that all the children of vertex~$1$ are higher-numbered
          than the root.
  	  If~$i=1$ and there is a descendant of $j$
          that is lower-numbered than the root, then $ij$ is irregular.
          (Note that in this case the root cannot be vertex~$1$;
           so this rule does not contradict rule (I3).)
  \item[(I5)] Suppose that vertex~$1$ has at least one child that is
     lower-numbered than the root $\rho$.
     (Note that this implies $\rho \neq 1$.)
     Let $T_1$ be the maximal increasing subtree of $T$ rooted at vertex~$1$,
     whose vertices are $1 = v_1<\cdots<v_{\ell+1}<v_{\ell+2}<\cdots<v_m$,
     where $v_{\ell+1} < \rho < v_{\ell+2}$ (of course $\rho \notin T_1$).
     Then:
 \begin{itemize}
   \item[(a)] all the edges on the path from vertex~$1$ to $v_{\ell+1}$
     are irregular; and
   \item[(b)] an edge $ij\in T_1$ with parent $i=v_s$ and child $j=v_t$ ($s<t$)
    is irregular in~case one of the following is satisfied:
 \begin{itemize}
   \item[(b1)] $\ell+2 \le s < t$ and there is a descendant of $v_{t}$
       in $T$ that is $< v_{s}$;
   \item[(b2)] $s \le \ell < \ell+2 \le t$ and there is a descendant of $v_{t}$
       in $T$ that is $< v_{s+1}$;
   \item[(b3)] $s = \ell+1 < \ell+2 \le t$ and there is a descendant of $v_{t}$
       in $T$ that is $< \rho$;
    \item[(b4)] $s < t \le \ell$ and there is a descendant $v_\tau$ of $v_t$
          in $T_1$ such that $v_{\tau+1} < \rho$
          and there is a descendant of $v_{\tau+1}$ in $T$ that is $<v_{s+1}$;
    \item[(b5)] $s < t \le \ell$ and there is a descendant $v_\tau$ of $v_t$
          in $T_1$ that is $> \rho$,
          and a descendant of $v_\tau$ in $T$ that is $<v_{s+1}$.
 \end{itemize}
 \end{itemize}
  \item[(I6)] All other edges are regular.
\end{itemize}
\end{definition}

\noindent
(We apologize for the complexity of this definition;
but these are the cases that seem to be needed.)

We recall that the polynomials $t_{n,k}(y,z)$
enumerate trees in Model~2 with a weight $y$ (resp.~$z$)
for each improper (resp.~proper) edge,
except that the $k$~proper edges emanating from vertex~1 are unweighted.
We now assert --- and this is the main result of this appendix ---
that the same polynomials $t_{n,k}(y,z)$
enumerate trees in Model~1 with a weight $y$ (resp.~$z$)
for each irregular (resp.~regular) edge,
except that the $k$~regular edges emanating from the root are unweighted:

\begin{proposition}
   \label{prop.tnkyz.model1}
The polynomials $t_{n,k}(y,z)$ defined in \reff{def.tnkyz} satisfy
\be
   t_{n,k}(y,z)
   \;=\;
   \sum_{T \in \scrt^\bullet_{n+1,k}} y^{\irreg(T)} z^{\reg(T) -k}
   \;.
 \label{def.tnkyz.model1}
\ee
\end{proposition}

To prove Proposition~\ref{prop.tnkyz.model1},
we will construct, for each fixed $n$ and $k$,
a bijection $\sigma$ from Model~2 (namely, the set $\scrt^{\<1;k\>}_{n+1}$)
to Model~1 (namely, the set $\scrt^\bullet_{n+1,k}$),
with the property that the number of proper (resp.~improper) edges in $T$
equals the number of regular (resp.~irregular) edges in $\sigma(T)$.
Moreover, we will be able to say which edge in $T$
corresponds to which edge in $\sigma(T)$:
that is, for each $T \in \scrt^{\<1;k\>}_{n+1}$
we will construct a bijection $\psi_T \colon\, E(T) \to E(\sigma(T))$
such that $e \in E(T)$ is proper (resp.~improper)
if and only if $\psi_T(e) \in E(\sigma(T))$ is regular (resp.~irregular).
We summarize this as follows:

\begin{proposition}
   \label{prop.bijection.model21}
There are bijections $(\sigma,\psi_T)$ from Model~2 to Model~1
that map proper (resp.~improper) edges in Model~2
to regular (resp.~irregular) edges in Model~1.
\end{proposition}

\noindent
The remainder of this appendix is devoted to proving
Proposition~\ref{prop.bijection.model21}.

\bigskip

Given a tree $T$ rooted at $r$ in Model~2,
let $v_1<v_2<\cdots<v_k$ be the $k$~children of the vertex~$1$.
If vertex~$1$ is the root,
then its $k$~children are obviously higher-numbered than the root.
In this situation we define $\sigma(T)=T$, which also belongs to Model~1;
and we define $\psi_T(e) = e$ for all $e \in E(T)$.

Now suppose that the vertex~$1$ is not the root.
First consider the case $k=0$.
Since we will use this special case as a tool in handling the general case,
we here denote the bijection $\phi$ instead of $\sigma$.

\begin{lemma}\label{lem-bijec1&2-k=0}
For $k=0$, there are bijections $(\phi,\psi_T)$ from Model~2 to Model~1
that map proper (resp.~improper) edges in Model~2
to regular (resp.~irregular) edges in Model~1.
\end{lemma}

The following construction is inspired by \cite[proof of Lemma~1]{Chauve_00}.

\begin{figure}[t]
  \begin{center}
  \begin{tabular}{cp{50pt}c}
  \begin{tikzpicture}[baseline=(current bounding box.base)]
    \node at (0,1) {$T$};
    \node (root) at (0,0) {};
    \fill (root) circle(.1) node[above] {$r$};
    \node (L) at (-1,-1) {};
    \node (DL) at (-1,-3) {};
    \node (DLtop) at (-1,-2.5) {};
    \draw (L) ellipse (.8 and .3) node {$L$};
    \draw (DL) ellipse (.8 and 1) node {$D_L$};
    \draw (root)--(L)--(DLtop);
    \node (H) at (1,-1) {};
    \node (DH) at (1,-3) {};
    \node (DHtop) at (1,-2.5) {};
    \draw (H) ellipse (.8 and .3) node {$H$};
    \draw (DH) ellipse (.8 and 1) node {$D_H$};
    \draw (root)--(H)--(DHtop);
  \end{tikzpicture}
  &&
  \begin{tikzpicture}[baseline=(current bounding box.base)]
    \node at (0,1) {$T_{\rm max}$};
    \node (root) at (0,0) {};
    \fill (root) circle(.1) node[above] {$r$};
    \node (H) at (0,-1) {};
    \node (DH') at (0,-3) {};
    \node (DHtop') at (0,-2.5) {};
    \draw (H) ellipse (.8 and .3) node {$H$};
    \draw (DH') ellipse (.5 and 1) node {$D'_H$};
    \draw (root)--(H)--(DHtop');
  \end{tikzpicture}
  \qquad
  \begin{tikzpicture}[baseline=(current bounding box.base)]
    \node at (0,1) {$T_0$};
    \node (root) at (0,0) {};
    \fill (root) circle(.1) node[above] {$r_0=r$};
    \node (L) at (0,-1) {};
    \node (DL) at (0,-3) {};
    \node (DLtop) at (0,-2.5) {};
    \draw (L) ellipse (.8 and .3) node {$L$};
    \draw (DL) ellipse (.8 and 1) node {$D_L$};
    \draw (root)--(L)--(DLtop);
  \end{tikzpicture}
  \qquad
  \begin{tikzpicture}[baseline=(current bounding box.base)]
    \node at (0,1) {$T_1$};
    \node (root) at (0,0) {};
    \fill (root) circle(.1) node[above] {$r_1$};
    \node (L1) at (0,-1) {};
    \node (DL1) at (0,-3) {};
    \node (DL1top) at (0,-2.5) {};
    \draw (L1) ellipse (.5 and .3) node {$L_1$};
    \draw (DL1) ellipse (.5 and .8) node {$D_{L_1}$};
    \draw (root)--(L1)--(DL1top);
  \end{tikzpicture}
  \raisebox{25pt}{$\cdots$}
  \begin{tikzpicture}[baseline=(current bounding box.base)]
    \node at (0,1) {$T_p$};
    \node (root) at (0,0) {};
    \fill (root) circle(.1) node[above] {$r_p$};
    \node (Lp) at (0,-1) {};
    \node (DLp) at (0,-3) {};
    \node (DLptop) at (0,-2.5) {};
    \draw (Lp) ellipse (.5 and .3) node {$L_p$};
    \draw (DLp) ellipse (.5 and .8) node {$D_{L_p}$};
    \draw (root)--(Lp)--(DLptop);
  \end{tikzpicture}
  \\
  (a) && (b)
  \end{tabular}
  \end{center}
  \caption{(a) The tree $T$ in Model~2 when $r \neq 1$. (b) Subtrees of $T$.}
  \label{fig-tree-model2}
\end{figure}
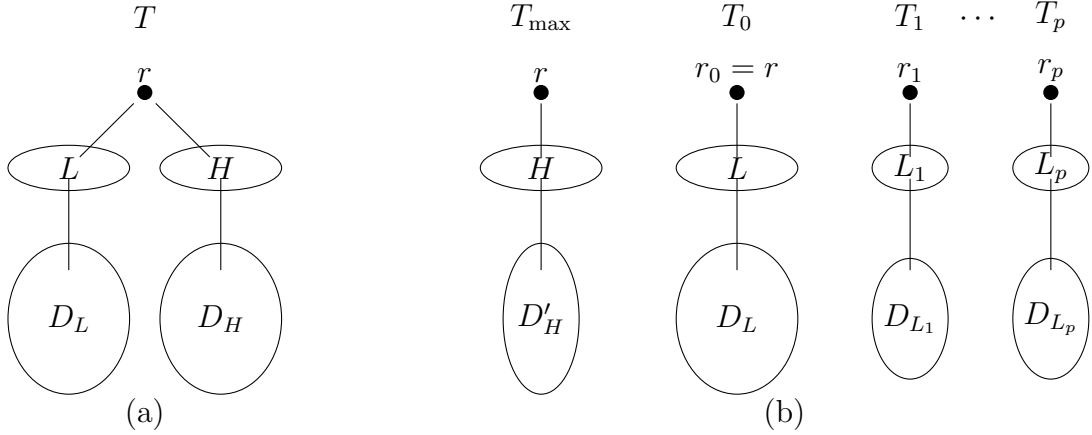

\proofof{Lemma~\ref{lem-bijec1&2-k=0}}
Let $T$ be a tree in Model~2,
in which $r$ is the root and vertex~$1$ is a leaf.
Since we have already handled the case $r=1$,
we assume henceforth that $r \neq 1$.
Let $L$ (resp.~$H$) denote the set of lower- (resp.~higher-) numbered children
of the root;
and let $D_L$ (resp.~$D_H$) denote the set of all descendants
of the vertices in $L$ (resp.~$H$),
excluding those in $L$ (resp.~$H$) itself.
See Figure~\ref{fig-tree-model2}(a).

Let $T_{\rm max}$ be the maximal increasing subtree of $T$ rooted at $r$,
and let $T_0, \ldots, T_p$ be the trees obtained from $T$
by deleting all the edges in $T_{\rm max}$.
Let $r_j$ be the root of $T_j$ for $0\le j\le p$.
In particular, we choose $r_0$ to be the root $r$.
Note that each $r_j$ is a vertex in $T_{\rm max}$
(otherwise it would not become a root
 when we delete the edges in $T_{\rm max}$);
and conversely, every vertex in $T_{\rm max}$ becomes a root $r_j$
(though its tree $T_j$ might be trivial).
Therefore, all of the higher-numbered children of $r_j$ belong to $T_{\rm max}$,
while all of the lower-numbered children of $r_j$ belong to $T_j$.
We denote the set of those lower-numbered children by $L_j$.
Of course, $L_0 = L$.
See Figure~\ref{fig-tree-model2}(b).
Furthermore, since $T_{\rm max}$ is an increasing tree,
it is rooted at its smallest label (namely, $r$);
therefore $r_j > r$ for $1 \le j \le p$.

{\bf Notation:} Let $S$ be a sequence of increasing numbers.
For $a\not\in S$ and $b \in S$,
define $\relab{a}{b}$ as an operator acting on $S$
such that $\relab{a}{b}(S):=S\cup\{a\}\backslash\{b\}$
is still increasingly ordered.
For example, $\relab{2}{5} (1,3,5,7)=(1,2,3,7)$.
We observe that the inverse of $\relab{a}{b}$ is $\relab{b}{a}$.
Further, if $T$ is a tree whose vertex set is $S$,
we write $\relab{a}{b}(T)$ to denote the tree
with vertex set $\relab{a}{b}(S)$ that is obtained from $T$
by relabeling the vertices according to the map $\relab{a}{b}$.

\begin{figure}[t]
  \begin{center}
  \begin{tabular}{cp{50pt}c}
  \begin{tikzpicture}[baseline=(current bounding box.base)]
    \node at (0,1) {$\phi(T)$};
    \node (root) at (0,0) {};
    \node (Ti) at (0,-1) {};
    \fill (root) circle(.1) node[above] {$r_i$};
    \draw (Ti) ellipse (.5 and 1) node {$T_i$};
    \node (1) at (0,-2) {};
    \fill (1) circle(.1) node[above] {$1$};
    \node (rhoTmax) at (0,-4) {};
    \draw (rhoTmax) ellipse (1.5 and 2) node {$\relab{1}{r_i}(T_{\rm max})$};
    \node (r) at (-1,-3) {};
    \node (T0) at (-2,-4) {};
    \fill (r) circle(.1) node[above] {$r$};
    \draw[rotate=-45] (T0) ellipse (.5 and 1.4) node {$T_0$};
    \draw (1)--(r);
    \node (r1) at (1,-3.2) {};
    \node (T1) at (1.9,-4) {};
    \fill (r1) circle(.1) node[left] {$r_1$};
    \draw[rotate=45] (T1) ellipse (.5 and 1.2) node {$T_1$};
    \node (rp) at (.5,-5.5) {};
    \node (Tp) at (.5,-6.5) {};
    \fill (rp) circle(.1) node[above] {$r_p$};
    \draw (Tp) ellipse (.5 and 1) node {$T_p$};
    \node[rotate=40] at (2,-6) {$\cdots$};
  \end{tikzpicture}
  &&
  \begin{tikzpicture}[baseline=(current bounding box.base)]
    \node at (0,1) {$\phi(T)$};
    \node (root) at (0,0) {};
    \node (T0) at (0,-1) {};
    \fill (root) circle(.1) node[above] {$r$};
    \draw (T0) ellipse (.5 and 1) node {$T_0$};
    \node (1) at (0,-2) {};
    \fill (1) circle(.1) node[above] {$1$};
    \node (H) at (0,-3) {};
    \node (DH) at (0,-5) {};
    \node (DHtop) at (0,-4.5) {};
    \draw (H) ellipse (.8 and .3) node {$H$};
    \draw (DH) ellipse (.8 and 1) node {$D_H$};
    \draw (1)--(H)--(DHtop);
  \end{tikzpicture}
  \\
  \\
  (a)  && (b)
  \end{tabular}
  \end{center}
  \vspace*{-3mm}
  \caption{The tree $\phi(T)$ in Model~1 when $r \neq 1$.
     (a) General case, but when $i=0$ the vertex~$r$ and subtree~$T_0$
         should be removed.
     (b) Redrawing of the special case $i=0$,
         where $\relab{1}{r_i}(T_{\rm max})$ is the
         entire subtree of $\phi(T)$ rooted at the vertex~$1$.
  }
  \label{fig-tree-model1}
\end{figure}
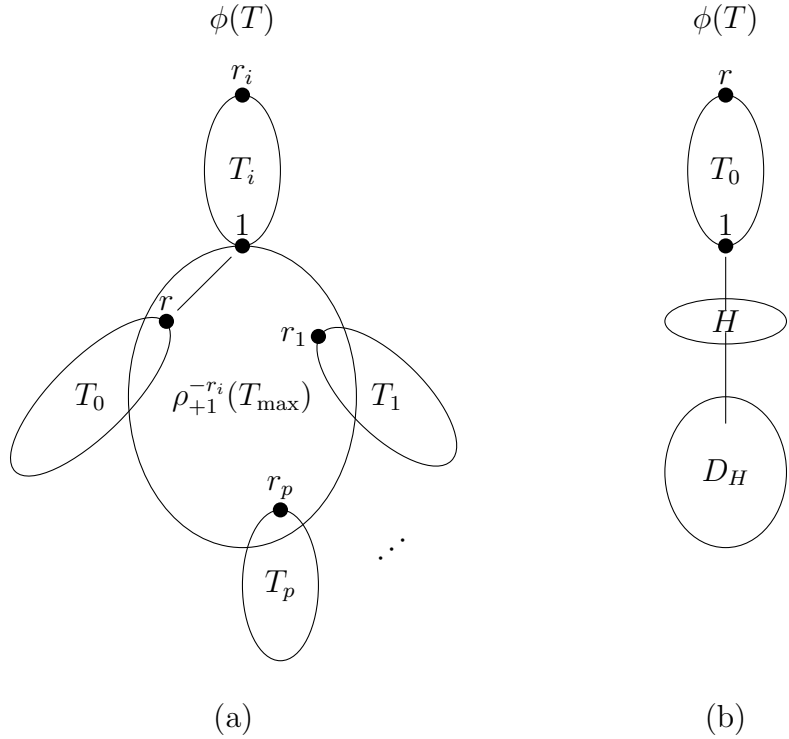

We now make the trivial observation that, since $r > 1$ by assumption,
the vertex~$1$ is not in $T_{\rm max}$.
Let $T_i$ be the tree containing vertex~$1$;
since $k=0$, the vertex~$1$ is a leaf in $T_i$.
The bijection $\phi$ is defined in three steps:

\begin{enumerate}[Step 1.]
  \item Take the tree $T_{\rm max}$, and relabel its vertices to obtain
    $\relab{1}{r_i}(T_{\rm max})$.
    Since $T_{\rm max}$ is an increasing tree,
    it is rooted at its smallest label (namely, $r$);
    therefore, the relabeled tree $\relab{1}{r_i}(T_{\rm max})$
    is rooted at its smallest label, which is the vertex~$1$.
    [If $i=0$, we relabeled $r \to 1$ and left all other labels unaffected.
     If $i \neq 0$, we relabeled $r \to 1$ and relabeled the
     second-smallest label of $T_{\rm max}$
     (that is, the lowest-numbered child of $r$)
     to $r$
     --- among other relabelings, the details of which
     will be worked out below.]
  \item Graft $\relab{1}{r_i}(T_{\rm max})$ onto $T_i$
        by identifying the two vertices $1$;
        call the result $T'_i$.
  \item Graft each tree $T_j~(j\neq i)$ onto $T'_i$ by
        identifying the two vertices $r_j$;
        call the result $\phi(T)$.
        See Figure~\ref{fig-tree-model1}(a).
\end{enumerate}

In this way we obtain a tree $\phi(T)$ rooted at $r_i$,
in which all the children of $r_i$ are lower-numbered,
and in which $\relab{1}{r_i}(T_{\rm max})$
is the maximal increasing subtree of $\phi(T)$ rooted at the vertex~$1$.
Furthermore, if $i \neq 0$,
then the lowest-numbered child of the vertex~$1$ is $r$,
which is smaller than the root~$r_i$ of $\phi(T)$;
while if $i=0$, then, as shown in Figure~\ref{fig-tree-model1}(b),
the children of vertex~$1$ are precisely the set $H$
(which did not undergo any relabeling in Step~1),
which are all larger than the root~$r$ of $\phi(T)$.

These observations allow us to obtain the inverse of $\phi$.
If the smallest-numbered child of vertex~$1$ is smaller than
the root of $\phi(T)$, then that child is $r$,
and we are in the case $i \neq 0$;
otherwise we are in the case $i=0$,
and the root of $\phi(T)$ is $r$.
If we delete the edges of $\relab{1}{r_i}(T_{\rm max})$ from $\phi(T)$,
we recover the trees $T_i$; this undoes Step~3.
We then undo Step~2 by separating the subtree rooted at vertex~$1$.
And finally, we undo Step~1 by relabeling $\relab{1}{r_i}(T_{\rm max})$
using the map $\relab{r_i}{1}$; this yields $T_{\rm max}$.
We can then reassemble the pieces
$T_{\rm max}$ and $T_0,\ldots,T_p$ to obtain $T$.

It is also clear how the map $\psi_T$ is defined,
since each edge of $T$ corresponds, via relabeling and grafting,
to a well-defined edge of $\phi(T)$.

We now look at how the map $\psi_T$ acts on proper and improper edges of $T$.
We first observe that the edges in $T_0,\ldots,T_p$
do not undergo any relabeling;
so an edge $e$ in one of these subtrees
is increasing (resp.~decreasing)
according as the edge $\psi_T(e)$ in $\phi(T)$
is increasing (resp.~decreasing).
Furthermore, the descendants in these subtrees
are the same as the descendants of their images in $\phi(T)$,
except that the vertex~$1$ in $\phi(T)$ has acquired extra descendants,
which are anyway higher-numbered and therefore do not affect
properness or improperness.
Therefore, an edge $e$ in one of these subtrees is proper (resp.~improper)
according as the edge $\psi_T(e)$ in $\phi(T)$ is proper (resp.~improper).
This means, by rules~(I1), (I2) and (I6) of Definition~\ref{def-irregular},
that an edge $e$ in one of these subtrees is proper (resp.~improper)
according as the edge $\psi_T(e)$ in $\phi(T)$ is regular (resp.~irregular).
[Note that rule~(I3) plays no role here, because $k=0$.
Rule~(I4) does not apply because all the edges in $\phi(T)$
emanating from vertex~$1$ lie in $\relab{1}{r_i}(T_{\rm max})$.
And rule~(I5) applies only within $\relab{1}{r_i}(T_{\rm max})$.]

We now need to consider the edges in $\relab{1}{r_i}(T_{\rm max})$.
We divide the proof into two cases:

\medskip

{\bf Case 1: $\bm{i=0}$.}
Let $r=u_1<u_2<\cdots<u_m$ be the vertices in $T_{\rm max}$;
then $\relab{1}{r_i}$ acts as follows:
$$
\relab{1}{r_i} \colon~(r,u_2,\ldots,u_m) \mapsto (1,u_2,\ldots,u_m).
$$
That is, as previously observed,
$\relab{1}{r_i}(T_{\rm max})$ is obtained from $T_{\rm max}$
by only relabeling the root~$r$ as the vertex~$1$;
and $\phi(T)$ is as shown in Figure~\ref{fig-tree-model1}(b).
Therefore, for all edges in $T_{\rm max}$
other than those emanating from the root,
properness/improperness in $T$ corresponds to
properness/improperness of their images in $\phi(T)$;
and by rules~(I1), (I2) and (I6),
this corresponds to regularity/irregularity of the images in $\phi(T)$.
[Rule~(I4) does not apply because the parent is not vertex~$1$;
and rule~(I5) does not apply because all the children of vertex~$1$
are higher-numbered than the root~$r$.]

Now consider an edge $e$ in $T$ that emanates from the root~$r$
to a higher-numbered child $h \in H$.
The bijection $\psi_T$ maps $e$ to an edge $e'$ in $\phi(T)$
that emanates from vertex~$1$ to $h \in H$.
The edge $e$ is improper in~case
there is a descendant of $h$ in $T$ that is $< r$;
and this is equivalent to the existence of a descendant of $h$ in $\phi(T)$
that is $< r$.
Since in this case $r$ is the root of $\phi(T)$,
and all the children of vertex~$1$ in $\phi(T)$ are $> r$,
rule~(I4) of Definition~\ref{def-irregular}
specifies that the edge $e'$ is irregular whenever $e$ is improper;
otherwise, by rule~(I6), it is regular.
[Once again, rule~(I5) does not apply here.]

\medskip

{\bf Case 2: $\bm{i \neq 0}$.}
Let $r=u_1<u_2<\cdots<u_{\ell}<u_{\ell+1}=r_i<u_{\ell+2} < \cdots<u_m$
be the vertices in $T_{\rm max}$;  then $\relab{1}{r_i}$ acts as follows:
\begin{eqnarray*}
\relab{1}{r_i} \colon~(r,u_2,u_3,\ldots,u_{\ell},r_i,u_{\ell+2},\ldots,u_m)
   & \mapsto & (1,r,u_2,\ldots,u_{\ell-1},u_{\ell},u_{\ell+2},\ldots,u_m)
          \\
   & := & (v_1,v_2,v_3,\ldots,v_{\ell},v_{\ell+1},v_{\ell+2},\ldots,v_m)
\end{eqnarray*}
Set $u_0:=1$.
Then $\relab{1}{r_i}(T_{\rm max})$ is obtained from $T_{\rm max}$
by relabeling each vertex $u_s$ that is $\le r_i$ by $u_{s-1}$,
and leaving all vertices $> r_i$ unchanged;
in other words, $v_s = u_{s-1}$ for $s \le \ell+1$
and $v_s = u_s$ for $s \ge \ell+2$.
Therefore, each edge $e = u_s u_t$ in $T_{\rm max} \subseteq T$
maps onto $\psi_T(e) = v_s v_t$
in $\relab{1}{r_i}(T_{\rm max}) \subseteq \phi(T)$;
and the descendants of $u_t$ in $T_{\rm max}$
map via $\psi_T$ onto the descendants of $v_t$
in $\relab{1}{r_i}(T_{\rm max})$.

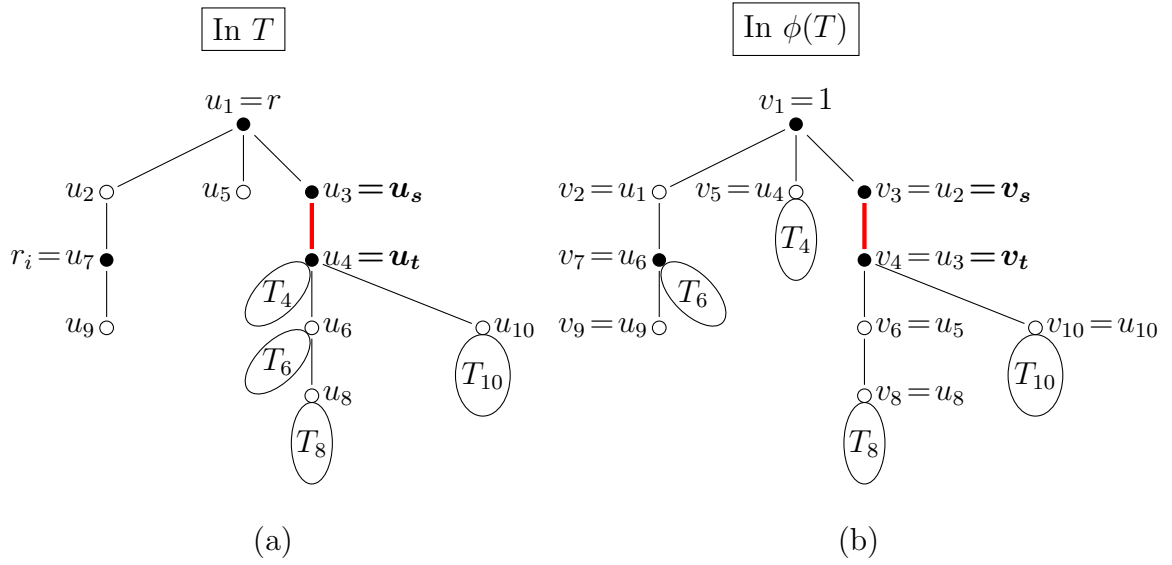
\begin{figure}[tp]
  \hspace*{-4mm}
  \begin{tabular}{cp{0pt}@{\hspace*{-4mm}}c}
  \begin{tikzpicture}[baseline = (current bounding box.base), scale = 0.9]
    \node[draw] at (0,1.4) {In $T$};
    \node (root) at (0,0) {};
    \fill (root) circle(.1) node[above] {$u_1 \!=\! r$};
    \node (u2) at (-2,-1) {};
    \node (u3) at (1,-1) {};
    \node (u4) at (1,-2) {};
    \node (u5) at (0,-1) {};
    \node (u6) at (1,-3) {};
    \node (u7) at (-2,-2) {};
    \node (u8) at (1,-4) {};
    \node (u9) at (-2,-3) {};
    \node (u10) at (3.5,-3) {};
    \foreach \x in {2,5,9}
    \draw (u\x) circle(.1) node[left] {$u_{\x}$};
    \foreach \x in {6,8,10}
    \draw (u\x) circle(.1) node[right] {$u_{\x}$};
    \fill (u3) circle(.1) node[right] {$u_3 \bm{\!=\! u_s}$};
    \fill (u4) circle(.1) node[right] {$u_4 \bm{\!=\! u_t}$};
    \fill (u7) circle(.1) node[left] {$r_i \!=\! u_7$};
    \draw (root)--(u2)--(u7)--(u9) (root)--(u5) (root)--(u3)--(u4)--(u6)--(u8) (u4)--(u10);
    \draw[red,ultra thick] (u3)--(u4);
    \node (T4) at (.5,-2.5) {};
    \node (T6) at (.5,-3.5) {};
    \node (T8) at (1,-4.7) {};
    \node (T10) at (3.5,-3.7) {};
    \draw[rotate = -45] (T4) ellipse (.3 and .6) node {$T_4$};
    \draw[rotate = -45] (T6) ellipse (.3 and .6) node {$T_6$};
    \draw (T8) ellipse (.3 and .6) node {$T_8$};
    \draw (T10) ellipse (.4 and .6) node {$T_{10}$};
  \end{tikzpicture}
  &&
  \begin{tikzpicture}[baseline = (current bounding box.base), scale = 0.9]
    \node[draw] at (0,1.4) {In $\phi(T)$};
    \node (root) at (0,0) {};
    \fill (root) circle(.1) node[above] {$v_1 \!=\! 1$};
    \node (v2) at (-2,-1) {};
    \node (v3) at (1,-1) {};
    \node (v4) at (1,-2) {};
    \node (v5) at (0,-1) {};
    \node (v6) at (1,-3) {};
    \node (v7) at (-2,-2) {};
    \node (v8) at (1,-4) {};
    \node (v9) at (-2,-3) {};
    \node (v10) at (3.5,-3) {};
    \foreach \x in {2,5}
    {
    \pgfmathtruncatemacro{\y}{\x-1}
    \draw (v\x) circle(.1) node[left] {$v_{\x} \!=\! u_{\y}$};
}
    \foreach \x in {8,10}
    \draw (u\x) circle(.1) node[right] {$v_{\x} \!=\! u_{\x}$};
    \fill (v3) circle(.1) node[right] {$v_3 \!=\! u_2 \bm{\!=\! v_s}$};
    \fill (v4) circle(.1) node[right] {$v_4 \!=\! u_3 \bm{\!=\! v_t}$};
    \fill (v7) circle(.1) node[left] {$v_7 \!=\! u_6$};
    \draw (v6) circle(.1) node[right] {$v_6 \!=\! u_5$};
    \draw (v9) circle(.1) node[left] {$v_9 \!=\! u_9$};
    \draw (root)--(u2)--(u7)--(u9) (root)--(u5) (root)--(u3)--(u4)--(u6)--(u8) (u4)--(u10);
    \draw[red,ultra thick] (u3)--(u4);
    \node (T4) at (0,-1.7) {};
    \node (T6) at (-1.5,-2.5) {};
    \node (T8) at (1,-4.7) {};
    \node (T10) at (3.5,-3.7) {};
    \draw (T4) ellipse (.3 and .6) node {$T_4$};
    \draw[rotate = 45] (T6) ellipse (.3 and .6) node {$T_6$};
    \draw (T8) ellipse (.3 and .6) node {$T_8$};
    \draw (T10) ellipse (.4 and .6) node {$T_{10}$};
  \end{tikzpicture}
  \\
  \\
  (a) && (b)
  \end{tabular}
  \caption{
      An example illustrating cases (b4,5)
      in the proof of Lemma~\ref{lem-bijec1&2-k=0},
      showing the trees $T_{\rm max} \subseteq T$
      and $\relab{1}{r_i}(T_{\rm max}) \subseteq \phi(T)$
      along with some of the trees $T_j$ hanging off them
      (namely, those trees attached to the descendants of $u_t$
       are shown).
      The edge $e = u_s u_t$ and its image $\psi_T(e) = v_s v_t$
      are shown in thick red.
      Note that tree $T_j$ is attached at vertex $u_j$,
      which equals $v_j = \psi_T(u_j)$ or $v_{j+1} = \psi_T(u_{j+1})$
      according as $u_j > r_i$ or $u_j < r_i$.
  }
  \label{fig-example-b4b5}
\end{figure}

Note that in $\phi(T)$, vertex~$1$ has at least one child (namely, $r$)
that is lower-numbered than the root $r_i$.
Therefore rule~(I5) applies, with $\rho = r_i$
and $T_1 = \relab{1}{r_i}(T_{\rm max})$:

\begin{itemize}
  \item[(a)] All the edges on the path from the root~$r$ to vertex~$r_i$
     in $T_{\rm max} \subseteq T$ are improper,
     since vertex~$1$ is a descendant of $r_i$ in $T$.
     These edges map, under the relabeling $\relab{1}{r_i}$,
     onto the path from vertex~$1$ to $v_{\ell+1}$ in $\phi(T)$.
     By rule~(I5a) of Definition~\ref{def-irregular},
     all the edges in this path are irregular.
\end{itemize}
The foregoing case needed to be treated separately,
because the vertices in the path from $r$ to $r_i$ in $T_{\rm max} \subseteq T$
have descendants (in particular, the vertex~$1$)
that do {\em not}\/ correspond (via the relabeling)
to descendants in their images in $\phi(T)$,
because the tree $T_i$ was moved from its position in $T$
to the root in $\phi(T)$.
This problem does not arise in the remaining cases:
\begin{itemize}
  \item[(b1)] Consider an edge $e = u_s u_t$ in $T_{\rm max}$,
     where $\ell+2 \le s < t$.
     These vertices do not get relabeled, so $\psi_T(e) = e$.
     This edge is improper in $T$ in~case there is a descendant of $u_t$
     in $T$ that is $< u_s$.
     By rule~(I5b1) of Definition~\ref{def-irregular},
     this edge is irregular in $\phi(T)$ in exactly the same situation.
  \item[(b2,3)] Now consider an edge $e = u_s u_t$ in $T_{\rm max}$,
     where $s \le \ell+1 < \ell+2 \le t$.
     Then vertex $u_s$ gets relabeled to $u_{s-1}$,
     while $u_t$ does not get relabeled;
     so $\psi_T(e) = v_s v_t = u_{s-1} u_t$.
     The edge $e$ is improper in $T$ in~case there is a descendant of $u_t$
     in $T$ that is $< u_s$.
     Now $u_s = v_{s+1}$ in~case $s \le \ell$,
     while $u_s = r_i = \rho$ in~case $s = \ell+1$.
     By rules~(I5b2,3) of Definition~\ref{def-irregular},
     the edge $v_s v_t$ is irregular in $\phi(T)$
     exactly when $e$ is improper in $T$.
\end{itemize}
Note that in~cases (b1--3), the descendants of $u_t$ in $T$
are the same as the descendants of $v_t$ in $\phi(T)$,
because the relevant trees $T_j$ were grafted in the same place
(since their roots $r_j$ did not get relabeled).
Things will be slightly more complicated in the remaining cases:
\begin{itemize}
  \item[(b4,5)] Consider an edge $e = u_s u_t$ in $T_{\rm max}$,
     where $s < t \le \ell$.
     Then vertices $u_s$ and $u_t$ both get relabeled,
     so $\psi_T(e) = v_s v_t = u_{s-1} u_{t-1}$.
     The edge $e$ is improper in $T$ in~case there is a descendant of $u_t$
     in $T$ that is $< u_s$.
     (Note that $u_s = v_{s+1}$ and $u_t = v_{t+1}$ because $s,t \le \ell$.)
     Such a descendant cannot lie in $T_{\rm max}$,
     because $T_{\rm max}$ is increasing,
     but it can lie in one of the trees $T_j$
     that is attached to $T_{\rm max}$.
     So consider all of the descendants $u_\tau$ of $u_t$ in $T_{\rm max}$.
     If one of these descendants is $r_i$,
     then we are in the already-treated case~(a);
     so we can assume that they are all either $< r_i$ or $> r_i$.
     The images of the vertices $u_\tau$ under $\relab{1}{r_i}$
     are the descendants $v_\tau = \relab{1}{r_i}(u_\tau)$ of $v_t$
     in $T_1 = \relab{1}{r_i}(T_{\rm max})$.
     Now consider the two cases $u_\tau < r_i$ and $u_\tau > r_i$
     (recalling that $r_i = \rho$):
\begin{itemize}
   \item[(b4)]  $u_\tau < \rho$ is equivalent to
     $u_\tau = v_{\tau+1} < \rho$.
     The edge $u_s u_t$ is improper in $T$ in~case
     there is a descendant of $u_\tau = v_{\tau+1}$ in $T$
     that is $< u_s = v_{s+1}$;
     and the descendants of  $u_\tau = r_j$ in the tree $T_j \subseteq T$
     are the same as the descendants of $v_{\tau+1} = r_j$
     in the tree $T_j \subseteq \phi(T)$.
     By rule~(I5b4) of Definition~\ref{def-irregular},
     the edge $v_s v_t$ is irregular in $\phi(T)$
     exactly when $e$ is improper in $T$.
   \item[(b5)]  $u_\tau > \rho$ is equivalent to
     $u_\tau = v_\tau > \rho$.
     The edge $u_s u_t$ is improper in $T$ in~case
     there is a descendant of $u_\tau = v_\tau$ in $T$
     that is $< u_s = v_{s+1}$;
     and the descendants of  $u_\tau = r_j$ in the tree $T_j \subseteq T$
     are the same as the descendants of $v_\tau = r_j$
     in the tree $T_j \subseteq \phi(T)$.
     By rule~(I5b5) of Definition~\ref{def-irregular},
     the edge $v_s v_t$ is irregular in $\phi(T)$
     exactly when $e$ is improper in $T$.
\end{itemize}
See Figure~\ref{fig-example-b4b5} for an example illustrating the cases (b4,5).
\end{itemize}
There is, {\em a priori}\/,
one additional case for an edge $e = u_s u_t$ in $T_{\rm max}$,
namely, $s < t = \ell+1$.
But this corresponds to the last edge on the path
from $r$ to $r_i$ in $T_{\rm max}$,
and hence was already treated on case~(a).

We have now considered all the cases in which an edge $e \in T$
can be improper;
so by rule~(I6) of Definition~\ref{def-irregular}, 
this completes the proof.
\qed

{\bf Remark.}  In case (b4,5) one might worry what happens
when $v_\tau < \rho$ while $v_{\tau+1} > \rho$,
which was not included in either (b4) or (b5).
This happens if and only~if $\tau = \ell+1$, i.e.~$u_\tau = r_i$,
in which case all the edges having $v_\tau$ as a descendant
are irregular by case~(a).
\myendremark

\bigskip

Now we consider the general case $k\ge1$.
The following construction is inspired by \cite[proof of Lemma~2]{Chauve_00}.

\proofof{Proposition~\ref{prop.bijection.model21}}
Given a tree $T$ rooted at $r$ in Model~2,
let $v_1<v_2<\cdots<v_k$ be the $k$ children of the vertex~$1$.
For any vertex $i$ other than the root,
we define its \textbfit{top~ancestor}
to be the ancestor of $i$ (possibly $i$ itself) that is a child of the root.
We construct the bijection $\sigma$ in the following three cases:
\begin{quote}
   Case~I:  $v_1 < r$.

   Case~II: $v_1 > r$ and the top ancestor of vertex~1 is
       lower-numbered than the root.

   Case~III: $v_1 > r$ and the top ancestor of vertex~1 is
       higher-numbered than the root.
\end{quote}

\begin{figure}[p]
  \begin{center}
  \begin{tikzpicture}[baseline=(current bounding box.base)]
    \node at (0,1) {$T$};
    \node (root) at (0,0) {};
    \node (Tr) at (0,-1) {};
    \fill (root) circle(.1) node[above] {$r$};
    \draw (Tr) ellipse(.5 and 1) node {$T_r$};
    \node (1) at (0,-2) {};
    \node (v1) at (-1,-3) {};
    \node (v2) at (1,-3) {};
    \node (vk) at (2,-3) {};
    \fill (1) circle(.1) node[above] {$1$};
    \fill (v1) circle(.1) node[left] {$v_1$};
    \draw (v2) circle(.1) node[left] {$v_2$};
    \draw (vk) circle(.1) node[right] {$v_k$};
    \draw (1)--(v1) (1)--(v2) (1)--(vk);
    \node at (1.5,-3) {$\cdots$};
    \node (T2) at (1,-4) {};
    \node (Tk) at (2,-4) {};
    \draw (T2) ellipse (.3 and 1) node {$T_2$};
    \draw (Tk) ellipse (.3 and 1) node {$T_k$};
    \node (L) at (-2,-4) {};
    \node (H) at (0,-4) {};
    \node (DL) at (-2,-6) {};
    \node (DH) at (0,-6) {};
    \node (DLtop) at (-2,-5.5) {};
    \node (DHtop) at (0,-5.5) {};
    \draw (L) ellipse (.5 and .3) node {$L$};
    \draw (H) ellipse (.5 and .3) node {$H$};
    \draw (DL) ellipse (.5 and 1) node {$D_L$};
    \draw (DH) ellipse (.5 and 1) node {$D_H$};
    \draw (v1)--(L)--(DLtop) (v1)--(H)--(DHtop);
    \draw[ultra thick, blue] (1)--(v1) (1)--(v2) (1)--(vk);
    \draw[ultra thick, red] (v1)--(H);
  \end{tikzpicture}
  \qquad \raisebox{-70pt}{$\longrightarrow$} \qquad
  \begin{tikzpicture}[baseline=(current bounding box.base)]
    \node at (0,1) {$\sigma(T)$};
    \node (v1) at (0,0) {};
    \node (r) at (0,-1) {};
    \node (v2) at (1,-1) {};
    \node (vk) at (2,-1) {};
    \fill (v1) circle(.1) node[above] {$v_1$};
    \fill (r) circle(.1) node[left] {$r$};
    \draw (v2) circle(.1) node[left] {$v_2$};
    \draw (vk) circle(.1) node[right] {$v_k$};
    \draw (v1)--(r) (v1)--(v2) (v1)--(vk);
    \node at (1.5,-1) {$\cdots$};
    \node (L) at (-1.5,-1) {};
    \node (DL) at (-1.5,-3) {};
    \node (DLtop) at (-1.5,-2.5) {};
    \draw (L) ellipse (.5 and .3) node {$L$};
    \draw (DL) ellipse (.5 and 1) node {$D_L$};
    \draw (v1)--(L)--(DLtop);
    \node (Tr) at (0,-2) {};
    \node (1) at (0,-3) {};
    \draw (Tr) ellipse(.5 and 1) node {$T_r$};
    \fill (1) circle(.1) node[above] {$1$};
    \node (T2) at (1,-2) {};
    \node (Tk) at (2,-2) {};
    \draw (T2) ellipse (.3 and 1) node {$T_2$};
    \draw (Tk) ellipse (.3 and 1) node {$T_k$};
    \node (H) at (0,-4) {};
    \node (DH) at (0,-6) {};
    \node (DHtop) at (0,-5.5) {};
    \draw (H) ellipse (.5 and .3) node {$H$};
    \draw (DH) ellipse (.5 and 1) node {$D_H$};
    \draw (1)--(H)--(DHtop);
    \draw[ultra thick, blue] (v1)--(r) (v1)--(v2) (v1)--(vk);
    \draw[ultra thick, red] (1)--(H);
  \end{tikzpicture}
  \end{center}
  \caption{
     The trees $T$ and $\sigma(T)$ in Case~I.
     The blue edges in $T$ [resp.~$\sigma(T)$] are proper (resp.~regular).
     The red edges in $T$ [resp.~$\sigma(T)$]
     could be either proper or improper (resp.~regular or irregular),
     depending on the behavior of their descendants.
  }
  \label{fig-bijec1&2-case1}
\end{figure}
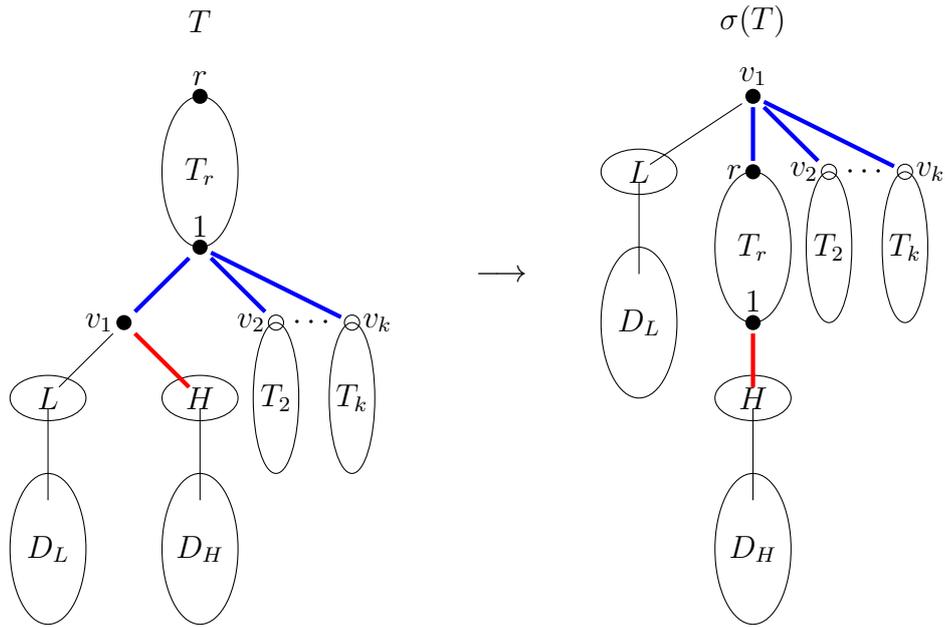

\medskip

\noindent
{\bf Case I: $\bm{v_1<r}$.}

Let $L$ (resp.~$H$) denote the lower (resp.~higher)-numbered children
of vertex $v_1$,
and let $D_L$ (resp.~$D_H$) denote their descendants
excluding those in $L$ (resp.~$H$) itself.
We construct the tree $\sigma(T)$ as follows:

\begin{enumerate}[Step 1.]
  \item Delete the $k$ edges emanating from vertex~$1$ to its children,
        and the edges from $v_1$ to its higher-numbered children $H$.
  \item Attach the trees rooted at the vertices of $H$ onto vertex~$1$
        via new edges.
  \item Attach the trees rooted at $r,v_2,\ldots,v_k$ onto vertex $v_1$
        via new edges.
\end{enumerate}
See Figure~\ref{fig-bijec1&2-case1}.
We obtain thereby a tree $\sigma(T)$ rooted at $v_1$
with $k$ higher-numbered children $r,v_2,\ldots,v_k$,
which also has the following properties:
\begin{itemize}
  \item[A1)] all the children of vertex~$1$ are higher-numbered than the root;
  \item[B1)] the top ancestor of vertex~$1$ is higher-numbered than the root.
\end{itemize}

The $k$ edges in $T$ that emanate from vertex~$1$ to its children
(shown in blue in Figure~\ref{fig-bijec1&2-case1}) are clearly proper.
These edges are mapped to the $k$ edges in $\sigma(T)$ that emanate from
the root $v_1$ to its higher-numbered children,
which are regular by rule~(I3) of Definition~\ref{def-irregular}.

An edge $e \in T$ that emanates from vertex $v_1$
to a higher-numbered child $h \in H$
(shown in red in Figure~\ref{fig-bijec1&2-case1})
is improper if (and only~if) there is a descendant of $h$ that is $< v_1$.
Such an edge is mapped to the edge $\psi_T(e) \in \sigma(T)$
that emanates from vertex~$1$ to its child $h \in H$;
and since all the children of vertex~$1$ are higher-numbered
than the root $v_1$,
rule~(I4) of Definition~\ref{def-irregular} applies
and says that the edge $\psi_T(e)$ is irregular if and only~if
$e$ is improper.

All other edges $e \in T$
(shown in black in Figure~\ref{fig-bijec1&2-case1})
have the property that $\psi_T(e) = e$.
Moreover, if $e \notin T_r$,
then the descendants of $e$ in $T$
are the same as its descendants in $\sigma(T)$,
so by rules~(I1), (I2) and (I6) of Definition~\ref{def-irregular},
the edge $e$ is proper/improper in $T$
exactly when it is regular/irregular in $\sigma(T)$.
Finally, all the edges $e \in T_r$ are improper in $T$
(because vertex~$1$ is a descendant),
and they are irregular in $\sigma(T)$ by rules~(I1) and (I2).

We have therefore shown that the bijection $\psi_T$
maps proper/improper edges in $T$ 
onto regular/irregular edges in $\sigma(T)$.

\bigskip

\begin{figure}[t]
  \begin{center}
  \begin{tikzpicture}[baseline=(current bounding box.base)]
    \node at (0,1) {$T$};
    \node (root) at (0,0) {};
    \node (L) at (-1,-1) {};
    \node (H) at (1,-1) {};
    \node (DL) at (-1,-3) {};
    \node (DH) at (1,-3) {};
    \node (DLtop) at (-1,-2.5) {};
    \node (DHtop) at (1,-2.5) {};
    \fill (root) circle(.1) node[above] {$r$};
    \draw (L) ellipse (.5 and .3) node {$L$};
    \draw (H) ellipse (.5 and .3) node {$H$};
    \draw (DL) ellipse (.5 and 1) node {$D_L$};
    \draw (DH) ellipse (.5 and 1) node {$D_H$};
    \draw (root)--(L)--(DLtop) (root)--(H)--(DHtop);
    \node (1) at (-1,-4) {};
    \node (v1) at (-2,-5) {};
    \node (vk) at (0,-5) {};
    \fill (1) circle(.1) node[above] {$1$};
    \draw (v1) circle(.1) node[left] {$v_1$};
    \draw (vk) circle(.1) node[right] {$v_k$};
    \draw (1)--(v1) (1)--(vk);
    \node at (-1,-5) {$\cdots$};
    \node (T1) at (-2,-6) {};
    \node (Tk) at (0,-6) {};
    \draw (T1) ellipse (.3 and 1) node {$T_1$};
    \draw (Tk) ellipse (.3 and 1) node {$T_k$};
    \draw[ultra thick, blue] (1)--(v1) (1)--(vk);
    \draw[ultra thick, red] (root)--(H);
  \end{tikzpicture}
  \qquad \raisebox{-70pt}{$\longrightarrow$} \qquad
  \begin{tikzpicture}[baseline=(current bounding box.base)]
    \node at (0,1) {$\sigma(T)$};
    \node (root) at (0,0) {};
    \node (v1) at (.5,-1) {};
    \node (vk) at (2,-1) {};
    \fill (root) circle(.1) node[above] {$r$};
    \draw (v1) circle(.1) node[left] {$v_1$};
    \draw (vk) circle(.1) node[right] {$v_k$};
    \draw (root)--(v1) (root)--(vk);
    \node at (1.25,-1) {$\cdots$};
    \node (L) at (-1,-1) {};
    \node (DL) at (-1,-3) {};
    \node (DLtop) at (-1,-2.5) {};
    \draw (L) ellipse (.5 and .3) node {$L$};
    \draw (DL) ellipse (.5 and 1) node {$D_L$};
    \draw (root)--(L)--(DLtop);
    \node (1) at (-1,-4) {};
    \node (H) at (-1,-5) {};
    \node (DH) at (-1,-7) {};
    \node (DHtop) at (-1,-6.5) {};
    \fill (1) circle(.1) node[above] {$1$};
    \draw (H) ellipse (.5 and .3) node {$H$};
    \draw (DH) ellipse (.5 and 1) node {$D_H$};
    \draw (1)--(H)--(DHtop);
    \node (T1) at (.5,-2) {};
    \node (Tk) at (2,-2) {};
    \draw (T1) ellipse (.3 and 1) node {$T_1$};
    \draw (Tk) ellipse (.3 and 1) node {$T_k$};
    \draw[ultra thick, blue] (root)--(v1) (root)--(vk);
    \draw[ultra thick, red] (1)--(H);
  \end{tikzpicture}
  \end{center}
  \caption{
     The trees $T$ and $\sigma(T)$ in Case~II.
     The blue edges in $T$ [resp.~$\sigma(T)$] are proper (resp.~regular).
     The red edges in $T$ [resp.~$\sigma(T)$]
     could be either proper or improper (resp.~regular or irregular),
     depending on the behavior of their descendants.
  }
  \label{fig-bijec1&2-case2}
\end{figure}
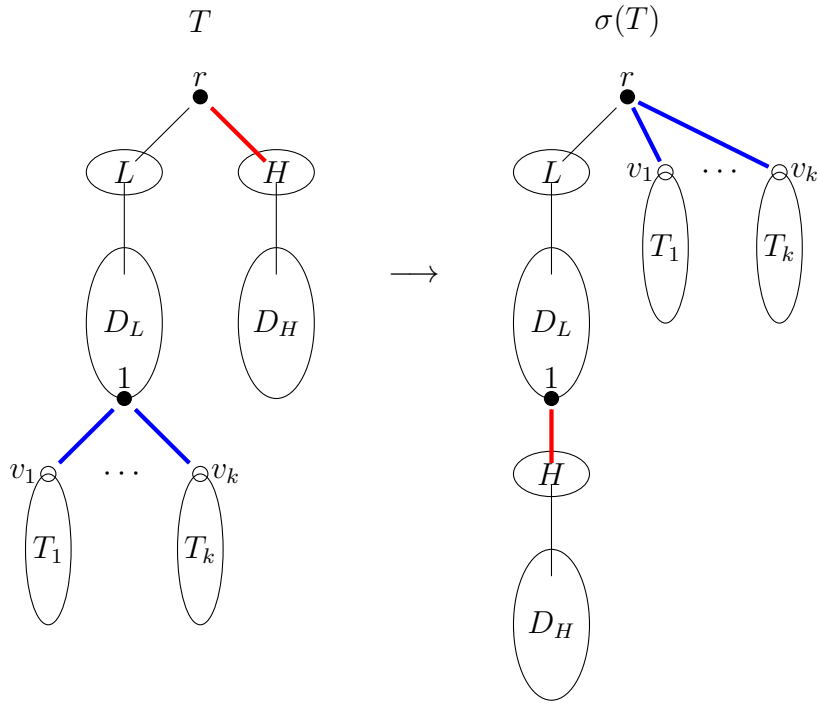

\noindent
{\bf Case II: $\bm{v_1 > r}$ and the top ancestor of vertex~1 is
   lower-numbered than the root.}

Let $L$ (resp.~$H$) denote the lower (resp.~higher)-numbered children
of the root~$r$,
and let $D_L$ (resp.~$D_H$) denote their descendants
excluding those in $L$ (resp.~$H$) itself.
We construct the tree $\sigma(T)$ as follows:

\begin{enumerate}[Step 1.]
  \item Delete the $k$ edges emanating from vertex~$1$ to its children,
        and the edges from $r$ to its higher-numbered children $H$.
  \item Attach the trees rooted at the vertices of $H$ onto vertex~$1$
        via new edges.
  \item Attach the trees rooted at $v_1,v_2,\ldots,v_k$ onto the root~$r$
        via new edges.
\end{enumerate}
(Note that this is identical to Case~I but with the roles of $v_1$ and $r$
 interchanged.)
See Figure~\ref{fig-bijec1&2-case2}.
We obtain thereby a tree $\sigma(T)$ rooted at $r$
with $k$ higher-numbered children $v_1,v_2,\ldots,v_k$,
which also has the following properties:
\begin{itemize}
  \item[A1)] all the children of vertex~$1$ are higher-numbered than the root;
  \item[B2)] the top ancestor of vertex~$1$ is lower-numbered than the root.
\end{itemize}

The $k$ edges in $T$ that emanate from vertex~$1$ to its children
(shown in blue in Figure~\ref{fig-bijec1&2-case2}) are clearly proper.
These edges are mapped to the $k$ edges in $\sigma(T)$ that emanate from
the root $r$ to its higher-numbered children,
which are regular by rule~(I3) of Definition~\ref{def-irregular}.

An edge $e \in T$ that emanates from the root~$r$
to a higher-numbered child $h \in H$
(shown in red in Figure~\ref{fig-bijec1&2-case2})
is improper if (and only~if) there is a descendant of $h$ that is $< r$.
Such an edge is mapped to the edge $\psi_T(e) \in \sigma(T)$
that emanates from vertex~$1$ to its child $h \in H$;
and since all the children of vertex~$1$ are higher-numbered
than the root $r$,
rule~(I4) of Definition~\ref{def-irregular} applies
and says that the edge $\psi_T(e)$ is irregular if and only~if
$e$ is improper.

All other edges $e \in T$
(shown in black in Figure~\ref{fig-bijec1&2-case1})
have the property that $\psi_T(e) = e$.
Moreover, if $e \notin T \restrict (\{r\} \cup L \cup D_L)$,
then the descendants of $e$ in $T$
are the same as its descendants in $\sigma(T)$,
so by rules~(I1), (I2) and (I6) of Definition~\ref{def-irregular},
the edge $e$ is proper/improper in $T$
exactly when it is regular/irregular in $\sigma(T)$.
Finally, all the edges $e \in T \restrict (\{r\} \cup L \cup D_L)$
are improper in $T$ (because vertex~$1$ is a descendant),
and they are irregular in $\sigma(T)$ by rules~(I1) and (I2).

We have therefore shown that the bijection $\psi_T$
maps proper/improper edges in $T$ 
onto regular/irregular edges in $\sigma(T)$.

\bigskip

\noindent
{\bf Case III: $\bm{v_1 > r}$ and the top ancestor of vertex~1 is
   higher-numbered than the root.}

Let $L$ (resp.~$H$) denote the lower- (resp.~higher-) numbered children
of the root $r$, and let $D_L$ (resp.~$D_H$) denote their descendants
excluding those in $L$ (resp.~$H$) itself.
Let $u$ be the vertex on the path from the root $r$ to vertex~$1$
such that the path from $r$ to $u$ is maximal increasing.
Clearly, $r<u$.
The first two steps in constructing the tree $\sigma(T)$ are as follows:

\begin{enumerate}[Step 1.]
  \item Delete the $k$ edges from vertex~$1$ to its children,
        and denote by $T_0$ the subtree rooted at $r$ in which $1$ is a leaf.
  \item Use the bijection $\phi$ constructed in Lemma~\ref{lem-bijec1&2-k=0}
        to yield a tree $\phi(T_0)$.
        Note that the vertex~$u$ plays the role of $r_i$
        in Lemma~\ref{lem-bijec1&2-k=0},
        so that the operator $\relab{1}{u}$
        acts on the maximal increasing subtree $(T_0)_{\rm max}$ of $T_0$
        rooted at $r$.
        Note also that $\phi(T_0)$ is rooted at $u$,
        which has no higher-numbered children.
        See Figure~\ref{fig-bijec1&2-case3}.
\end{enumerate}

\begin{figure}[p]
  \begin{center}
  \begin{tikzpicture}[baseline=(current bounding box.base)]
    \node at (0,1) {$T$};
    \node (root) at (0,0) {};
    \node (L) at (-1,-1) {};
    \node (H) at (1,-1) {};
    \node (DL) at (-1,-3) {};
    \node (DH) at (1,-3) {};
    \node (DLtop) at (-1,-2.5) {};
    \node (DHtop) at (1,-2.5) {};
    \fill (root) circle(.1) node[above] {$r$};
    \draw (L) ellipse (.5 and .3) node {$L$};
    \draw (H) ellipse (.5 and .3) node {$H$};
    \draw (DL) ellipse (.5 and 1) node {$D_L$};
    \draw (DH) ellipse (.5 and 1) node {$D_H$};
    \draw (root)--(L)--(DLtop) (root)--(H)--(DHtop);
    \node (1) at (1,-4) {};
    \node (v1) at (0,-5) {};
    \node (v2) at (1,-5) {};
    \node (vk) at (3,-5) {};
    \fill (1) circle(.1) node[above] {$1$};
    \draw (v1) circle(.1) node[left] {$v_1$};
    \draw (v2) circle(.1) node[left] {$v_2$};
    \draw (vk) circle(.1) node[right] {$v_k$};
    \draw (1)--(v1) (1)--(v2) (1)--(vk);
    \node at (2,-5) {$\cdots$};
    \node (T1) at (0,-6) {};
    \node (T2) at (1,-6) {};
    \node (Tk) at (3,-6) {};
    \draw (T1) ellipse (.3 and 1) node {$T_1$};
    \draw (T2) ellipse (.3 and 1) node {$T_2$};
    \draw (Tk) ellipse (.3 and 1) node {$T_k$};
    \node (u) at (1.25,-2.5) {};
    \node (un) at (2,-2) {};
    \fill  (u) circle(.1);
    \draw[->] (u)--(un) node {$u$};
  \end{tikzpicture}
  \hspace*{3cm}
  \begin{tikzpicture}[baseline=(current bounding box.base)]
    \node at (0,1) {$T_0$};
    \node (root) at (0,0) {};
    \node (L) at (-1,-1) {};
    \node (H) at (1,-1) {};
    \node (DL) at (-1,-3) {};
    \node (DH) at (1,-3) {};
    \node (DLtop) at (-1,-2.5) {};
    \node (DHtop) at (1,-2.5) {};
    \fill (root) circle(.1) node[above] {$r$};
    \draw (L) ellipse (.5 and .3) node {$L$};
    \draw (H) ellipse (.5 and .3) node {$H$};
    \draw (DL) ellipse (.5 and 1) node {$D_L$};
    \draw (DH) ellipse (.5 and 1) node {$D_H$};
    \draw (root)--(L)--(DLtop) (root)--(H)--(DHtop);
    \node (u) at (1.25,-2.5) {};
    \node (un) at (2,-2) {};
    \fill  (u) circle(.1);
    \draw[->] (u)--(un) node {$u$};
    \node (1) at (1,-4) {};
    \fill (1) circle(.1) node[below] {$1$};
  \end{tikzpicture}
  \\[2cm]
  \begin{tikzpicture}[baseline=(current bounding box.base)]
    \node at (0,1) {$\phi(T_0)$};
    \node (root) at (0,0) {};
    \node (Tu) at (0,-1) {};
    \fill (root) circle(.1) node[above] {$u$};
    \draw (Tu) ellipse (.5 and 1) node {$T_u$};
    \node (1) at (0,-2) {};
    \fill (1) circle(.1) node[above] {$1$};
    \node (rhoT0max) at (0,-3.5) {};
    \draw (rhoT0max) ellipse (1.5 and 1.5) node {$\relab{1}{u}((T_0)_{\rm max})$};
    \node (r) at (-1,-3) {};
    \node (L) at (-2,-4) {};
    \node (DL) at (-2,-6) {};
    \node (DLtop) at (-2,-5.5) {};
    \fill (r) circle(.1) node[above] {$r$};
    \draw (L) ellipse (.5 and .3) node {$L$};
    \draw (DL) ellipse (.5 and 1) node {$D_L$};
    \draw (1)--(r)--(L)--(DLtop);
    \node (r1) at (1.3,-3.6) {};
    \node (T1) at (2.2,-4.4) {};
    \fill (r1) circle(.1);
    \draw[rotate=45] (T1) ellipse (.5 and 1.2);
    \node (rp) at (.5,-4.5) {};
    \node (Tp) at (.5,-5.5) {};
    \fill (rp) circle(.1);
    \draw (Tp) ellipse (.5 and 1);
    \node[rotate=40] at (1.8,-5.5) {$\cdots$};
  \end{tikzpicture}
  \end{center}
  \caption{
     The trees $T$, $T_0$ and $\phi(T_0)$ in Case~III.
     Here the vertex $u$ can be either in $H$ or in $D_H$.
  }
  \label{fig-bijec1&2-case3}
\end{figure}

Then we distinguish two subcases, according as $u < v_1$ or $u > v_1$:

\bigskip

\noindent
{\bf Case III(a): $\bm{u<v_1}$.}

\begin{figure}[tp]
  \begin{center}
  \begin{tikzpicture}[baseline=(current bounding box.base)]
    \node at (0,1) {$T$};
    \node (root) at (0,0) {};
    \node (T0) at (0,-1.5) {};
    \fill (root) circle(.1) node[above] {$r$};
    \draw (T0) ellipse (1 and 1.5) node {$T_0$};
    \node (1) at (0,-3) {};
    \node (v1) at (-1,-4) {};
    \node (v2) at (0,-4) {};
    \node (vk) at (2,-4) {};
    \fill (1) circle(.1) node[above] {$1$};
    \draw (v1) circle(.1) node[left] {$v_1$};
    \draw (v2) circle(.1) node[left] {$v_2$};
    \draw (vk) circle(.1) node[right] {$v_k$};
    \draw[ultra thick, blue] (1)--(v1) (1)--(v2) (1)--(vk);
    \node at (1,-4) {$\cdots$};
    \node (T1) at (-1,-5) {};
    \node (T2) at (0,-5) {};
    \node (Tk) at (2,-5) {};
    \draw (T1) ellipse (.3 and 1) node {$T_1$};
    \draw (T2) ellipse (.3 and 1) node {$T_2$};
    \draw (Tk) ellipse (.3 and 1) node {$T_k$};
    \node (u) at (.75,-2) {};
    \node (un) at (1.5,-2) {};
    \fill  (u) circle(.1);
    \draw[->] (u)--(un) node {$u$};
  \end{tikzpicture}
  \qquad \raisebox{-70pt}{$\longrightarrow$} \qquad
  \begin{tikzpicture}[baseline=(current bounding box.base)]
    \node at (0,1) {$\sigma(T)$};
    \node (root) at (0,0) {};
    \node (phiT0) at (0,-1.5) {};
    \fill (root) circle(.1) node[above] {$u$};
    \draw (phiT0) ellipse (1 and 1.5);
    \node at (0,-2) {$\phi(T_0)$};
    \node (new1) at (0,-0.5) {};
    \fill (new1) circle(.1) node[right] {$1$};
    \node (newr) at (-0.5,-1.2) {};
    \fill (newr) circle(.1) node[below] {$r$};
    \draw (new1)--(newr);
    \node (v1) at (2,-1) {};
    \node (v2) at (3,-1) {};
    \node (vk) at (5,-1) {};
    \node (T1) at (2,-2) {};
    \node (T2) at (3,-2) {};
    \node (Tk) at (5,-2) {};
    \draw (v1) circle(.1) node[left] {$v_1$};
    \draw (v2) circle(.1) node[left] {$v_2$};
    \draw (vk) circle(.1) node[right] {$v_k$};
    \draw (T1) ellipse (.3 and 1) node {$T_1$};
    \draw (T2) ellipse (.3 and 1) node {$T_2$};
    \draw (Tk) ellipse (.3 and 1) node {$T_k$};
    \draw[ultra thick, blue] (root)--(v1) (root)--(v2) (root)--(vk);
    \node at (4,-1) {$\cdots$};
  \end{tikzpicture}
  \end{center}
  \caption{The trees $T$ and $\sigma(T)$ in Case III(a).}
  \label{fig-bijec1&2-case3a}
\end{figure}
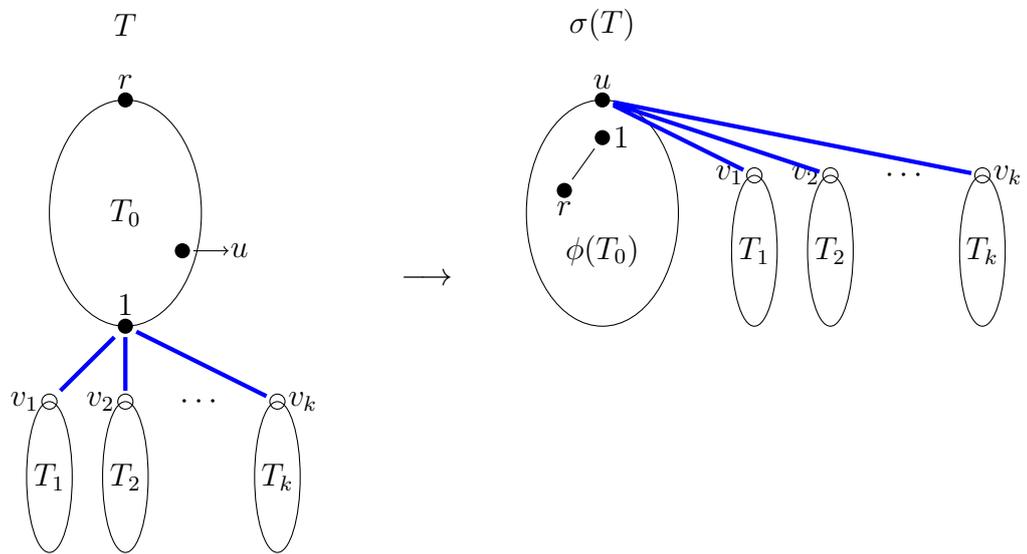

\begin{enumerate}[Step 3.]
  \item Attach the trees rooted at $v_1,v_2,\ldots,v_k$ onto the root $u$
     of $\phi(T_0)$ via new edges.
\end{enumerate}
See Figure~\ref{fig-bijec1&2-case3a}.
We obtain thereby a tree $\sigma(T)$ rooted at $u$
with $k$ higher-numbered children $v_1,v_2,\ldots,v_k$,
which also has the following properties:
\begin{itemize}
  \item[A2)] vertex~$1$ has at least one child that is lower-numbered
     than the root;
  \item[B2)] the top ancestor of vertex~$1$ is lower-numbered than the root.
\end{itemize}

The $k$ edges in $T$ that emanate from vertex~$1$ to its children
(shown in blue in Figure~\ref{fig-bijec1&2-case3a}) are clearly proper.
These edges are mapped to the $k$ edges in $\sigma(T)$ that emanate from
the root $u$ to its higher-numbered children,
which are regular by rule~(I3) of Definition~\ref{def-irregular}.
By the discussion in Lemma~\ref{lem-bijec1&2-k=0},
the proper/improper edges in $T_0$
are all mapped to regular/irregular edges in $\phi(T_0)$.
Finally, the proper/improper edges in $T_1,\ldots,T_k \subseteq T$
are mapped to regular/irregular edges in $T_1,\ldots,T_k \subseteq \sigma(T)$
by rules~(I1), (I2) and (I6).
We have therefore shown that the bijection $\psi_T$
maps proper/improper edges in $T$ 
onto regular/irregular edges in $\sigma(T)$.

\bigskip

\noindent
{\bf Case III(b): $\bm{u>v_1}$.}

\begin{figure}[p]
  \begin{center}
  \begin{tikzpicture}[baseline=(current bounding box.base)]
    \node at (0,1) {$T$};
    \node (root) at (0,0) {};
    \node (L) at (-1,-1) {};
    \node (H) at (1,-1) {};
    \node (DL) at (-1,-3) {};
    \node (DH) at (1,-3) {};
    \node (DLtop) at (-1,-2.5) {};
    \node (DHtop) at (1,-2.5) {};
    \fill (root) circle(.1) node[above] {$r$};
    \draw (L) ellipse (.5 and .3) node {$L$};
    \draw (H) ellipse (.5 and .3) node {$H$};
    \draw (DL) ellipse (.5 and 1) node {$D_L$};
    \draw (DH) ellipse (.5 and 1) node {$D_H$};
    \draw (root)--(L)--(DLtop) (root)--(H)--(DHtop);
    \node (1) at (1,-4) {};
    \node (v1) at (0,-5) {};
    \node (v2) at (2.5,-5) {};
    \node (vk) at (4.5,-5) {};
    \fill (1) circle(.1) node[above] {$1$};
    \draw (v1) circle(.1) node[left] {$v_1$};
    \draw (v2) circle(.1) node[left] {$v_2$};
    \draw (vk) circle(.1) node[right] {$v_k$};
    \draw (1)--(v1) (1)--(v2) (1)--(vk);
    \node at (3.5,-5) {$\cdots$};
    \node (L1) at (-1,-6) {};
    \node (H1) at (1,-6) {};
    \node (DL1) at (-1,-8) {};
    \node (DH1) at (1,-8) {};
    \node (DL1top) at (-1,-7.5) {};
    \node (DH1top) at (1,-7.5) {};
    \draw (L1) ellipse (.8 and .3) node {$L_1$};
    \draw (H1) ellipse (.8 and .3) node {$H_1$};
    \draw (DL1) ellipse (.8 and 1) node {$D_{L_1}$};
    \draw (DH1) ellipse (.8 and 1) node {$D_{H_1}$};
    \draw (v1)--(L1)--(DL1top) (v1)--(H1)--(DH1top);
    \node (T2) at (2.5,-6) {};
    \node (Tk) at (4.5,-6) {};
    \draw (T2) ellipse (.3 and 1) node {$T_2$};
    \draw (Tk) ellipse (.3 and 1) node {$T_k$};
    \node (u) at (1.25,-2.5) {};
    \node (un) at (2,-2) {};
    \fill  (u) circle(.1);
    \draw[->] (u)--(un) node {$u$};
    \node (w) at (.5,-6) {};
    \node (wn) at (0,-6.7) {};
    \fill  (w) circle(.1);
    \draw[->] (w)--(wn) node {$w$};
  \end{tikzpicture}
  \hspace*{2cm}
  \begin{tikzpicture}[baseline=(current bounding box.base)]
    \node at (0,1) {$T_1$};
    \node (root) at (0,0) {};
    \node (L1) at (-1,-1) {};
    \node (H1) at (1,-1) {};
    \node (DL1) at (-1,-3) {};
    \node (DH1) at (1,-3) {};
    \node (DL1top) at (-1,-2.5) {};
    \node (DH1top) at (1,-2.5) {};
    \fill (root) circle(.1) node[above] {$v_1$};
    \draw (L1) ellipse (.8 and .3) node {$L_1$};
    \draw (H1) ellipse (.8 and .3) node {$H_1$};
    \draw (DL1) ellipse (.8 and 1) node {$D_{L_1}$};
    \draw (DH1) ellipse (.8 and 1) node {$D_{H_1}$};
    \draw (root)--(L1)--(DL1top) (root)--(H1)--(DH1top);
    \node (w) at (.5,-1) {};
    \node (wn) at (0,-1.7) {};
    \fill  (w) circle(.1);
    \draw[->] (w)--(wn) node {$w$};
  \end{tikzpicture}
  \\[2cm]
  \begin{tikzpicture}[baseline=(current bounding box.base)]
    \node at (0,1) {$\phi'(T_1)$};
    \node (root) at (0,0) {};
    \node (L1) at (-1.5,-1) {};
    \node (DL1) at (-1.5,-3) {};
    \node (DL1top) at (-1.5,-2.5) {};
    \fill (root) circle(.1) node[above] {$v_1$};
    \draw (L1) ellipse (.8 and .3) node {$L_1$};
    \draw (DL1) ellipse (.8 and 1) node {$D_{L_1}$};
    \draw (root)--(L1)--(DL1top);
    \node (w) at (1,-1) {};
    \node (rhoT1max) at (1,-2.5) {};
    \fill (w) circle(.1) node[above] {$w$};
    \draw[blue, ultra thick] (root)--(w);
    \draw (rhoT1max) ellipse (1.5 and 1.5) node {$\relab{u}{v_1}((T_1)_{\rm max})$};
    \node (u) at (1.7,-3.5) {};
    \node (un) at (2.5,-4.2) {};
    \fill  (u) circle(.1);
    \draw[->] (u)--(un) node {$u$};
    \node (r1) at (2,-1.8) {};
    \node (T1) at (3,-2.5) {};
    \fill (r1) circle(.1);
    \draw[rotate=50] (T1) ellipse (.5 and 1.2);
    \node (rp) at (.5,-3.5) {};
    \node (Tp) at (.5,-4.5) {};
    \fill (rp) circle(.1);
    \draw (Tp) ellipse (.5 and 1);
  \end{tikzpicture}
  \end{center}
  \caption{
     The trees $T$, $T_1$ and $\phi'(T_1)$ in Case III(b).
     Note that here $T$ is the same as in Figure~\ref{fig-bijec1&2-case3},
     but the subtree $T_1$ is shown in more detail.
  }
  \label{fig-bijec1&2-case3b-phi'}
\end{figure}
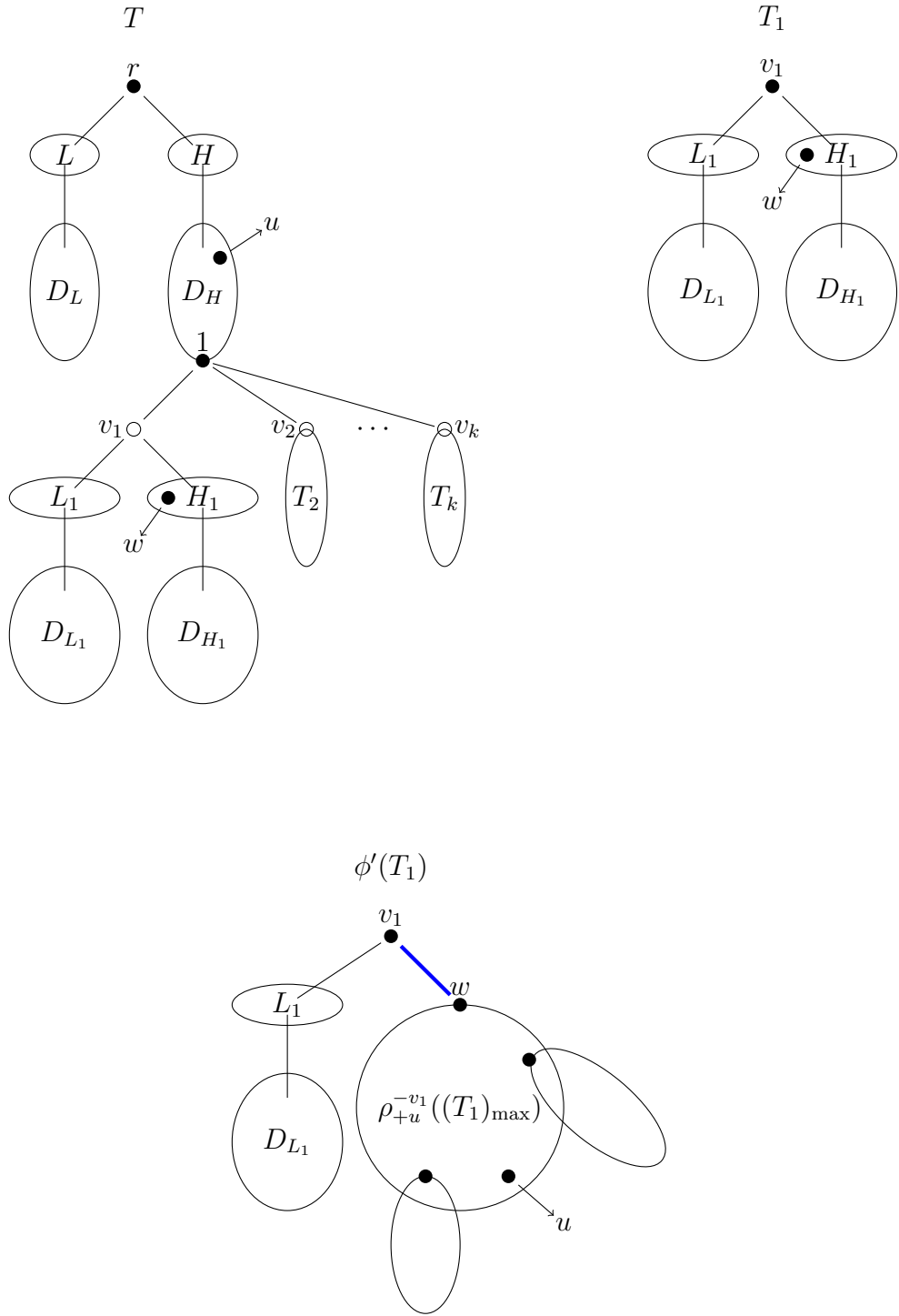

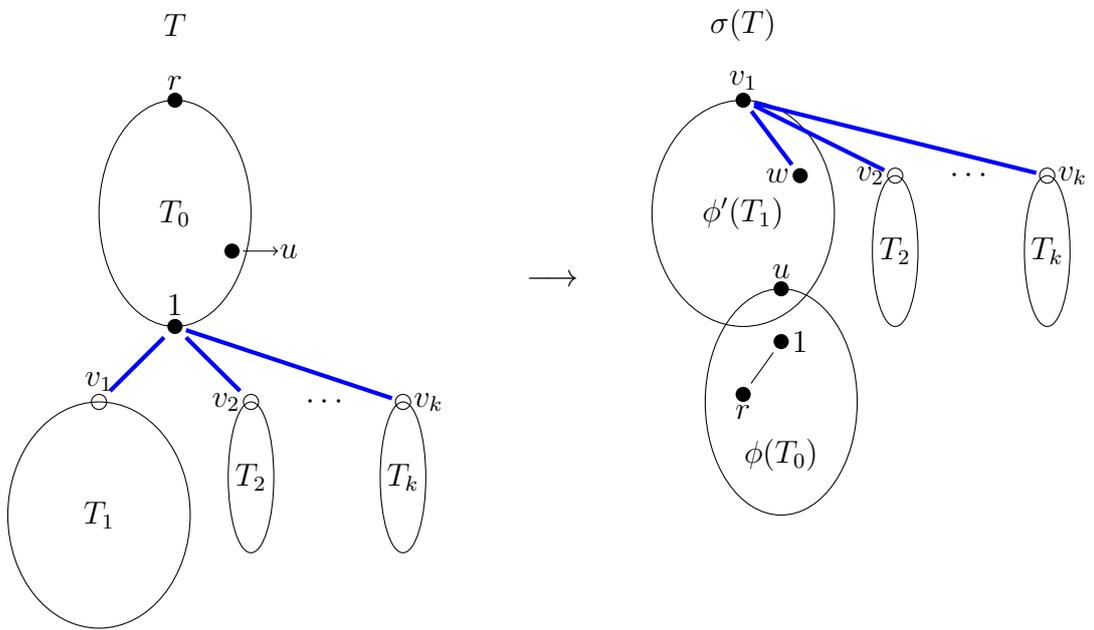
\begin{figure}[tp]
  \begin{center}
  \begin{tikzpicture}[baseline=(current bounding box.base)]
    \node at (0,1) {$T$};
    \node (root) at (0,0) {};
    \node (T0) at (0,-1.5) {};
    \fill (root) circle(.1) node[above] {$r$};
    \draw (T0) ellipse (1 and 1.5) node {$T_0$};
    \node (1) at (0,-3) {};
    \node (v1) at (-1,-4) {};
    \node (v2) at (1,-4) {};
    \node (vk) at (3,-4) {};
    \fill (1) circle(.1) node[above] {$1$};
    \draw (v1) circle(.1) node[above] {$v_1$};
    \draw (v2) circle(.1) node[left] {$v_2$};
    \draw (vk) circle(.1) node[right] {$v_k$};
    \draw[ultra thick, blue] (1)--(v1) (1)--(v2) (1)--(vk);
    \node at (2,-4) {$\cdots$};
    \node (T1) at (-1,-5.5) {};
    \node (T2) at (1,-5) {};
    \node (Tk) at (3,-5) {};
    \draw (T1) ellipse (1.2 and 1.5) node {$T_1$};
    \draw (T2) ellipse (.3 and 1) node {$T_2$};
    \draw (Tk) ellipse (.3 and 1) node {$T_k$};
    \node (u) at (.75,-2) {};
    \node (un) at (1.5,-2) {};
    \fill  (u) circle(.1);
    \draw[->] (u)--(un) node {$u$};
  \end{tikzpicture}
  \qquad \raisebox{-70pt}{$\longrightarrow$} \qquad
  \begin{tikzpicture}[baseline=(current bounding box.base)]
    \node at (0,1) {$\sigma(T)$};
    \node (root) at (0,0) {};
    \node (phi'T1) at (0,-1.5) {};
    \fill (root) circle(.1) node[above] {$v_1$};
    \draw (phi'T1) ellipse (1.2 and 1.5) node {$\phi'(T_1)$};
    \node (v2) at (2,-1) {};
    \node (vk) at (4,-1) {};
    \node (T2) at (2,-2) {};
    \node (Tk) at (4,-2) {};
    \draw (v2) circle(.1) node[left] {$v_2$};
    \draw (vk) circle(.1) node[right] {$v_k$};
    \draw (T2) ellipse (.3 and 1) node {$T_2$};
    \draw (Tk) ellipse (.3 and 1) node {$T_k$};
    \draw[ultra thick, blue] (root)--(v2) (root)--(vk);
    \node at (3,-1) {$\cdots$};
    \node (w) at (.75,-1) {};
    \fill (w) circle(.1) node[left] {$w$};
    \draw[ultra thick, blue] (root)--(w);
    \node (u) at (.5,-2.5) {};
    \fill  (u) circle(.1) node[above] {$u$};
    \node (phiT0) at (.5,-4) {};
    \draw (phiT0) ellipse (1 and 1.5);
    \node at (.5,-4.7) {$\phi(T_0)$};
    \node (new1) at (0.5,-3.2) {};
    \fill (new1) circle(.1) node[right] {$1$};
    \node (newr) at (0,-3.9) {};
    \fill (newr) circle(.1) node[below] {$r$};
    \draw (new1)--(newr);
  \end{tikzpicture}
  \end{center}
  \caption{The trees $T$ and $\sigma(T)$ in Case III(b).}
  \label{fig-bijec1&2-case3b}
\end{figure}

Consider the subtree $T_1 \subseteq T$ rooted at $v_1$.
Let $L_1$ (resp.~$H_1$) denote the lower- (resp.~higher-) numbered children
of the vertex~$v_1$,
and let $D_{L_1}$ (resp.~$D_{H_1}$) denote their descendants
excluding those in $L_1$ (resp.~$H_1$) itself.
Let $w$ be the smallest vertex in $H_1$.
Clearly, $w>v_1$.
We then proceed as follows:

\begin{itemize}
  \item[Step 3.] Let $(T_1)_{\rm max}$ be the maximal increasing
        subtree of $T_1$ rooted at $v_1$.
        Relabel its vertices to obtain $\relab{u}{v_1}((T_1)_{\rm max})$,
        which is a tree rooted at $w$,
        since $w$ is the second-smallest vertex in $(T_1)_{\rm max}$.
  \item[Step 4.] Create a new tree $T'_1$, rooted at $v_1$, as follows:
\begin{itemize}
   \item Attach the subtrees in $L_1 \cup D_{L_1}$ to vertex~$v_1$
       just as they are in $T_1$.
   \item Attach the tree $\relab{u}{v_1}((T_1)_{\rm max})$
        onto $v_1$ via a new edge from $v_1$ to $w$.
\end{itemize}
  \item[Step 5.] Create a new tree $\phi'(T_1)$, rooted at $v_1$,
      by grafting the remaining trees in
      ${T_1 \setminus (L_1 \cup D_{L_1} \cup (T_1)_{\rm max})}$
      onto $T'_1$ by identifying vertices with the same label.
      Note that, in $\phi'(T_1)$, the root $v_1$ has only one
      higher-numbered child, namely $w$.
      See Figure~\ref{fig-bijec1&2-case3b-phi'}.
  \item[Step 6.] Graft $\phi(T_0)$ onto $\phi'(T_1)$
      by identifying vertex $u$ to obtain $T'$.
  \item[Step 7.] Attach the trees rooted at $v_2,\ldots,v_k$
      onto the root $v_1$ of $T'$ via new edges, to obtain $\sigma(T)$.
      See Figure~\ref{fig-bijec1&2-case3b}.
\end{itemize}
Note that $\phi'(T_1)$ has one more vertex than $T_1$;
but one vertex is lost in Step~6
in identifying the vertices $u$ of $\phi(T_0)$ and $\phi'(T_1)$,
so $\sigma(T)$ has the same number of vertices as $T$.

We obtain thereby a tree $\sigma(T)$ rooted at $v_1$
with $k$ higher-numbered children $w,v_2,\ldots,v_k$,
which also has the following properties:
\begin{itemize}
  \item[A2)] vertex~$1$ has at least one child (namely, $r$)
      that is lower-numbered than the root $v_1$;
  \item[B1)] the top ancestor of vertex~$1$ (namely, $w$)
      is higher-numbered than the root $v_1$.
\end{itemize}

The $k$ edges in $T$ that emanate from vertex~$1$ to its children
(shown in blue in Figure~\ref{fig-bijec1&2-case3b}) are clearly proper.
These edges are mapped to the $k$ edges in $\sigma(T)$ that emanate from
the root $v_1$ to its higher-numbered children,
which are regular by rule~(I3) of Definition~\ref{def-irregular}.
By the discussion in Lemma~\ref{lem-bijec1&2-k=0},
the proper/improper edges in $T_0$
are all mapped to regular/irregular edges in $\phi(T_0)$.
By similar reasoning,
the proper/improper edges in $T_1$
are all mapped to regular/irregular edges in $\phi'(T_1)$.
Finally, the proper/improper edges in $T_2,\ldots,T_k \subseteq T$
are mapped to regular/irregular edges in $T_2,\ldots,T_k \subseteq \sigma(T)$
by rules~(I1), (I2) and (I6).
We have therefore shown that the bijection $\psi_T$
maps proper/improper edges in $T$ 
onto regular/irregular edges in $\sigma(T)$.

\bigskip

We complete the proof of Proposition~\ref{prop.bijection.model21}
by remarking that the map $\sigma$ can be reversed,
since every tree in Model~1 must satisfy one of the four properties
consisting of a pair (A1)/(A2) and (B1)/(B2)
stated in Case~I, Case~II, Case~III(a) and Case~III(b).
\qed


\clearpage

\begin{thebibliography}{99}

\bibitem{Aigner_99}  M. Aigner, Catalan-like numbers and determinants,
   J. Combin. Theory A {\bf 87}, 33--51 (1999).



\bibitem{Asner_70}  B.A. Asner, Jr., On the total nonnegativity of the
   Hurwitz matrix, SIAM J. Appl. Math. {\bf 18}, 407--414 (1970).


\bibitem{Barry_16}  P. Barry, {\em Riordan Arrays: A Primer}\/
   (Logic Press, County Kildare, Ireland, 2016).

\bibitem{Barry_17}  P. Barry, Constructing exponential Riordan arrays
   from their $A$ and $Z$ sequences,
   J. Integer. Seq. {\bf 17}, article 14.2.6 (2014). 

\bibitem{Bergeron_98}  F. Bergeron, G. Labelle and P. Leroux,
      {\em Combinatorial Species and Tree-Like Structures}\/
      (Cambridge University Press, Cambridge--New York, 1998).

\bibitem{Borodin_17}  A. Borodin and G. Olshanski,
   {\em Representations of the Infinite Symmetric Group}\/
   (Cambridge University Press, Cambridge, 2017).

\bibitem{Bouwkamp_86}  C.J. Bouwkamp,
   Solution to Problem 85-16: A conjectured definite integral,
   SIAM Rev. {\bf 28}, 568--569 (1986).

\bibitem{Brenti_89}  F. Brenti,
   Unimodal, log-concave and P\'olya frequency sequences in combinatorics,
   Mem. Amer. Math. Soc. {\bf 81}, no.~413 (1989).

\bibitem{Brenti_95}  F. Brenti, Combinatorics and total positivity,
   J. Combin. Theory A {\bf 71}, 175--218 (1995).

\bibitem{Brenti_96}  F. Brenti, The applications of total positivity
   to combinatorics, and conversely.
   In:  {\em Total Positivity and its Applications}\/,
   edited by M.~Gasca and C.A.~Micchelli
   (Kluwer, Dordrecht, 1996), pp.~451--473.

\bibitem{Brumfiel_79}  G.W. Brumfiel,
    {\em Partially Ordered Rings and Semi-Algebraic Geometry}\/,
    London Mathematical Society Lecture Note Series \#37
    (Cambridge University Press, Cambridge--New York, 1979).


\bibitem{Chauve_99}  C. Chauve, S. Dulucq and O. Guibert,
   Enumeration of some labelled trees,
   Research Report RR-1226-99, LaBRI, Universit\'e Bordeaux~I (1999).
   Available on-line at
   \url{http://www.cecm.sfu.ca/~cchauve/Publications/RR-1226-99.ps}

\bibitem{Chauve_00}  C. Chauve, S. Dulucq and O. Guibert,
   Enumeration of some labelled trees,
   in {\em Formal Power Series and Algebraic Combinatorics}\/
   (FPSAC'00, Moscow, June 2000),
   edited by D.~Krob, A.A.~Mikhalev and A.V.~Mikhalev
   (Springer-Verlag, Berlin, 2000), pp.~146--157.

\bibitem{Chen_93}   W.Y.C. Chen, Context-free grammars, differential operators
   and formal power series,
   Theoret. Comput. Sci. {\bf 117}, 113--129 (1993).

\bibitem{Chen_21}  W.Y.C. Chen and H.R.L. Yang,
   A context-free grammar for the Ramanujan--Shor polynomials,
   Adv. Appl. Math. {\bf 126}, 101908 (2021), 24 pp.

\bibitem{Chen_15a}  X. Chen, H. Liang and Y. Wang,
   Total positivity of Riordan arrays,
   European J. Combin. {\bf 46}, 68--74 (2015).

\bibitem{Chen_15b}  X. Chen, H. Liang and Y. Wang,
   Total positivity of recursive matrices,
   Lin. Alg. Appl. {\bf 471}, 383--393 (2015).

\bibitem{Chen_19}  X. Chen and Y. Wang,
   Notes on the total positivity of Riordan arrays,
   Lin. Alg. Appl. {\bf 569}, 156--161 (2019).



\bibitem{Corless_96}  R.M. Corless, G.H. Gonnet, D.E.G. Hare,
   D.J. Jeffrey and D.E. Knuth,
   On the Lambert $W$ function,
   Adv. Comput. Math. {\bf 5}, 329--359 (1996).

\bibitem{Critzer_12_OEIS}  G. Critzer, 8~January 2012,
   contribution to \cite[A071207]{OEIS}.

\bibitem{Curtis_98}  E.B. Curtis, D. Ingerman and J.A. Morrow,
   Circular planar graphs and resistor networks,
   Lin. Alg. Appl. {\bf 283}, 115--150 (1998).

\bibitem{Deutsch_04} E. Deutsch and L. Shapiro, Exponential Riordan arrays,
   handwritten lecture notes, Nankai University, 26~February 2004,
   available on-line at
   \url{http://www.combinatorics.net/ppt2004/Louis%20W.%20Shapiro/shapiro.pdf}

\bibitem{Deutsch_05}  E. Deutsch, L. Ferrari and S. Rinaldi,
   Production matrices,  Adv. Appl. Math. {\bf 34}, 101--122 (2005).

\bibitem{Deutsch_09}  E. Deutsch, L. Ferrari and S. Rinaldi,
   Production matrices and Riordan arrays,
   Ann. Comb. {\bf 13}, 65--85 (2009).

\bibitem{Ding_22}  M.-J. Ding, L. Mu and B.-X. Zhu,
   Coefficientwise total positivity and $\gamma$-positivity of the
   generalized Eulerian polynomials,
   preprint (2022).

\bibitem{Dumont_96}  D. Dumont and A. Ramamonjisoa,
   Grammaire de Ramanujan et arbres de Cayley,
   Electron. J. Combin. {\bf 3}, no.~2, \#R17 (1996).

\bibitem{Dyachenko_14}  A. Dyachenko, Total nonnegativity of infinite
   Hurwitz matrices of entire and meromorphic functions,
   Complex Anal. Oper. Theory {\bf 8}, 1097--1127 (2014).

\bibitem{Fallat_11}  S.M. Fallat and C.R. Johnson,
   {\em Totally Nonnegative Matrices}\/
   (Princeton University Press, Princeton NJ, 2011).

\bibitem{Fomin_01}  S. Fomin, Loop-erased walks and total positivity,
   Trans. Amer. Math. Soc. {\bf 353}, 3563--3583 (2001).

\bibitem{Fomin_10}  S. Fomin, Total positivity and cluster algebras,
   in {\em Proceedings of the International Congress of Mathematicians}\/,
   vol.~II,
   edited by R.~Bhatia, A.~Pal, G.~Rangarajan, V.~Srinivas and M.~Vanninathan
   (Hindustan Book Agency, New Delhi, 2010), pp.~125--145.

\bibitem{Fomin_forthcoming}  S. Fomin, L. Williams and A. Zelevinsky,
   {\em Introduction to Cluster Algebras}\/,
   forthcoming book;
   preliminary draft of Chapters~1--7 posted at
   arXiv:1608.05735 [math.CO], arXiv:1707.07190 [math.CO],
   arXiv:2008.09189 [math.AC] and arXiv:2106.02160 [math.CO]
   at arXiv.org.

\bibitem{Fomin_99}  S. Fomin and A. Zelevinsky,
   Double Bruhat cells and total positivity,
   J. Amer. Math. Soc. {\bf 12}, 335--380 (1999).

\bibitem{Fomin_00}  S. Fomin and A. Zelevinsky,
   Total positivity: tests and parametrizations,
   Math. Intelligencer {\bf 22}, no.~1, 23--33 (2000).


\bibitem{Gantmacher_02}  F.R. Gantmacher and M.G. Krein,
   {\em Oscillation Matrices and Kernels and Small Vibrations of
       Mechanical Systems}\/
   (AMS Chelsea Publishing, Providence RI, 2002).
   Based on the second Russian edition, 1950.

\bibitem{Gantmakher_37}  F. Gantmakher and M. Krein,
   Sur les matrices compl\`etement non n\'egatives et oscillatoires,
   Compositio Math. {\bf 4}, 445--476 (1937).

\bibitem{Gao_non-triangular_transforms}
   A.L.L. Gao, M. P\'etr\'eolle, A.D. Sokal, A.L.B. Yang and B.-X. Zhu,
   Total positivity of a class of Riordan-like matrices,
      implying a class of non-triangular linear transforms
      that preserve Hankel-total positivity,
   in preparation.

\bibitem{Gasca_96}  M. Gasca and C.A. Micchelli, eds.,
   {\em Total Positivity and its Applications}\/
   (Kluwer, Dordrecht, 1996).

\bibitem{Gessel_16}  I.M. Gessel, Lagrange inversion,
   J. Combin. Theory A {\bf 144}, 212--249 (2016).

\bibitem{Gilmore_21}  T. Gilmore,
   Trees, forests, and total positivity: I.~$q$-trees and $q$-forests matrices,
   Electron. J. Combin. {\bf 28}, no.~3, \#P3.54 (2021).

\bibitem{Guo_18}  V.J.W. Guo, A bijective proof of the Shor recurrence,
   European J. Combin. {\bf 70}, 92--98 (2018).


\bibitem{Guo_07}  V.J.W. Guo and J. Zeng,  A generalization of the
   Ramanujan polynomials and plane trees,
   Adv. Appl. Math. {\bf 39}, 96--115 (2007).

\bibitem{He_15}  T.-X. He, Matrix characterizations of Riordan arrays,
   Lin. Alg. Appl. {\bf 465}, 15--42 (2015).

\bibitem{Holtz_03}  O. Holtz, Hermite--Biehler, Routh--Hurwitz,
   and total positivity, Lin. Alg. Appl. {\bf 372}, 105--110 (2003).

\bibitem{Josuat-Verges_15}  M. Josuat-Verg\`es,
   Derivatives of the tree function,
   Ramanujan J. {\bf 38}, 1--15 (2015).

\bibitem{Kalugin_12b}  G.A. Kalugin, D.J. Jeffrey and R.M. Corless,
   Bernstein, Pick, Poisson and related integral expressions for Lambert $W$,
   Integral Transforms Spec. Funct. {\bf 23}, 817--829 (2012).

\bibitem{Karlin_68}  S. Karlin, {\em Total Positivity}\/
   (Stanford University Press, Stanford CA, 1968).

\bibitem{Karlin_59}  S. Karlin and J. McGregor, Coincidence probabilities,
    Pacific J. Math. {\bf 9}, 1141--1164 (1959).

\bibitem{Kemperman_82}  J.H.B. Kemperman,
   A Hurwitz matrix is totally positive,
   SIAM J. Math. Anal. {\bf 13}, 331--341 (1982).

\bibitem{Labelle_81}  G. Labelle, Une nouvelle d\'emonstration combinatoire
   des formules d'inversion de Lagrange,
   Adv. Math. {\bf 42}, 217--247 (1981).

\bibitem{Lam_84}  T.Y. Lam, An introduction to real algebra,
   Rocky Mountain J. Math. {\bf 14}, 767--814 (1984).

\bibitem{Liang_16}  H. Liang, L. Mu and Y. Wang,
   Catalan-like numbers and Stieltjes moment sequences,
   Discrete Math. {\bf 339}, 484--488 (2016).

\bibitem{Lin_14}  Z. Lin and J. Zeng, Positivity properties of Jacobi--Stirling
   numbers and generalized Ramanujan polynomials,
   Adv. Appl. Math. {\bf 53}, 12--27 (2014).

\bibitem{Lis_17}  M. Lis, The planar Ising model and total positivity,
   J. Stat. Phys. {\bf 166}, 72--89 (2017).

\bibitem{Lusztig_94}  G. Lusztig, Total positivity in reductive groups,
   in {\em Lie Theory and Geometry}\/,
   edited by J.-L.~Brylinski, R.~Brylinski, V.~Guillemin and V.~Kac
   (Birkh\"auser Boston, Boston MA, 1994), pp.~531--568.

\bibitem{Lusztig_98}  G. Lusztig, Introduction to total positivity,
   in {\em Positivity in Lie Theory: Open Problems}\/,
   edited by J.~Hilgert, J.D.~Lawson, K.-H.~Neeb and E.B.~Vinberg
   (de Gruyter, Berlin, 1998), pp.~133--145.

\bibitem{Lusztig_08}  G. Lusztig, A survey of total positivity,
   Milan J. Math. {\bf 76}, 125--134 (2008).

\bibitem{Marshall_08}  M. Marshall, {\em Positive Polynomials and
    Sums of Squares}\/, Mathematical Surveys and Monographs \#146
    (American Mathematical Society, Providence RI, 2008).

\bibitem{Merlini_00}  D. Merlini and M.C. Verri,
   Generating trees and proper Riordan arrays,
   Discrete Math. {\bf 218}, 167--183 (2000).

\bibitem{Moon_70} J.W. Moon, {\em Counting Labelled Trees}\/
   (Canadian Mathematical Congress, Montreal, 1970).


\bibitem{Mu_20}  L. Mu and Y. Wang, private communication (2020).


\bibitem{OEIS}  The On-Line Encyclopedia of Integer Sequences,
   published electronically at \url{http://oeis.org}

\bibitem{latpath_lah}  M. P\'etr\'eolle and A.D. Sokal,
   Lattice paths and branched continued fractions, II:
   Multivariate Lah polynomials and Lah symmetric functions,
   European J. Combin. {\bf 92}, 103235 (2021).

\bibitem{latpath_SRTR}  M. P\'etr\'eolle, A.D. Sokal and B.-X. Zhu,
   Lattice paths and branched continued fractions:
       An infinite sequence of generalizations
       of the Stieltjes--Rogers and Thron--Rogers polynomials,
       with coefficientwise Hankel-total positivity,
   preprint (2018), arXiv:1807.03271 [math.CO] at arXiv.org,
   to appear in Memoirs Amer. Math. Soc.

\bibitem{Pinkus_10}  A. Pinkus, {\em Totally Positive Matrices}\/
   (Cambridge University Press, Cambridge, 2010).


\bibitem{Prestel_01}  A. Prestel and C.N. Delzell,
   {\em Positive Polynomials: From Hilbert's 17th Problem to Real Algebra}\/
   (Springer-Verlag, Berlin, 2001).

\bibitem{Randazzo_21}  L. Randazzo, Arboretum for a generalisation of
   Ramanujan polynomials, Ramanujan J. {\bf 54}, 591--604 (2021).

\bibitem{Riordan_68}  J. Riordan, {\em Combinatorial Identities}\/
   (Wiley, New York, 1968).
   [Reprinted with corrections by Robert E.~Krieger Publishing Co.,
    Huntington NY, 1979.]





\bibitem{Schoenberg_53}  I.J. Schoenberg and A. Whitney,
   On P\'olya frequency functions. III. The positivity of translation
   determinants with an application to the interpolation problem
   by spline curves,
   Trans. Amer. Math. Soc. {\bf 74}, 246--259 (1953).

\bibitem{Shapiro_91}  L.W. Shapiro, S. Getu, W.J. Woan and L.C. Woodson,
   The Riordan group,
   Discrete Appl. Math. {\bf 34}, 229--239 (1991).

\bibitem{Shapiro_22}  L. Shapiro, R. Sprugnoli, P. Barry, G.-S. Cheon,
   T.-X. He, D. Merlini and W. Wang,
   {\em The Riordan Group and Applications}\/
   (Springer, Cham, 2022).


\bibitem{Shor_95}  P.W. Shor, A new proof of Cayley's formula for counting
   labeled trees, J. Combin. Theory A {\bf 71}, 154--158 (1995).


\bibitem{Skandera_03}  M. Skandera, Introductory notes on total positivity
   (June 2003), available at
   \url{http://www.math.lsa.umich.edu/~fomin/565/intp.ps}

\bibitem{Sokal_flajolet}  A.D. Sokal, Coefficientwise total positivity
   (via continued fractions) for some Hankel matrices of combinatorial
   polynomials, talk at the S\'eminaire de Combinatoire Philippe Flajolet,
   Institut Henri Poincar\'e, Paris, 5 June 2014;
   transparencies available at
   \url{http://semflajolet.math.cnrs.fr/index.php/Main/2013-2014}

\bibitem{Sokal_OPSFA}  A.D. Sokal, Coefficientwise Hankel-total positivity,
   talk at the 15th International Symposium on Orthogonal Polynomials,
   Special Functions and Applications (OPSFA 2019),
   Hagenberg, Austria, 23 July 2019;
   transparencies available at
   \url{https://www3.risc.jku.at/conferences/opsfa2019/talk/sokal.pdf}

\bibitem{Sokal_trees_enumeration}  A.D. Sokal,
   A remark on the enumeration of rooted labeled trees,
   Discrete Math. {\bf 343}, 111865 (2020).

\bibitem{forests_totalpos}  A.D. Sokal,
   Total positivity of some polynomial matrices
   that enumerate labeled trees and forests,
   I: Forests of rooted labeled trees,
   Monatsh. Math. {\bf 200}, 389--452 (2023).

\bibitem{Sokal_totalpos}  A.D. Sokal, Coefficientwise total positivity
   (via continued fractions) for some Hankel matrices of combinatorial
   polynomials, in preparation.

\bibitem{Sprugnoli_94}  R. Sprugnoli, Riordan arrays and combinatorial sums,
   Discrete Math. {\bf 132}, 267--290 (1994).

\bibitem{Stanley_99}  R.P. Stanley, {\em Enumerative Combinatorics}\/,
      vol.~2 (Cambridge University Press, Cambridge--New York, 1999).

\bibitem{Stembridge_91}  J.R. Stembridge,
   Immanants of totally positive matrices are nonnegative,
   Bull. London Math. Soc. {\bf 23}, 422--428 (1991).

\bibitem{Stieltjes_1889}  T.J. Stieltjes, Sur la r\'eduction en fraction
   continue d'une s\'erie proc\'edant selon les puissances descendantes
   d'une variable,
   Ann. Fac. Sci. Toulouse {\bf 3}, H1--H17 (1889).

\bibitem{Stieltjes_1894}  T.J. Stieltjes, Recherches sur les fractions
    continues, Ann. Fac. Sci. Toulouse {\bf 8}, J1--J122 (1894)
    and {\bf 9}, A1--A47 (1895).
    [Reprinted, together with an English translation,
     in T.J. Stieltjes, {\em \OE{}uvres Compl\`etes/Collected Papers}\/
     (Springer-Verlag, Berlin, 1993), vol.~II, pp.~401--566 and 609--745.]

\bibitem{Thoma_64}  E. Thoma, Die unzerlegbaren, positiv-definiten
   Klassenfunktionen der abz\"ahlbar unendlichen, symmetrischen Gruppe,
   Math. Z. {\bf 85}, 40--61 (1964).


\bibitem{Zeng_99}  J. Zeng, A Ramanujan sequence that refines the Cayley
   formula for trees, Ramanujan J. {\bf 3}, 45--54 (1999).

\bibitem{Zhu_13}  B.-X. Zhu, Log-convexity and strong $q$-log-convexity
   for some triangular arrays,
   Adv. Appl. Math. {\bf 50}, 595--606 (2013).

\bibitem{Zhu_14}  B.-X. Zhu, Some positivities in certain triangular arrays,
    Proc. Amer. Math. Soc. {\bf 142}, 2943--2952 (2014).


\bibitem{Zhu_21a}  B.-X. Zhu, Total positivity from the exponential
   Riordan arrays, SIAM J. Discrete Math. {\bf 35}, 2971--3003 (2021).

\bibitem{Zhu_21b}  B.-X. Zhu, Stieltjes moment properties and
   continued fractions from combinatorial triangles,
   Adv. Appl. Math. {\bf 130}, 102232 (2021).

\bibitem{Zhu_22}  B.-X. Zhu, Coefficientwise Hankel-total positivity of
   row-generating polynomials for the $m$-Jacobi-Rogers triangle,
   preprint (March~2022), arXiv:2202.03793 [math.CO].

\end{thebibliography}
\end{document}